\documentclass{jssc}

\def\textsubscript#1%
{$_{\text{#1}}$}
\newcommand*\supercite[1]{\textsuperscript{\cite{#1}}}
\def\cdd{\mbox{\boldmath$\cdot$}~}

\newcommand{\rulex}{\hfill\rule{1mm}{3mm}}

\input{amssym.def}
\usepackage{graphicx}
\usepackage{amssymb}
\usepackage{multirow}
\usepackage{multicol}
\usepackage{bm}
\usepackage{color}

\makeatletter
\def\@oddfoot{\hfill}
\newcount\shumeicount
\def\setshumei#1#2#3{%
  \shumeicount=\count0
  \def\@oddhead{%
    \raise-5pt\hbox to0pt{\vrule width\hsize height 0pt depth 0.4pt\hss}\relax
    \ifnum \shumeicount=\count0
      \raise-7pt\hbox to0pt{\vrule width\hsize height 0pt depth 0.4pt\hss}\relax
      #1
    \else
      \ifodd\count0
        #2
      \else
        #3
       \fi
     \fi
  }%
}
\makeatother
\makeatletter
\def\@oddfoot{\hfill}
\newcount\shujiaocount
\def\setshujiao{%
  \shujiaocount=\count0
  \def\@oddfoot{%
      \ifodd\count0
      \else
      \fi
  }%
}
\makeatother
\def\title#1#2#3#4{{
  \vspace*{0.3cm}
  \begin{flushleft} \Large\bf #1\end{flushleft}
  \vspace*{-0.2cm}
      \begin{flushleft}
      \bf #2
      \end{flushleft}
      \footnotetext{\hspace{-6mm} #3\\ #4}}}

\def\dshm#1#2#3#4
{\setshumei{J Syst Sci Complex (#1) #2:
{\thepage--\pageref{LastPage}}\hfill}
            {\hfill {\small #3}\hfill\hbox to0pt{\hss\thepage}}
            {\hbox to0pt{\thepage\hss}\hfill {\small #4}\hfill
            }
            \setshujiao}
\def\drd#1#2
{{\vskip 1cm\small \begin{flushleft}
 #1 \\
 #2\\
\copyright The Editorial Office of JSSC \&  Springer-Verlag GmbH
Germany 2021
\end{flushleft}}}

\def\hat{\widehat}

\def\bar{\overline}
\def\epsilon{\varepsilon}
\def\proof{\vspace{1mm}\indent {\it Proof}\quad}

\begin{document}

\title{NSNO: Neumann Series Neural Operator for Solving Helmholtz Equations in Inhomogeneous Medium $^*$}
{\uppercase{CHEN} Fukai \cdd \uppercase{Liu} Ziyang \cdd \uppercase{LIN}
Guochang \cdd \uppercase{Chen} Junqing \cdd \uppercase{Shi} Zuoqiang}
{\uppercase{Chen} Fukai \\
Department of Mathematical Sciences, Tsinghua University, Beijing, $100084$, China;  \\
Email: cfk19@mails.tsinghua.edu.cn.  \\  
\uppercase{Liu} Ziyang\\
Yau Mathematical Sciences Center, Tsinghua University, Beijing, $100084$, China; \\
Email: zy-liu20@mails.tsinghua.edu.cn.  \\ 
\uppercase{Lin} Guochang \\
Yau Mathematical Sciences Center, Tsinghua University, Beijing, $100084$, China; \\
Email: lingc19@mails.tsinghua.edu.cn. \\
\uppercase{Chen} Junqing\\
Department of Mathematical Sciences, Tsinghua University, Beijing, $100084$, China;  \\ 
Email: jqchen@tsinghua.edu.cn. \\
\uppercase{Shi} Zuoqiang (Corresponding author)\\
Yau Mathematical Sciences Center, Tsinghua University, Beijing, $100084$, China;  Yanqi Lake Beijing Institute of Mathematical Sciences and Applications, Beijing 101408, China.  Email: zqshi@tsinghua.edu.cn.
   } 
{$^*$This research was supported by NSFC under Grant No. 92370125 and National Key  R\&D Program of China 2019YFA0709600, 2019YFA0709602.}


\drd{DOI: }{Received: x x 20xx}{ / Revised: x x 20xx}


\dshm{20XX}{XX}{\uppercase{Neumann Series Neural Operator}}{\uppercase{Chen Fukai} $\cdd$ \uppercase{Liu Ziyang} $\cdd$ \uppercase{Lin Guochang} $\cdd$ \uppercase{Chen Junqing} $\cdd$
\uppercase{Shi Zuoqiang}}

\Abstract{
In this paper, we propose Neumann Series Neural Operator (NSNO) to learn the solution operator of Helmholtz equation from inhomogeneity coefficients and source terms to solutions. Helmholtz equation is a crucial partial differential equation (PDE) with applications in various scientific and engineering fields. However, efficient solver of Helmholtz equation is still a big challenge especially in the case of high wavenumber. Recently, deep learning has shown great potential in solving PDEs especially in learning solution operators. Inspired by Neumann series in Helmholtz equation, we design a novel network architecture in which U-Net is embedded inside to capture the multiscale feature. Extensive experiments show that the proposed NSNO significantly outperforms the state-of-the-art FNO with at least 60\% lower relative $L^2$-error, especially in the large wavenumber case, and has 50\% lower computational cost and less data requirement. Moreover, NSNO can be used as the surrogate model in inverse scattering problems. Numerical tests show that NSNO is able to give comparable results with traditional finite difference forward solver while the computational cost is reduced tremendously.}      




\section{Introduction}

Helmholtz equation is a fundamental partial differential equation (PDE) describing wave propagation in many areas of physics and engineering such as acoustics, electromagnetics \supercite{colton1998inverse}, and medical imaging \supercite{arridge1999optical}. One of the most common scenarios is the propagation of waves in medium with spatially varying properties, which frequently occurs in the scattering of acoustic and electromagnetic waves\supercite{colton2000recent}, such that the governing equation is given by
\begin{equation}\label{eq1}
    \Delta u+k^2(1+q(x))u=f(x),
\end{equation}
where $u$ is the scalar field representing the wave, $f(x)$ is the source term, $q(x)$ is the coefficient representing the inhomogeneity, which is usually compactly supported.

Traditional numerical methods for solving Helmholtz equation in inhomogeneous medium numerically such as finite difference method (FDM)\supercite{singer1998high} and finite element method (FEM)\supercite{babuvska1995generalized} require a fine grid or mesh to capture the high-frequency components of the solution when the wavenumber is large, which leads to
a large scale indefinite linear system. 
Modern methods such as Krylov subspace methods\supercite{saad1986gmres} for solving large indefinte linear systems are known to have a slow convergence rate for Helmholtz problem, especially in inhomogeneous medium case\supercite{ernst2011difficult}. Moreover, in inverse scattering problems, the Helmholtz equation \eqref{eq1} has to be solved repeatedly with different coefficients $q$ and source terms $f$, which will be computationally intolerable if these equations are solved separately. Therefore, solving the inhomogeneous Helmholtz equation still remains a challenging task and it is of significant importance to seek for the solution operator mapping the coefficient $q$ and source term $f$ simultaneously to the solution to the Helmholtz equation \eqref{eq1}.

Benefiting from the attractive capability of neural networks in approximating functions\supercite{hornik1989multilayer}, deep learning has made remarkable progress in areas such as image recognition and natural language processing\supercite{goodfellow2016deep, young2018recent}. Therefore, deep learning is recognized as a promising solution to solve partial differential equations and the solution operator from the parameter space to the solution space\supercite{beck2020overview, kovachki2021neural}. By using neural networks with numerous learnable parameters as an ansatz to represent the solution or the solution operator, and training neural networks based on properly designed loss functions, the neural network is able to approximate the required solution or the solution operator.

In the literature, both solving a single PDE and PDE solution operators have attracted wide attention. In \cite{raissi2019physics}, the physics-informed neural network (PINN) taking the residual of the PDEs as the loss function to train the neural network is presented. Based on PINN, a neural network based solver for Helmholtz equations in homogeneous background has been proposed in \cite{alkhalifah2021wavefield}, while plane wave activation functions are utilized in \cite{cui2022efficient} to further improve accuracy. For PDE solution operators, two general neural operator framework named DeepONet and Fourier neural operator (FNO) have been developed in \cite{lu2021learning, li2020fourier, lu2022comprehensive}. In DeepONet, the features in spatial coordinates and differential equation parameters are extracted separately by the trunk net and the branch net, and are then combined by dot product, while in FNO, the mapping from parameter space to solution space is formulated as an iterative integral, where the integral kernel is parameterized in Fourier space. Furthermore, physics-informed loss is incorporate with DeepONet and FNO in \cite{wang2021learning} and \cite{li2021physics}, respectively, to explore the possibility of unsupervised solution operator learning. An solution operator to the heterogeneous Helmholtz equation mapping the sound speed distribution to the acoustic wavefield is approximated by neural networks based on an iterative scheme in \cite{stanziola2021helmholtz}.

However, two important issues are not fully considered in existing neural network based methods. For one thing, most existing solution operators either map the coefficient in the differential operator or the source term only to the solution. Since the coefficient $q$ and the source term $f$ in the Helmholtz equation \eqref{eq1} belongs to different function spaces, it is a non-trival task to learn the solution operator mapping $q$ and $f$ simultaneously to $u$. For the other, fitting the high oscillations in the solution to Helmholtz equation with large wavenumber always leads to unstable or even divergent training process\supercite{BINet}. Therefore, new network architecture should be designed to capture the multi-scale features in the solution to Helmholtz equations.

To address the two issues above, in the paper, we propose Neumann series neural operator (NSNO). Specifically, the solution to Helmholtz equation is rewritten in the form of a Neumann series. Each term in the Neumann series is the solution to a Helmholtz equation in homogeneous medium subject to different source terms. For the Helmholtz equation corresponding to the first term, the source term is exactly $f$, while $q$ only appears in the source terms of the Helmholtz equations corresponding to the remaining terms, such that $q$ and $f$ are fully decoupled. Moreover, based on the Neumann series, we only need to solve the operator from the source term to the solution of Helmholtz equation in homogeneous medium. Instead of using the Fourier neural operator (FNO) directly to approximate the operator, we propose a novel network architecture combining FNO and U-Net\supercite{UNet} named UNO to capture the multi-scale property of the solution. By extracting and fusing information from different spatial resolutions, the proposed UNO is able to approximate the solution accurately even in cases of large wavenumbers. Our main contributions can be summarized as follows:  

\begin{itemize}
    \item We propose NSNO, a novel framework for neural operators mapping inhomogeneity coefficients and source terms simultaneously to the solution to Helmholtz equation in inhomogeneous medium based on Neumann series. To capture the multi-scale properties of Helmholtz equation, a novel network architecture UNO combining FNO and U-Net is also proposed as the building blocks of NSNO.
    \item Extensive numerical experiments show that NSNO significantly outperforms the state-of-the-art FNO with at least 60\% lower $L^2$-error. The results also show that compared with FNO, the proposed UNO architecture is more efficient (50\% less in computational cost), more accurate on multi-scale problems and has less training data requirement.
    \item We use NSNO as the surrogate model for the forward operator in inverse scattering problem where the scatterer is sampled from the MNIST dataset. NSNO is able to give comparable results with the traditional finite difference forward solver with over 20 times faster in speed.
\end{itemize}

The rest of this paper is organized as follows. The operator learning problem is setup in Section \ref{setup}. In Section \ref{ns}, we show the Neumann series reformulation and analyze the convergence of the Neumann series. The detailed network architecture and training process is presented in Section \ref{net}. Numerical experiments are shown in Section \ref{exp}, and the application of NSNO in inverse scattering problem is given in Section \ref{invp}. At last, conclusion remarks are made in Section \ref{con}.

\section{Problem Setup}\label{setup}

Consider the 2-dimensional Helmholtz equation in inhomogeneous medium subject to the Sommerfeld radiation condition at infinity\supercite{colton1998inverse}
\begin{align}
\Delta u+k^2(1+q(x))u &= f(x), \quad x\in\mathbb{R}^2, \label{total}\\
\lim\limits_{r\rightarrow \infty} \sqrt{r}(\frac{\partial u}{\partial r}-iku)&=0, \quad r=|x|.
\end{align}
where $k$ is the wavenumber, $i=\sqrt{-1}$ is the imaginary unit, $q(x)$ is the compactly supported coefficient representing the inhomogeneity, and $f(x)$ is the source term. In practice, the problem is usually reduced to a bounded domain $\it\Omega$ containing the support of $q$ by introducing an artificial surface. For simplicity, we take $\it\Omega$ as a rectangle and employ the first-order absorbing boundary condition\supercite{jin2015finite}
\begin{equation}
\frac{\partial u}{\partial n}-iku=0, \quad \text{on} \quad \partial \it\Omega.
\end{equation}

Note that in applications such as the inverse scattering problem, we normally needs multiple incident fields to reconstruct the scatterer $q$ better, which can be computational expensive if the forward equations are solved case-by-case. Therefore, in this paper, we aim to learn the following operator:
\begin{equation}\label{operator}
\begin{split}
\mathcal{S}: L^{\infty}(\it\Omega)\times L^2(\it\Omega) &\rightarrow L^2(\it\Omega)\\
(q,f)&\mapsto u,
\end{split}
\end{equation}
where $u$ is the solution to
\begin{align}
\Delta u+k^2(1+q(x))u=f(x),  &\quad \text{in} \quad \it\Omega, \label{op1}\\
\frac{\partial u}{\partial n}-iku=0, &\quad \text{on} \quad  \partial \it\Omega. \label{bd1}
\end{align}

\section{Neumann Series}\label{ns}

Note that for the solution operator $\mathcal{S}$, the input $q$ and $f$ belonging to different function spaces are coupled together, making it difficult to solve the solution operator mapping $q$ and $f$ simultaneously to $u$. Therefore, we consider utilizing the Neumann series to decouple $q$ and $f$, such that the solution operator $\mathcal{S}$ is transformed into another operator $G$ mapping the source term only to the solution. In this section, we first present some preliminary results on the variational formulation and stability estimate of Helmholtz equation, based on which the Neumann series is defined and its convergence property is analyzed.

\subsection{Preliminaries}
Suppose $\it\Omega$ is a rectangle in $\mathbb{R}^2$. Consider the Helmholtz equation
\begin{equation}\label{helm}
\begin{split}
\Delta u+k^2u=g(x),  &\quad \text{in} \quad \it\Omega, \\
\frac{\partial u}{\partial n}-iku=0, &\quad \text{on} \quad  \partial \it\Omega.
\end{split}
\end{equation}
Its variational formulation is defined as
\begin{equation}\label{variation}
    \begin{split}
        &\text{Find} \quad u\in H^1(\it\Omega) \quad \text{such that}\\
        &a(u,v)=(g,v), \quad \forall v\in H^1(\it\Omega),
    \end{split}
\end{equation}
where \begin{equation}
a(u,v)=(\nabla u, \nabla v)-k^2(u, v)-ik\langle u, v\rangle,
\end{equation}
$(\cdot, \cdot)$ and $\langle\cdot, \cdot\rangle$ denote the $L^2$-inner product on $\it\Omega$ and $\partial \it\Omega$, respectively.

The existence of the solution to the variational problem \eqref{variation} and the corresponding stability estimate are presented in the following theorem.
\begin{theorem}\label{thm2.1} The variation problem \eqref{variation} has a unique solution in $H^1(\it\Omega)$. Moreover, if the wavenumber $k>1$, there exists a constant $C$, which only depends on the domain $\it\Omega$, such that for any $f\in L^2(\it\Omega)$, the solution to problem \eqref{variation} satisfies
\begin{equation}
k\Vert u\Vert_{L^2(\it\Omega)}\leq C \Vert g\Vert_{L^2(\it\Omega)}.
\end{equation}
\end{theorem}
The proof of Theorem \ref{thm2.1} can be found in \cite{cms/1188405673}. Based on Theorem \ref{thm2.1}, we can define an operator $G$ by
\begin{equation}
\begin{split}
G: L^2(\it\Omega) &\rightarrow H^1(\it\Omega)\subset L^2(\it\Omega)\\
g & \mapsto u,
\end{split}
\end{equation}
such that $u$ is the solution to the variational problem \eqref{variation}, i.e., the weak solution to the original Helmholtz equation. It can be seen that $G$ is linear and bounded in $L^2(\it\Omega)$ with $\Vert G\Vert_{L^2(\it\Omega)}\leq C/k$.

\subsection{Neumann Series}\label{NS_theory}
Note that the Helmholtz equation \eqref{op1} can be rewritten as $\Delta u+k^2 u=f-k^2qu$, which inspires the following iterative scheme
\begin{equation}\label{iteration}
\Delta u_{n+1}+k^2 u_{n+1}=f-k^2qu_n.
\end{equation}
By the definition and linearity of $G$, the iterative scheme can be rewritten as
\begin{equation}
u_{n+1}=G(f-k^2qu_n)=G(f)+G(-k^2qu_n)=u_0-(k^2Gq)u_n,
\end{equation}
where $u_0=G(f)$. By using the iterative scheme \eqref{iteration} recursively for $N$ steps, we obtain the Neumann series
\begin{equation}\label{Neumann series}
u_N=u_0-(k^2Gq)u_0+(k^2Gq)(k^2Gq)u_0\cdots+(-k^2Gq)^Nu_0.
\end{equation}

The convergence of this Neumann series is presented in the following theorem.

\begin{theorem}[Convergence of Neumann series] \label{thm2.2} When $\Vert q\Vert_{L^{\infty}(\it\Omega)}$ is sufficiently small, the Neumann series \eqref{Neumann series} converges in $L^2(\it\Omega)$ as $N\rightarrow\infty$.
\end{theorem}

\proof It suffices to prove $\Vert -k^2Gq\Vert_{L^2(\it\Omega)}<1$. Note that
\begin{equation}\label{norm}
\Vert -k^2Gq\Vert_{L^2(\it\Omega)}=\sup\limits_{g\in L^2(\it\Omega)}\frac{\Vert -k^2G(qg)\Vert_{L^2(\it\Omega)}}{\Vert g\Vert_{L^2(\it\Omega)}}\leq \sup\limits_{g\in L^2(\it\Omega)}\frac{Ck \Vert qg\Vert_{L^2(\it\Omega)}}{\Vert g\Vert_{L^2(\it\Omega)}}\leq Ck\Vert q\Vert_{L^{\infty}(\it\Omega)}.
\end{equation}
It can be seen that the Neumann series converges as long as $\Vert q\Vert_{L^{\infty}(\it\Omega)}<\frac{1}{Ck}$.
\rulex 

Therefore, based on the Neumann series \eqref{Neumann series}, solving the solution operator $\mathcal{S}$ from $(q,f)$ to $u$ can be transformed into solving $G$, which is the solution operator from the source term only to the solution of Helmholtz equation.

\section{Network Architecture and Training Process}\label{net}

In this section, we first give an overall illustration of the overall network architecture based on the Neumann series \eqref{Neumann series}. Following that, we further elaborate on the approximation of operator $G$ by neural networks. Finally, we discuss the training process, where physical-informed loss is utilized to further improve the performance of the learned operator.

\subsection{Network Architecture}

\begin{figure}
\centering
\includegraphics[width=13cm]{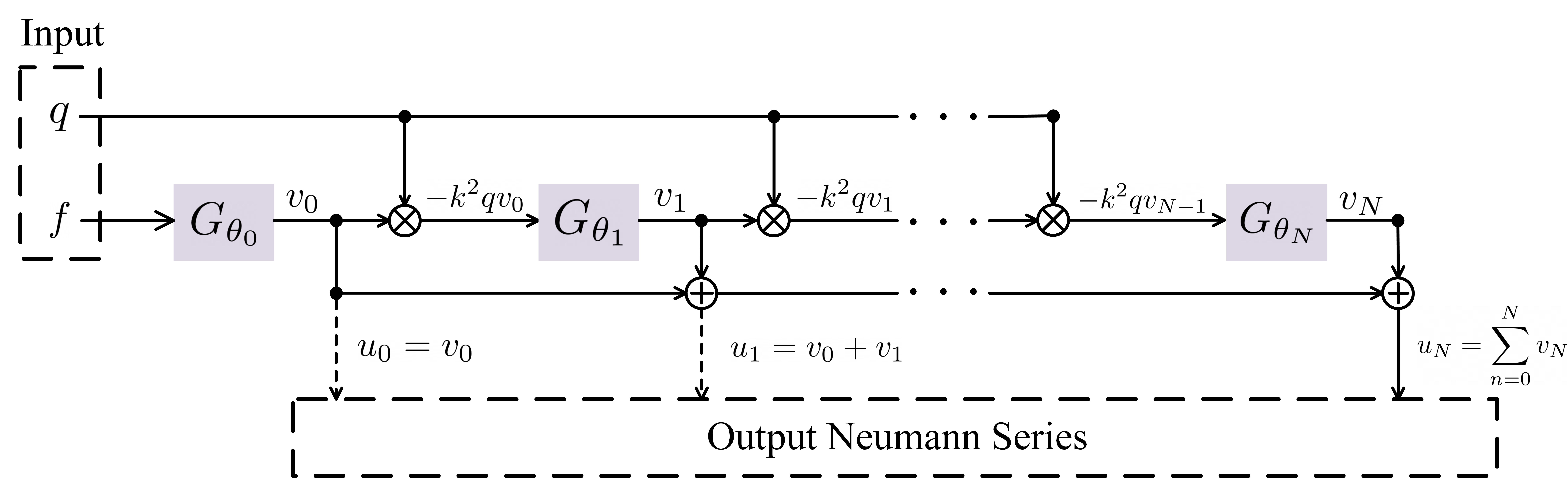}
\caption{Overall Network Architecture of NSNO}
\label{model1}
\end{figure}

Fig. \ref{model1} shows the overall network architecture of the proposed Neumann series neural operator, referred to as NSNO, where the Neumann series is truncated to $N+1$ items\footnote{In the experiments, we have found that $N=2$, i.e., only three items to construct the Neumann series, is sufficient to give an accurate result, as will be discussed in Section \ref{nitem}.}. The operator $G$ is approximated by $N+1$ separate neural networks $G_{\bm{\theta}}=\{G_{\theta_0}, G_{\theta_1}, \cdots, G_{\theta_N}\}$ with learnable parameters $\bm{\theta}=\{\theta_0, \theta_1, \cdots, \theta_N\}$, respectively. Note that although the neural networks $\{G_{\theta_0}, G_{\theta_1}, \cdots, G_{\theta_N}\}$ are all approximation of the same operator $G$, we still use different parameters since compared with the neural network reusing the same parameter multiple times, using different parameters gives the neural network stronger representation ability. The first neural network takes $f$ as the input and outputs the approximation of the first item $v_0\approx u_0$ in Neumann series \eqref{Neumann series}. For the rest $N$ neural networks, $G_{\theta_n}$ takes the multiplication of $-k^2q$ and the output of the last neural network $v_{n-1}$ as the input, and outputs the approximation of the $n+1$-th item $v_n\approx (-k^2Gq)^nu_0$ in Neumann series \eqref{Neumann series}. 

In practice, the domain $\it\Omega\in\mathbb{R}^2$ is discretized by a regularly spaced Cartesian grid of dimension $H\times W$ such that functions defined on $\it\Omega$ can be represented by tensors. The input and output of $G_{\bm{\theta}}$ has the same spatial dimensions as the grid and has two channels in most cases representing the real and imaginary parts. For $G_{\theta_0}$ the input has only one channel if $f$ is real. The detailed structure of $G_{\bm{\theta}}$ will be specified in section \ref{G}.

In summary, NSNO takes the tuple $(q,f)$ as input, where $f$ will only be fed into $G_{\theta_0}$, while $-k^2q$ will be multiplied with the outputs of $G_{\theta_0}, G_{\theta_1}, \cdots, G_{\theta_{N-1}}$. In this way, $q$ and $f$ belonging to different function spaces can be fully decoupled. The output of NSNO is the sum of the outputs of the $N+1$ neural networks $G_{\bm{\theta}}$. Besides, the intermediate results $u_0, u_1, \cdots u_{N-1}$ can also be output if required, adding more flexibility to NSNO. 

\subsection{Network Architecture of $G_{\bm{\theta}}$}\label{G}

In this subsection, we introduce two types of network architecture for $G_{\bm{\theta}}$, i.e., the Fourier neural operator (FNO), which is a direct application of \cite{li2020fourier}, and the U-shaped neural operator, which is a combination of the UNet architecture\supercite{UNet} and the FNO.

\subsubsection{Fourier Neural Operator}

\begin{figure}
\centering
\includegraphics[width=14cm]{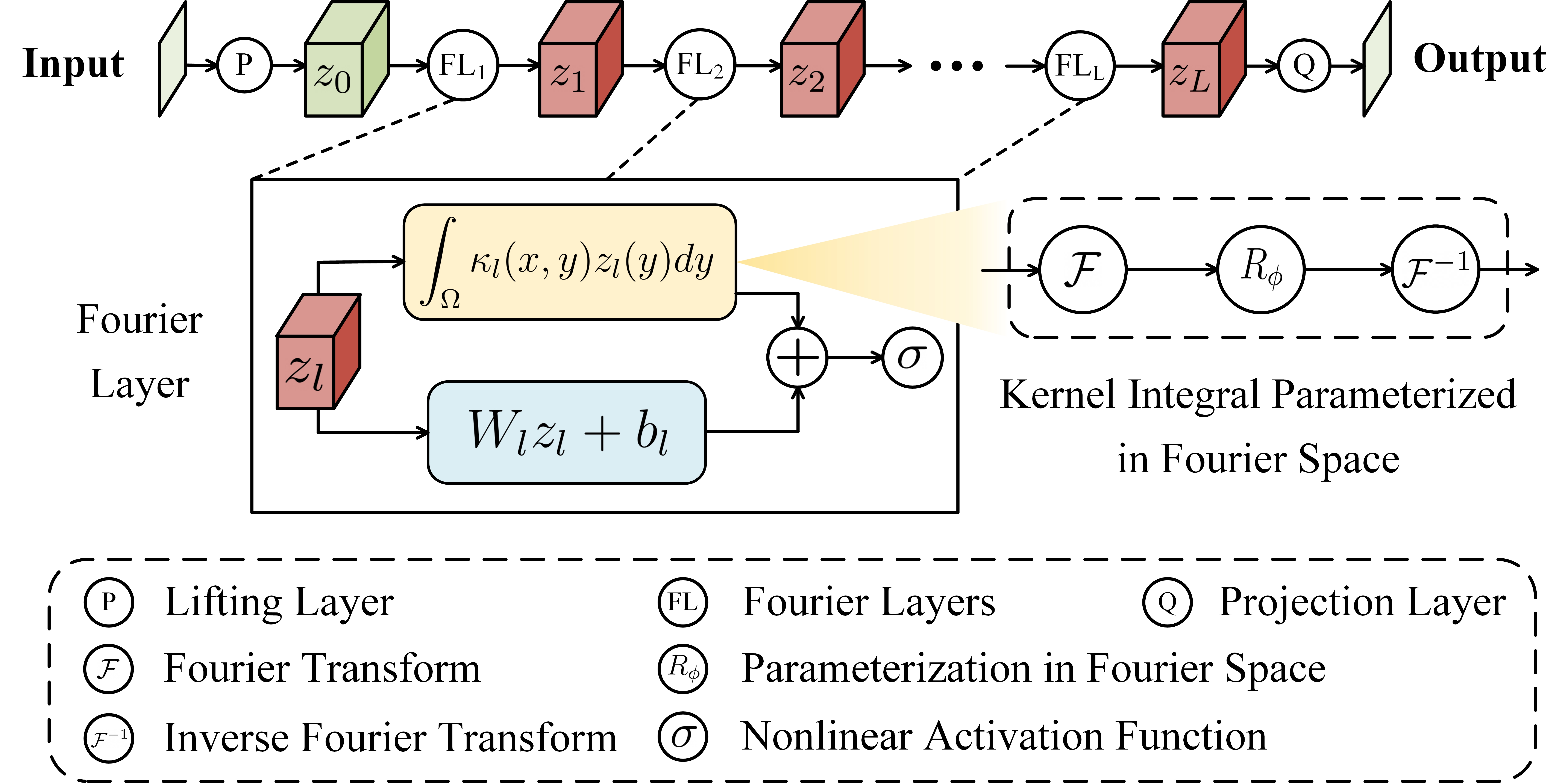}
\caption{Network architecture of FNO}
\label{model4}
\end{figure}

One natural idea is to choose $G_{\bm{\theta}}$ as the Fourier neural operator (FNO) directly. As is shown in Fig. \ref{model4}, the architecture of FNO starts with a lifting layer, followed by a series of Fourier layers, and is ended with a projection layer. 

The lifting layer lifts the input to a higher dimension by a fully connected neural network $P$ such that $z_0$ takes values in $\mathbb{R}^C$. A series of iterative Fourier layers is then applied to $z_0$, resulting in a sequence of functions $z_0\mapsto z_1\mapsto\cdots\mapsto z_L$ taking values in $\mathbb{R}^C$. Specifically, the iterative scheme is given by
\begin{equation}\label{FNO update}
z_{l+1}(x)=\sigma \left(\int_{\it\Omega} \kappa_l(x,y)z_l(y)dy+W_lz_l(x)+b_l\right), \quad l=0, 1, \cdots, L-1,
\end{equation}
where $\kappa_l$ is an integral kernel, $W_l$ is a linear transform, $b_l$ is the bias, which are all learnable, and $\sigma$ is a nonlinear activation function. By letting $\kappa_l(x,y)=\kappa_l(x-y)$ and using the convolution therorem, the kernel integral in \eqref{FNO update} can be rewritten as
\begin{equation}
\int_{\it\Omega} \kappa_l(x,y)z_l(y)dy=\mathcal{F}^{-1}\left(\mathcal{F}(\kappa_l)\cdot\mathcal{F}(z_l)\right)(x):=\mathcal{F}^{-1}\left(R_{\phi}\cdot\mathcal{F}(z_l)\right)(x),
\end{equation}
where $\mathcal{F}$ and $\mathcal{F}^{-1}$ is the Fourier transform and its inverse, respectively, $R_{\phi}:= \mathcal{F}(\kappa_l)$ is the parameterization of $\kappa_l$ in the Fourier space with parameters $\phi$. In the discrete case, the Fourier series is truncated such that higher modes are filtered and only $k_{\text{max}}$ modes are reserved. Therefore, $R_{\phi}$ can be directly parameterized as a complex-valued $k_{\text{max}}\times C\times C$ tensor. Finally, $z_L$ is projected back to the output space with required dimension with another fully connected neural network $Q$.

\subsubsection{U-shaped Neural Operator}

The architecture of FNO enables its superior accuracy in various applications, such as Burger's equation, Darcy flow, and so on. However, new challenges arise in solving Helmholtz equations with multi-scale features, especially in the high wavenumber regime due to the couple between the high frequency wave solutions and the numerical grid\supercite{liu2022ht}. Building on global Fourier transform filtering higher frequency modes, FNO tends to learn an over-smooth solution. In multi-scale problems, FNO fails in capturing the intrinsic multis-cale features in the solution, leading to poor performance.

\begin{figure}
\centering
\includegraphics[width=14cm]{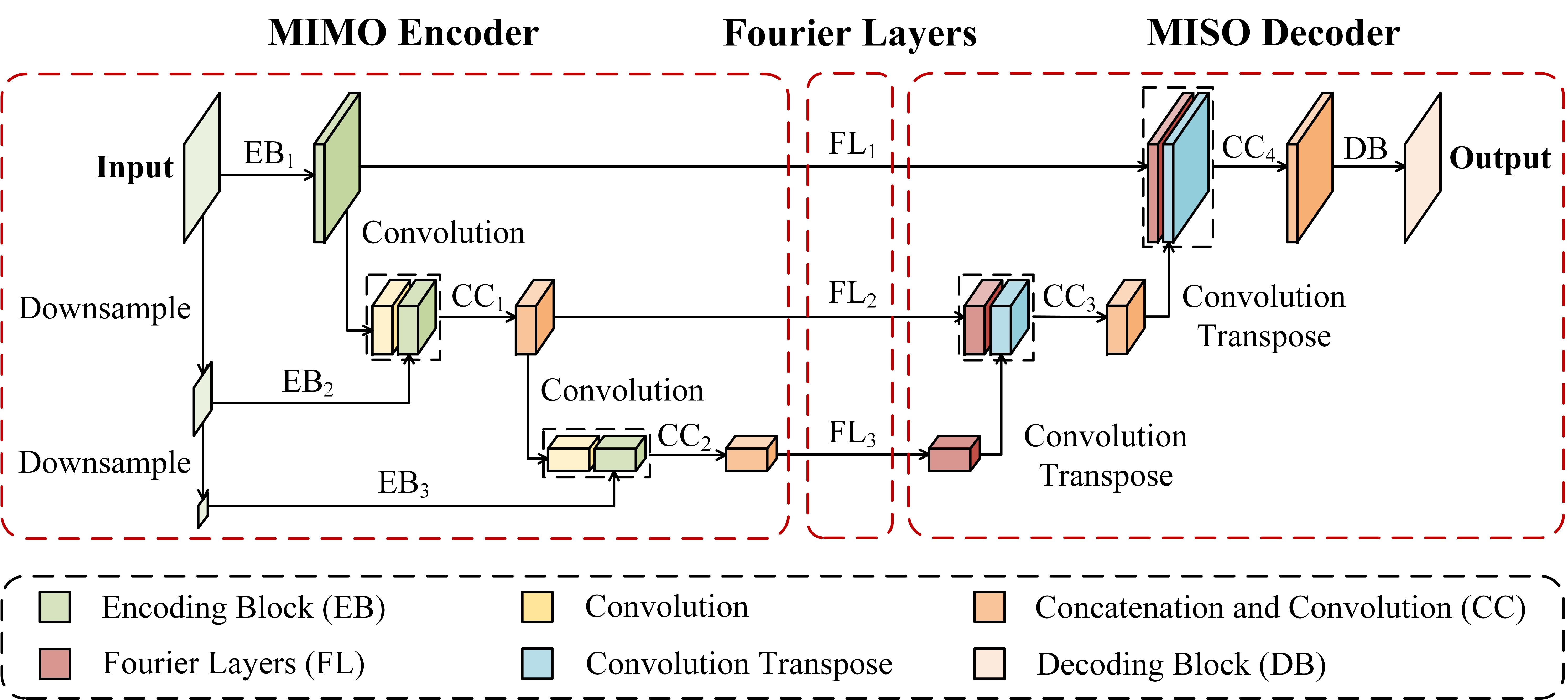}
\caption{Network architecture of UNO}
\label{model2}
\end{figure}

To address this issue, we utilize the U-Net structure widely used in image processing and computer vision\supercite{zhou2018unet++, li2018h}. By using a hierarchical encoder-decoder structure and skip connections, U-Net is able to effectively capture and merge information at different scales, thus having capability to solve multi-scale problems. Therefore, we design a novel network architecture for $G_{\bm{\theta}}$ combining both U-Net and FNO structure termed UNO, where U-Net exploits the multi-scale structure of the solutions, while FNO achieves the transformation from the source term space to the solution space. Specifically, as is shown in Fig. \ref{model2}, UNO consists of three modules, i.e., the multiple-input multiple-output (MIMO) encoder, Fourier layers as skip connections and the multiple-input single-output (MISO) decoder, which will be detailed below.

\subparagraph{Multiple-Input Multiple-Output (MIMO) Encoder}

\begin{figure}
\centering
\includegraphics[width=14cm]{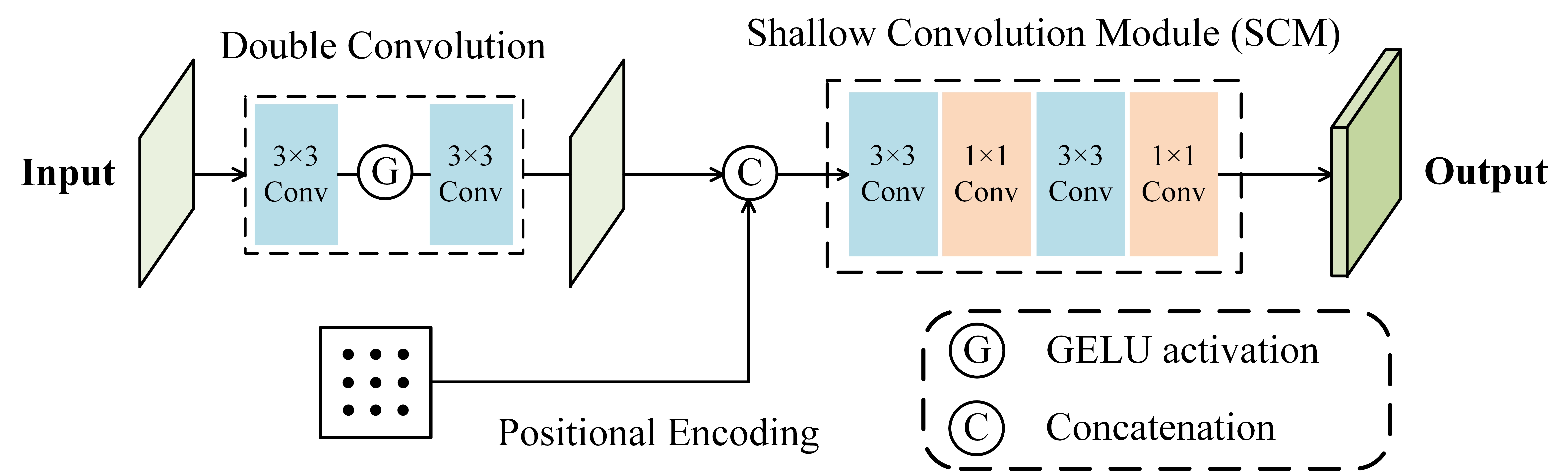}
\caption{Network architecture of the encoding block in UNO}
\label{model3}
\end{figure}

In the MIMO encoder module, we firstly downsample the input twice to spatial dimensions $\frac{H}{2}\times \frac{W}{2}$ and $\frac{H}{4}\times \frac{W}{4}$, respectively. The original input and the two downsampled inputs are fed into three encoding blocks to extract features from different scales. As is shown in Fig. \ref{model3}, in the encoding blocks, the input is firstly fed into a double convolution module consists of two $3\times 3$ convolution layers interleaved by a non-linear activation function, for which we choose GELU\supercite{hendrycks2016gaussian} in this paper. The output of the double convolution module is concatenated with the positional encoding, for which we simply take the Cartesian coordinates of the grid, and then passed to a shallow convolution module (SCM), in which we use two stacks of $3\times 3$ and $1\times 1$ convolutional layers\supercite{cho2021rethinking}.

Moreover, the output of the EB block is combined with the feature extracted from the downsampled input, realizing the fusion of features in different scales. Specifically, the output of the EB block at the $k-1$-th level $\text{EB}_{k-1}^{\text{out}}$ is fed to a convolution layer with a stride of 2 and a doubled number of output channel, resulting in $(\text{EB}_{k-1}^{\text{out}})^{\downarrow}$, which is of the same size as $\text{EB}_{k}^{\text{out}}$ at the $k$-th level. The two tensors $(\text{EB}_{k-1}^{\text{out}})^{\downarrow}$ and $\text{EB}_{k}^{\text{out}}$ are thus concatenated and passed to a $1\times 1$ convolution layer, realizing the integration of the feature extracted from both scales.

\subparagraph{Fourier Layers}

Three separate Fourier layers defined in \eqref{FNO update} serve as the skip connections in UNO, transforming the outputs of the MIMO encoder at different scales from the source term space to the solution space. Note that the channel dimension $C$ differs in the three scales. The scale with finer mesh corresponds to fewer channels, which further improves efficiency by reducing the memory usage and computation cost in larger scales.

\subparagraph{Multiple-Input Single-Output (MISO) Decoder}

The MISO decoder module takes the three outputs from the Fourier layers at different scales as input. The output of the Fourier layers at the $k$-th level $\text{FL}_{k}^{\text{out}}$ is fed to a convolution transpose layer with a stride of 2 and a halved number of output channel, resulting in $(\text{FL}_{k-1}^{\text{out}})^{\uparrow}$, which is of the same size as $\text{FL}_{k-1}^{\text{out}}$ at the $k-1$-th level. The two tensors $(\text{FL}_{k}^{\text{out}})^{\uparrow}$ and $\text{FL}_{k-1}^{\text{out}}$ are then concatenated and passed to a $1\times 1$ convolution layer. In this way, the features extracted at different scales are merged together and finally passed to a decoding block, for which we simply choose a $3\times 3$ convolution layer in this paper.

\subsection{Training Process}

In the training process, suppose the training set contains $N_{\text{train}}$ tuples of $(q, f, u)$ sampled from the same distribution. The data loss function evaluates the discrepancy between the exact solution and the solution obtained by NSNO, which is given by
\begin{equation}
    \mathcal{L}_{\text{data}}=\frac{1}{N_{\text{train}}}\sum_{i=1}^{N_{\text{train}}} \Vert u_i-\hat{u}_i\Vert_{L^2(\it\Omega)},
\end{equation}
where $u_i=\mathcal{S}(q_i, f_i)$ is the exact solution to \eqref{op1} corresponding to $(q_i, f_i)$, and $\hat{u}_i=\text{NSNO}(q_i, f_i)$ is the numerical solution obtained by the proposed NSNO. Furthermore, to avoid overfitting and improve the generalization ability of the model, we also introduce the physics-informed loss to minimize the violation of the Helmholtz equation \eqref{op1}:
\begin{equation}
    \mathcal{L}_{\text{pde}}=\frac{1}{N_{\text{train}}}\sum_{i=1}^{N_{\text{train}}} \Vert \Delta \hat{u}_i+k^2(1+q_i)\hat{u}_i-f_i\Vert_{L^2(\it\Omega)},
\end{equation}
where $\Delta\hat{u}_i$ is computed by a five-point finite difference scheme. The total loss function for training the propsed NSNO is defined as
\begin{equation}\label{loss}
\mathcal{L}_{\text{total}}=\mathcal{L}_{\text{data}}+\lambda \mathcal{L}_{\text{pde}},
\end{equation}
where $\lambda$ is the weight balancing the two loss functions.

\section{Experiments}\label{exp}

In this section, a series of experiments are conducted under various data distributions and parameters settings to evaluate the performance of the proposed NSNO. Firstly we present the experiment setup and dataset generation. Following that, we show the benchmark results on the Helmholtz equation solution operator learning problem to demonstrate the effectiveness of the proposed NSNO. Besides, we discuss the necessity of introducing the physics-informed loss and the convergence properties of the Neumann series.

\subsection{Experiment Setup}\label{exsetup}
The basic settings of our experiments are listed below.
\begin{itemize}
    \item \textbf{Domain discretization:} In the experiments, the spatial domain $\it\Omega$ is set as $\it\Omega$$=[0, 1]^2$ and is discretized uniformly to a $256\times 256$ grid.
    \item \textbf{Model hyperparameters:} Unless otherwise specified, for the Fourier layers in FNO, we set $k_{\max}=12$ and channel dimension $C=32$ as is in the original FNO paper, and the number of iterations is set as $L=4$. For the three Fourier layers at different scales in UNO, we also set $k_{\max}=12$ while the channel dimension is set as $C=8, 16, 32$ from the fine level to the coarse level, respectively. Besides, the number of iterations is set as $L=3$ such that UNO has similar number of parameters as FNO. We have also found that using only three Neumann iterations steps is sufficient to give a satisfactory result.
    \item \textbf{Training hyperparamters:} Unless otherwise specified, the training set has 1000 instances while the test set has 100 testing instances. We use the Adam optimizer\supercite{kingma2014adam} to train the neural network for 500 epochs with an initial learning rate of 0.001 that is halved every 100 epochs. The batchsize is set as 20. The weight $\lambda$ in the loss function \eqref{loss} is set as 0.05. All the experiments are conducted on a single Nvidia V100 GPU with 32GB memory.
    \item \textbf{Performance evaluation:} We use the average relative $L^2$-error on test sets defined as $\frac{1}{N_{\text{test}}}\sum_{i=1}^{N_{\text{test}}}\frac{\Vert u_i-\hat{u}_i\Vert_{L^2(\it\Omega)}}{\Vert u_i\Vert_{L^2(\it\Omega)}}$ to evaluate the performance of the neural operators, where $u_i=\mathcal{S}(q_i, f_i)$ is the exact solution and $\hat{u}_i=\text{NSNO}(q_i, f_i)$ is the numerical solution.
\end{itemize}

\subsection{Benchmark Models and Dataset Generation}

\subsubsection{Benchmark Models}
Note that instead of using the Neumann series to decouple $q$ and $f$, another way is to directly use FNO or UNO to solve the solution operator mapping the tuple $(q, f)$ to $u$. Therefore, by choosing whether to use the Neumann series and the proposed UNO architecture, the following four models are considered.
\begin{itemize}
    \item \textbf{FNO:} Directly use FNO to learn the mapping from $(q, f)$ to $u$.
    \item \textbf{UNO:} Directly use UNO to learn the mapping from $(q, f)$ to $u$.
    \item \textbf{NS-FNO:} NSNO with $G_{\bm{\theta}}$ chosen as FNO.
    \item \textbf{NS-UNO:} NSNO with $G_{\bm{\theta}}$ chosen as UNO.
\end{itemize}

\subsubsection{Dataset generation}
We generate various datasets for the coefficient $q$ and source field $f$ to give a thorough evaluation on the performance of the proposed models. For each dataset, we generate $q$ and $f$ separately from one of the following distributions.

\begin{figure}
\centering
\subfigure[Illustration of the T-shaped $q$]{
\label{illustrationT}
\includegraphics[width=0.3\textwidth]{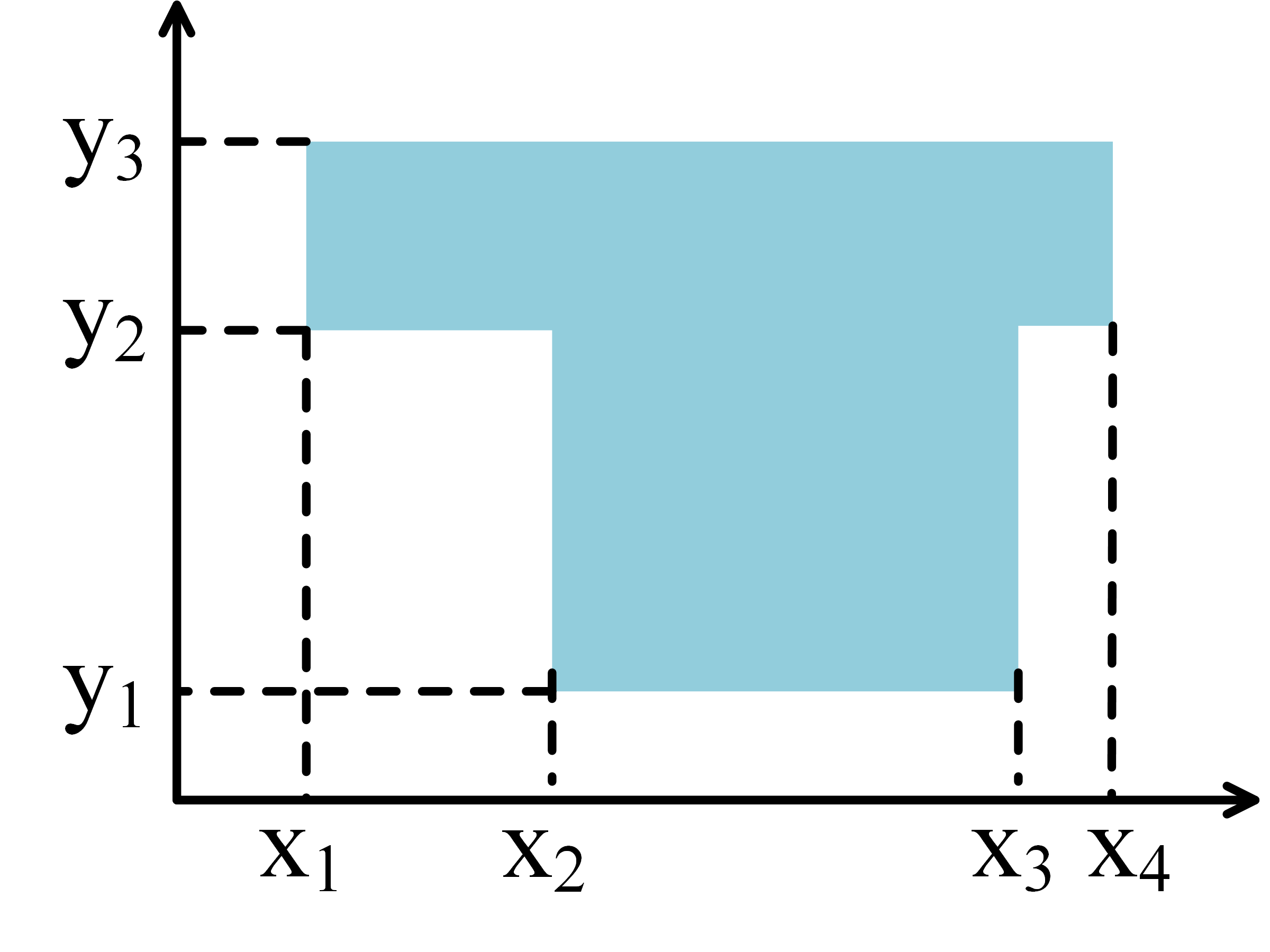}}
\hspace{0.04in}
\centering
\subfigure[T-shaped, example 1 \quad]{
\includegraphics[width=0.3\textwidth]{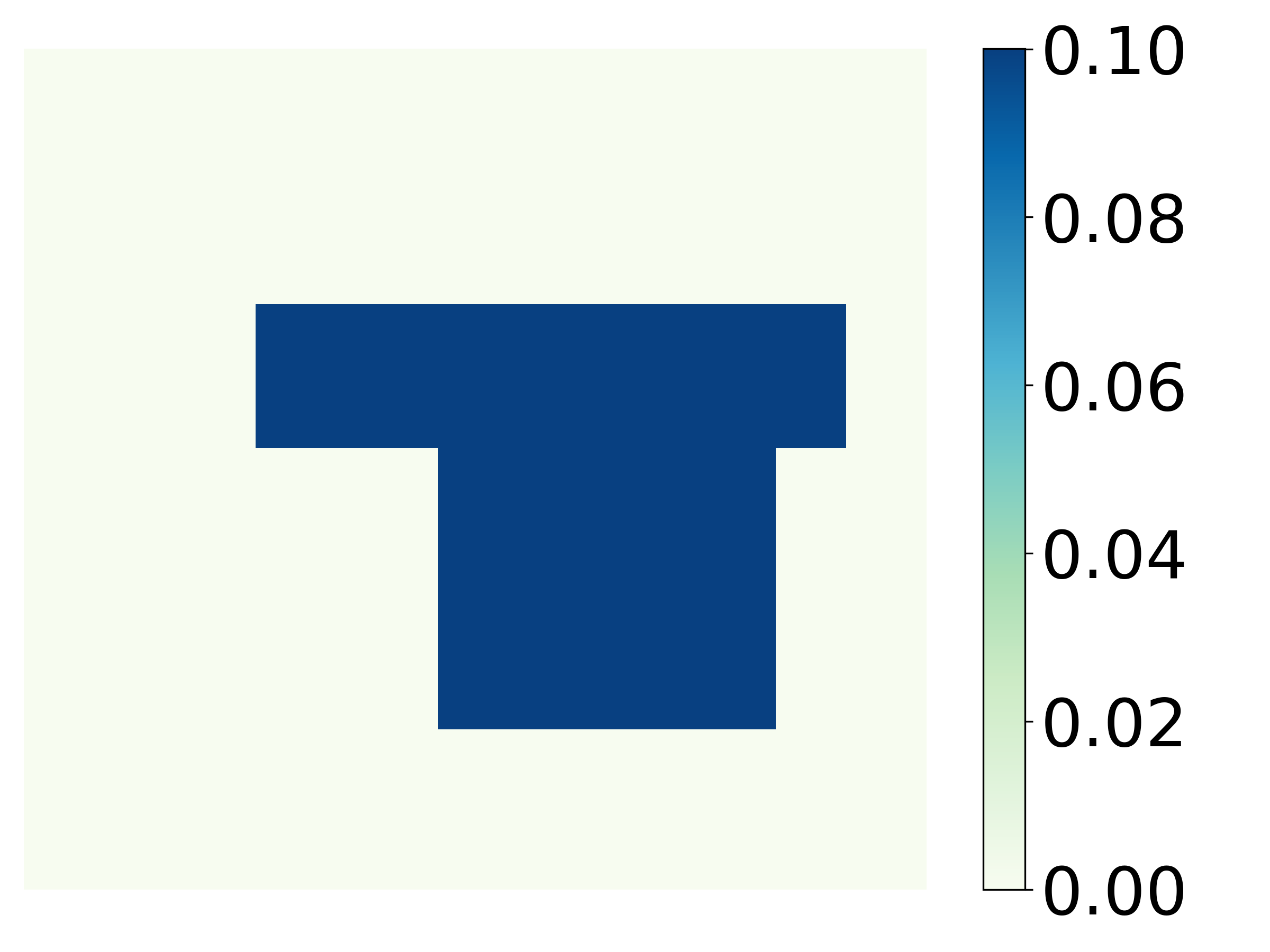}}
\hspace{0.04in}
\centering
\subfigure[T-shaped, example 2 \quad]{
\includegraphics[width=0.3\textwidth]{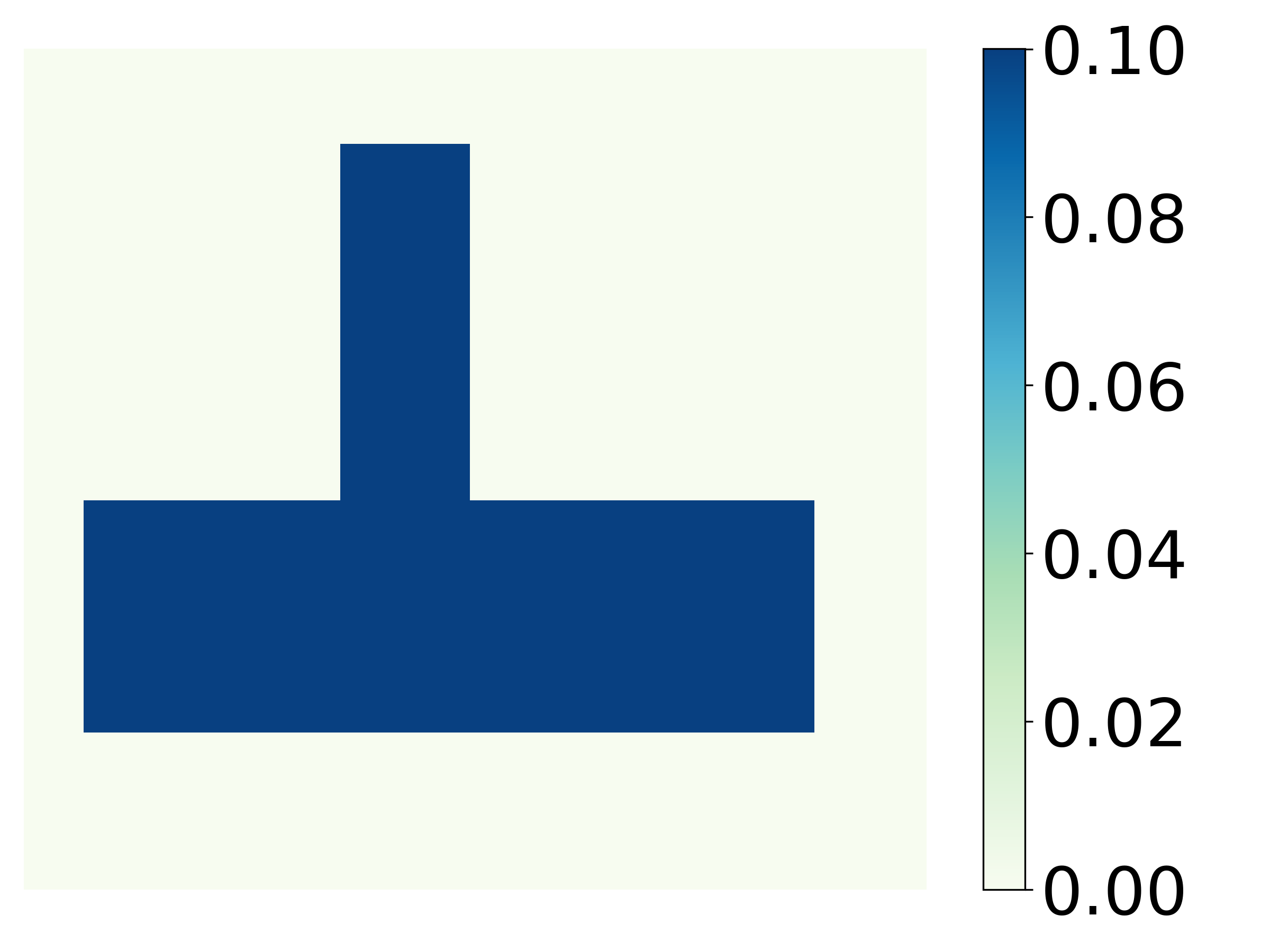}}

\subfigure[Random Circle, example 1]{
\includegraphics[width=0.3\textwidth]{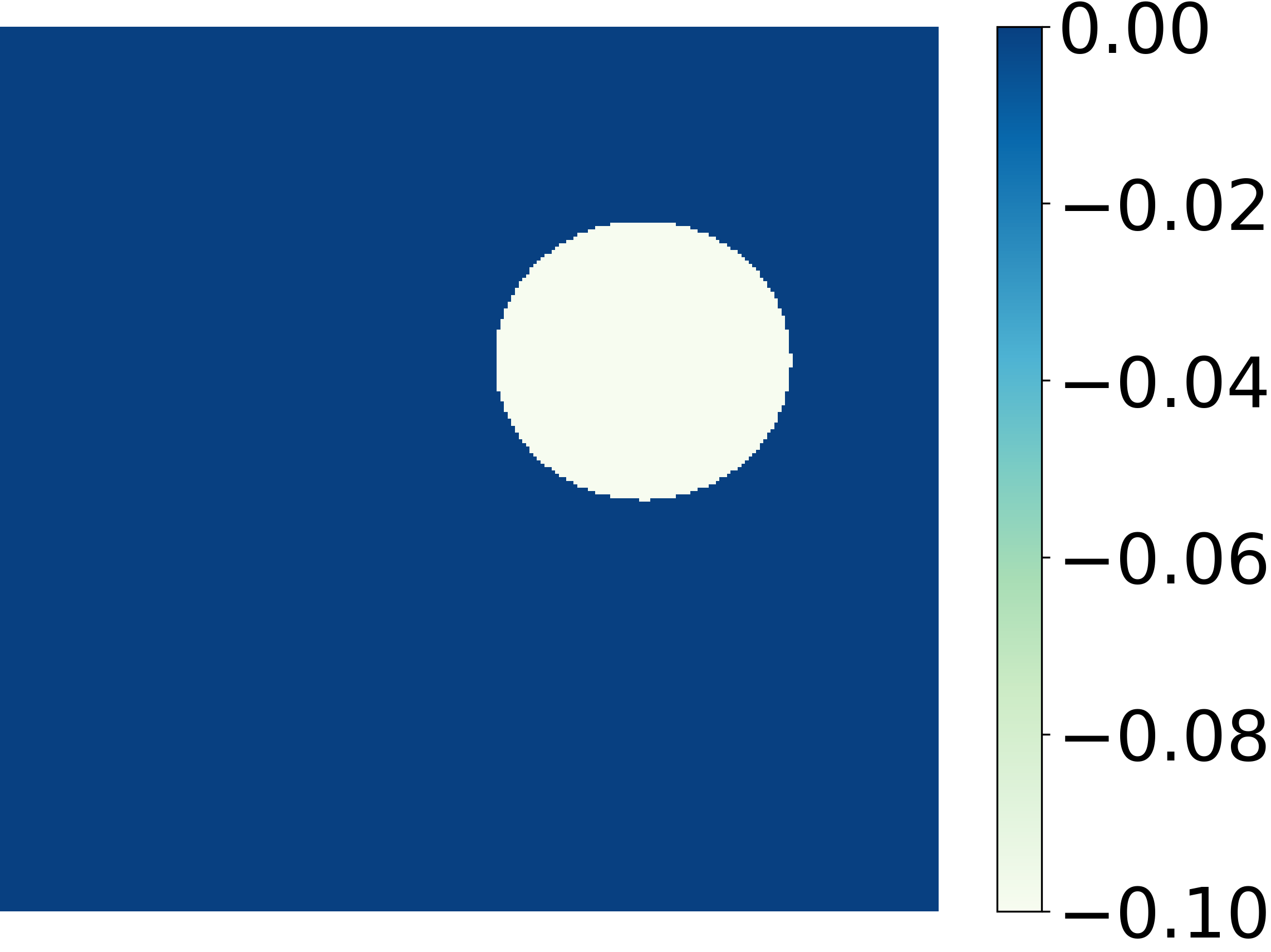}}
\hspace{0.04in}
\centering
\subfigure[Random Circle, example 2]{
\includegraphics[width=0.3\textwidth]{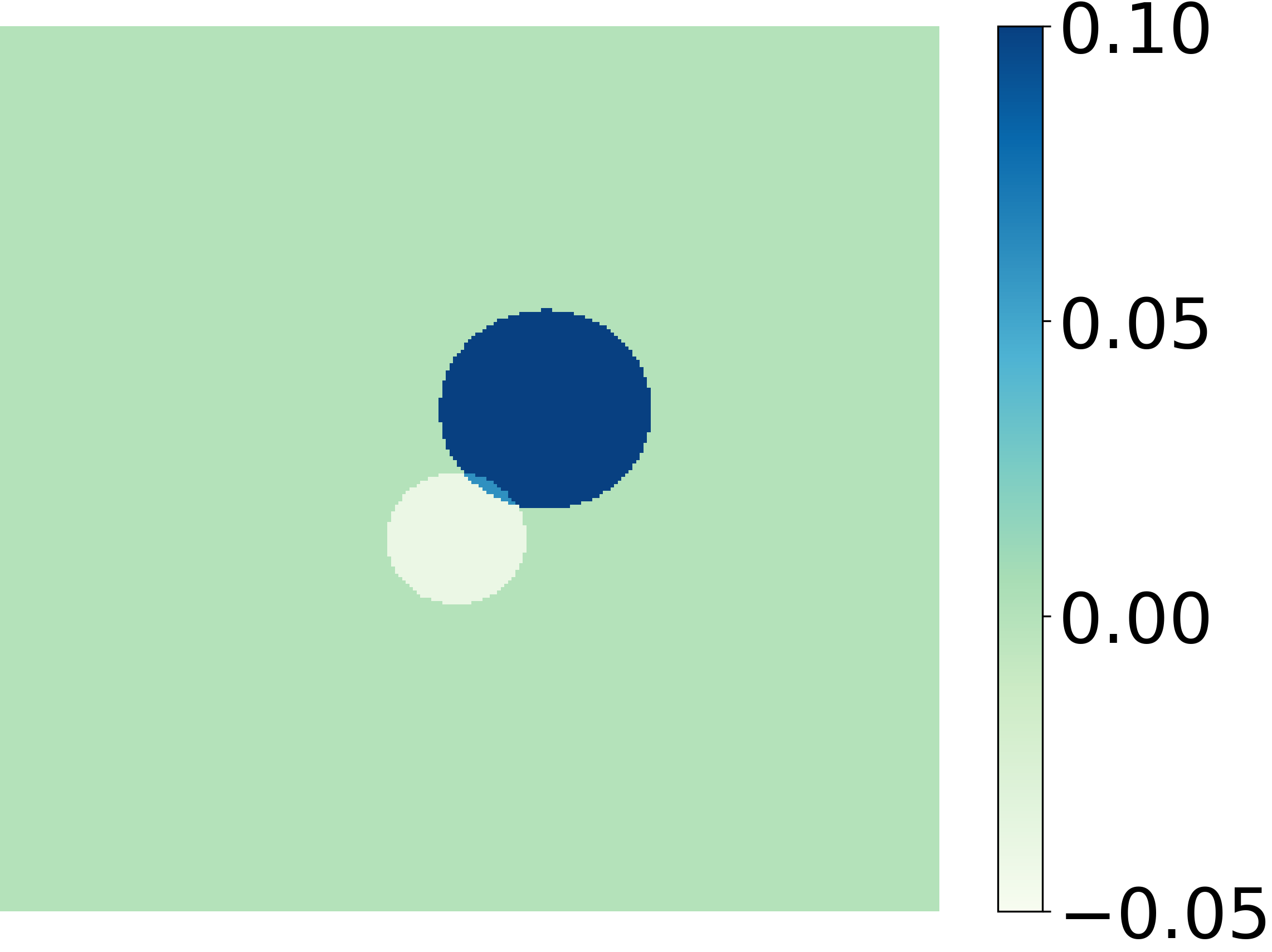}}
\hspace{0.04in}
\centering
\subfigure[Random Circle, example 3]{
\includegraphics[width=0.3\textwidth]{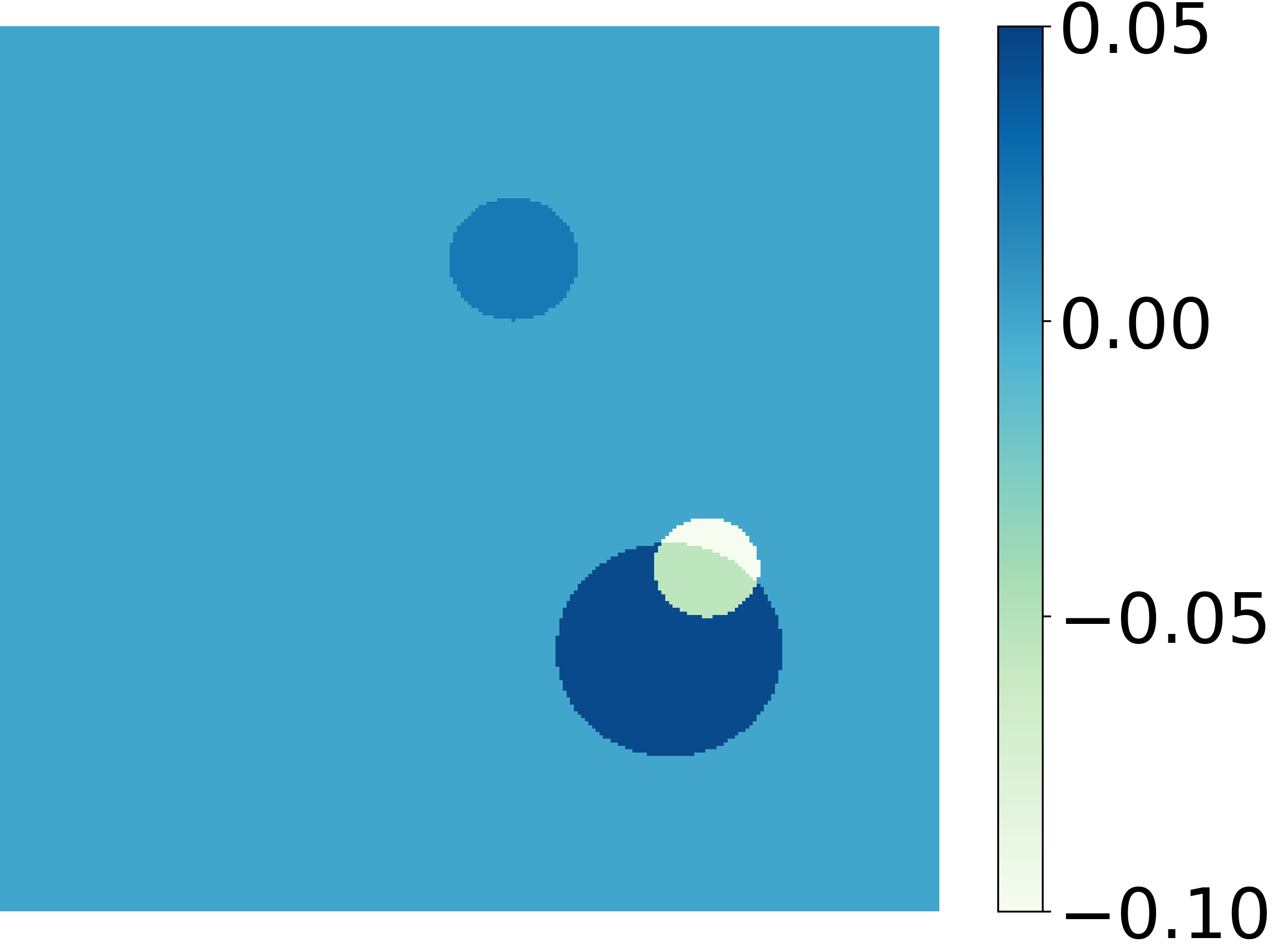}}

\subfigure[Smoothed Circle, example 1]{
\includegraphics[width=0.3\textwidth]{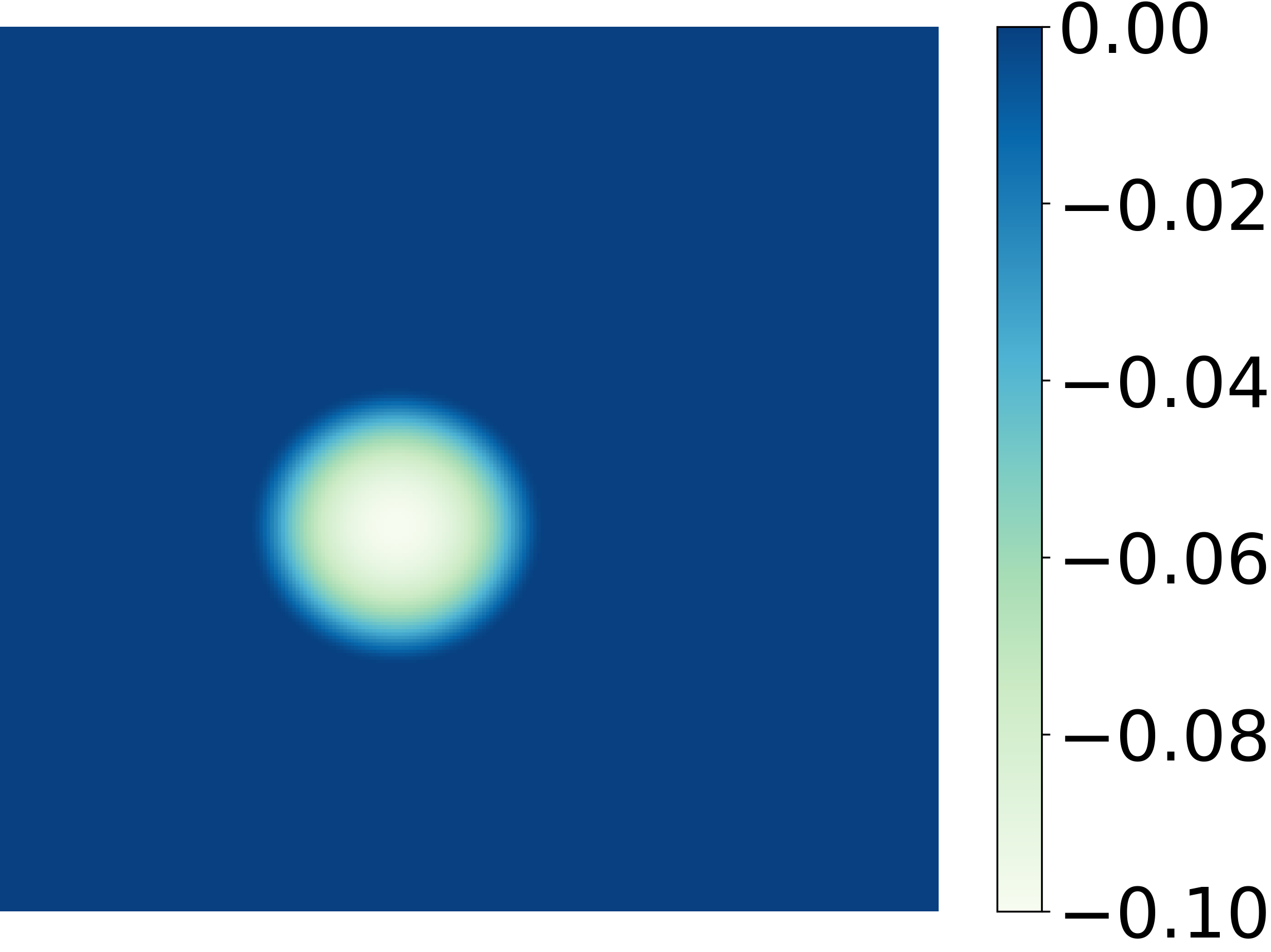}}
\hspace{0.04in}
\centering
\subfigure[Smoothed Circle, example 2]{
\includegraphics[width=0.3\textwidth]{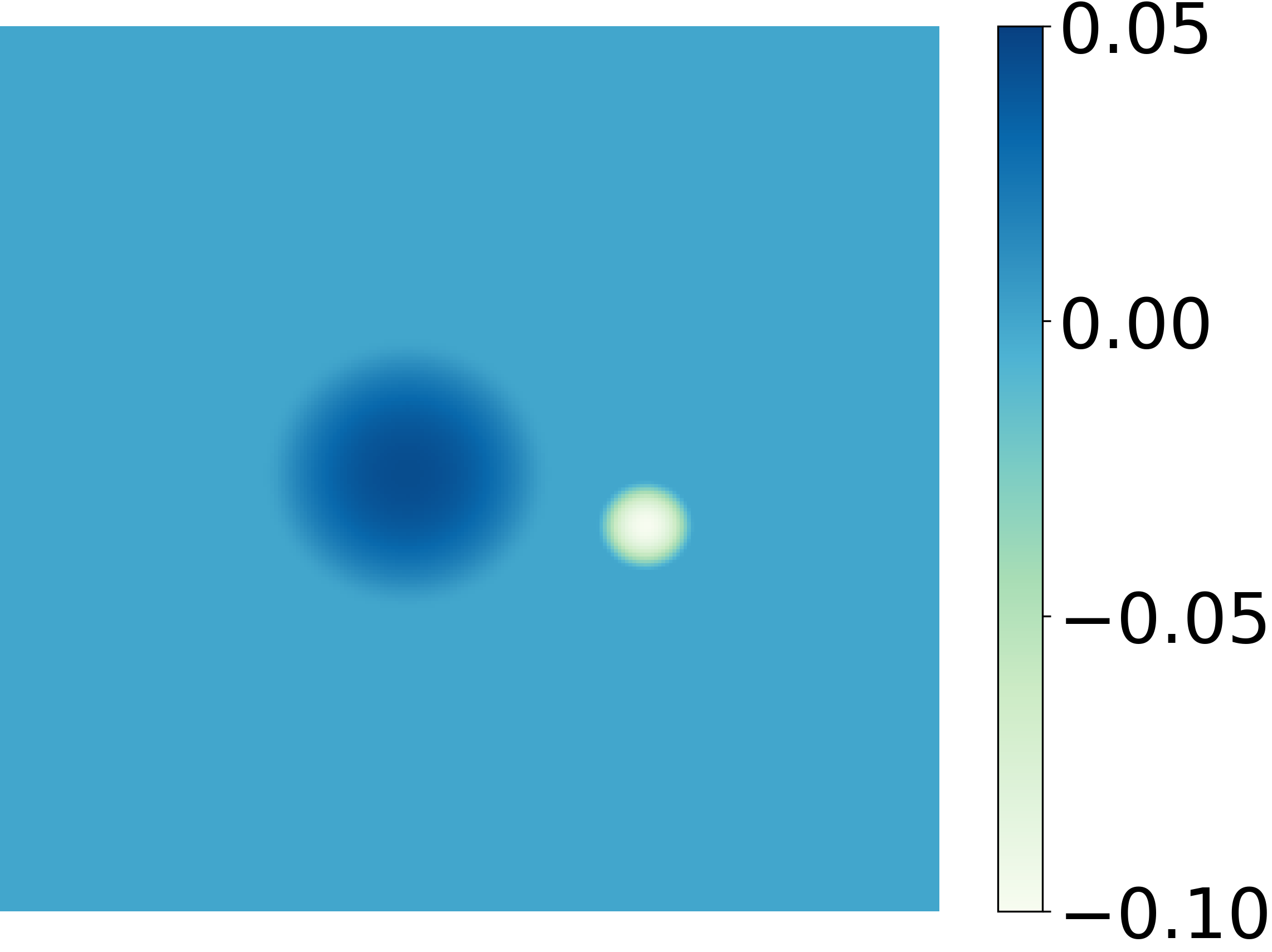}}
\hspace{0.04in}
\centering
\subfigure[Smoothed Circle, example 3]{
\includegraphics[width=0.3\textwidth]{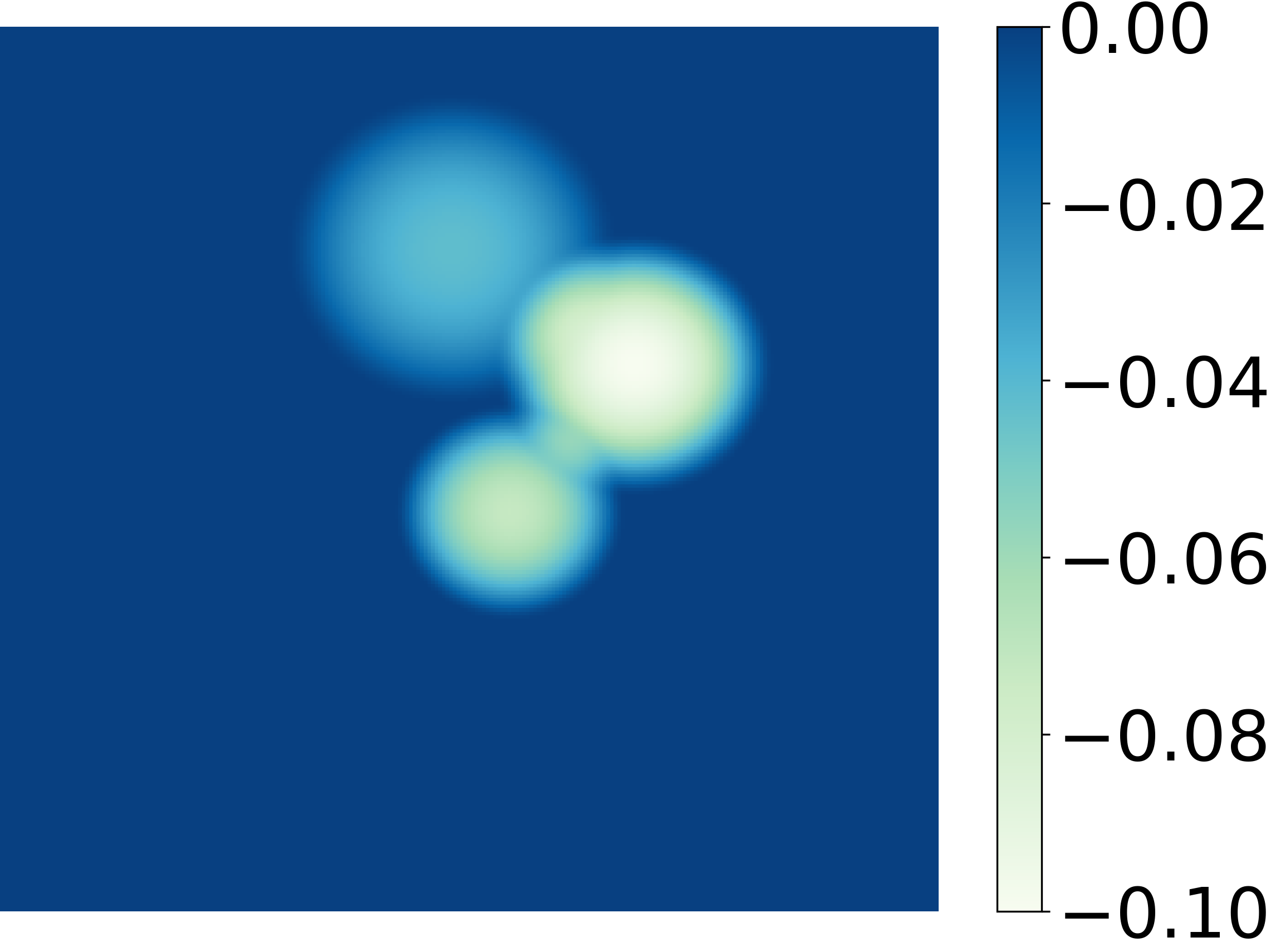}}
\caption{Examples of $q$. (a)-(c): T-shaped distribution. (d)-(f): Random circle distribution with the number of circles from 1 to 3. (g)-(i): Smoothed random circle distribution with the number of circles from 1 to 3.}
\label{example_q}
\end{figure}

The three distributions for $q$ are listed below, along with examples sampled from each distribution shown in Fig. \ref{example_q}.
\begin{itemize}
    \item \textbf{T-shaped:} As is illustrated in Fig. \ref{illustrationT}, for T-shaped distribution, $q$ is compactly supported with randomly generated T-shaped support. Specifically, four points are generated uniformly in [0.05, 0.95] and sorted, denoted as $x_1<x_2<x_3<x_4$ as the coordinates of boundary points in the $x$-direction, then three points are generated uniformly in [0.05, 0.95] and sorted, denoted as $y_1<y_2<y_3$ as the coordinates of boundary points in the $y$-direction. The support of $q$ is then taken as $[x_2, x_3]\times[y_1, y_2]\cup [x_1, x_4]\times[y_2, y_3]$. Finally, we rotate $q$ by an angle chosen randomly in $\{0, \frac{\pi}{2}, \pi, \frac{3\pi}{2}\}$. The function value in the T-shaped support is fixed as 0.1.
    \item \textbf{Random circle\supercite{wei2018deep}:} For the random circle distribution, $q$ is piecewise constant with the support taken as the union of randomly 1-3 circles. Specifically,
    \begin{equation}
        q=\sum_{i=1}^{N_c} \mu_i\chi_{D_i}, \ D_i=\{(x,y): (x-x_i)^2+(y-y_i)^2\leq r_i^2\}, \ N_c\in\{1, 2, 3\},
    \end{equation}
    where $\chi_{D_i}$ is the indicator function on $D_i$, $x_i, y_i\sim U[0.2, 0.8]$, $r_i\sim U[0.05, 0.2]$, $\mu_i\sim U[-1,1]$. Finally $q$ is normalized such that $\Vert q\Vert_{L^{\infty}(\it\Omega)}=0.1$. Note that the circles are allowed to overlap with each other.
    \item \textbf{Smoothed circle:} The smoothed circle distribution is the smoothed version of the random circle distribution such that $q\in C_{0}^{\infty}(\it\Omega)$. Specifically,
    \begin{equation}
        q=\sum_{i=1}^{N_c} \mu_i\chi_{D_i}\exp\left[-\frac{1}{1-\frac{(x-x_i)^2+(y-y_i)^2}{r_i^2}}\right],
    \end{equation}
    where $D_i, N_c, \chi_{D_i}, x_i, y_i, r_i$ and $\mu_i$ are defined the same as the random circle distribution. Finally $q$ is also normalized such that $\Vert q\Vert_{L^{\infty}(\it\Omega)}=0.1$.
\end{itemize}

\begin{figure}
\centering
\setcounter{subfigure}{0}
\subfigure[Gaussian(10) \quad]{
\includegraphics[width=0.3\textwidth]{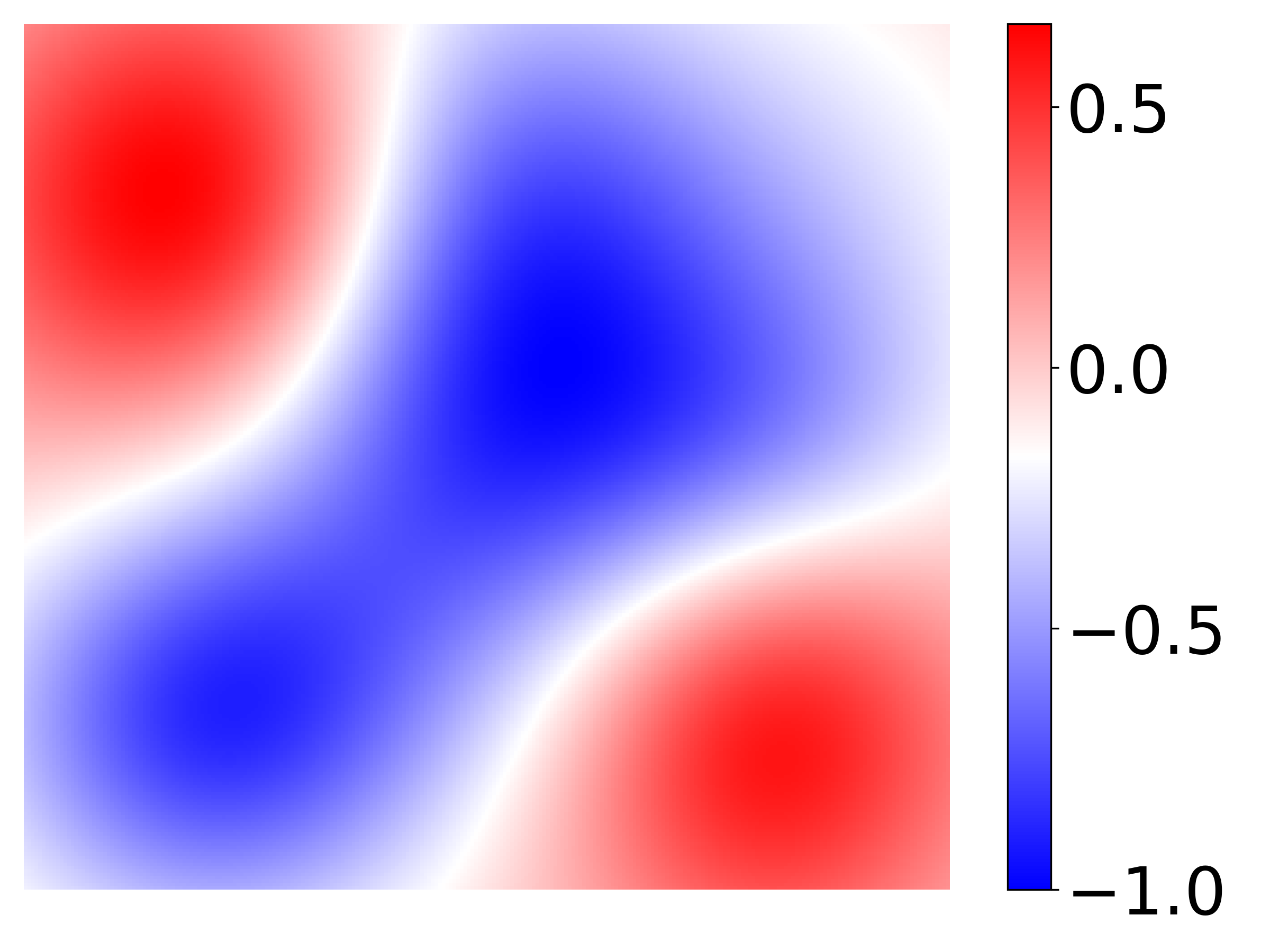}}
\hspace{0.04in}
\centering
\subfigure[Gaussian(30) \quad ]{
\includegraphics[width=0.3\textwidth]{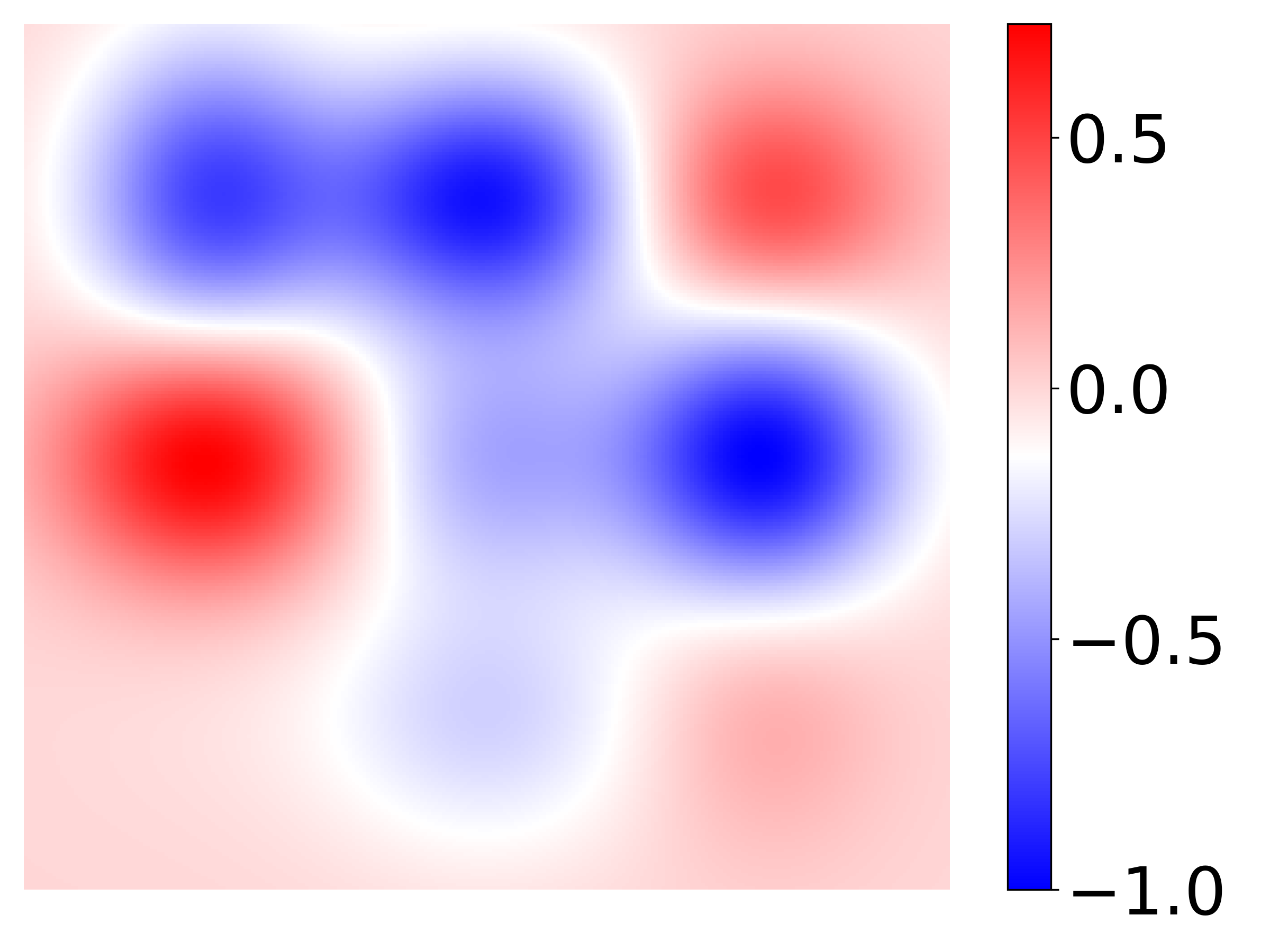}}
\hspace{0.04in}
\centering
\subfigure[Gaussian(50) \quad]{
\includegraphics[width=0.3\textwidth]{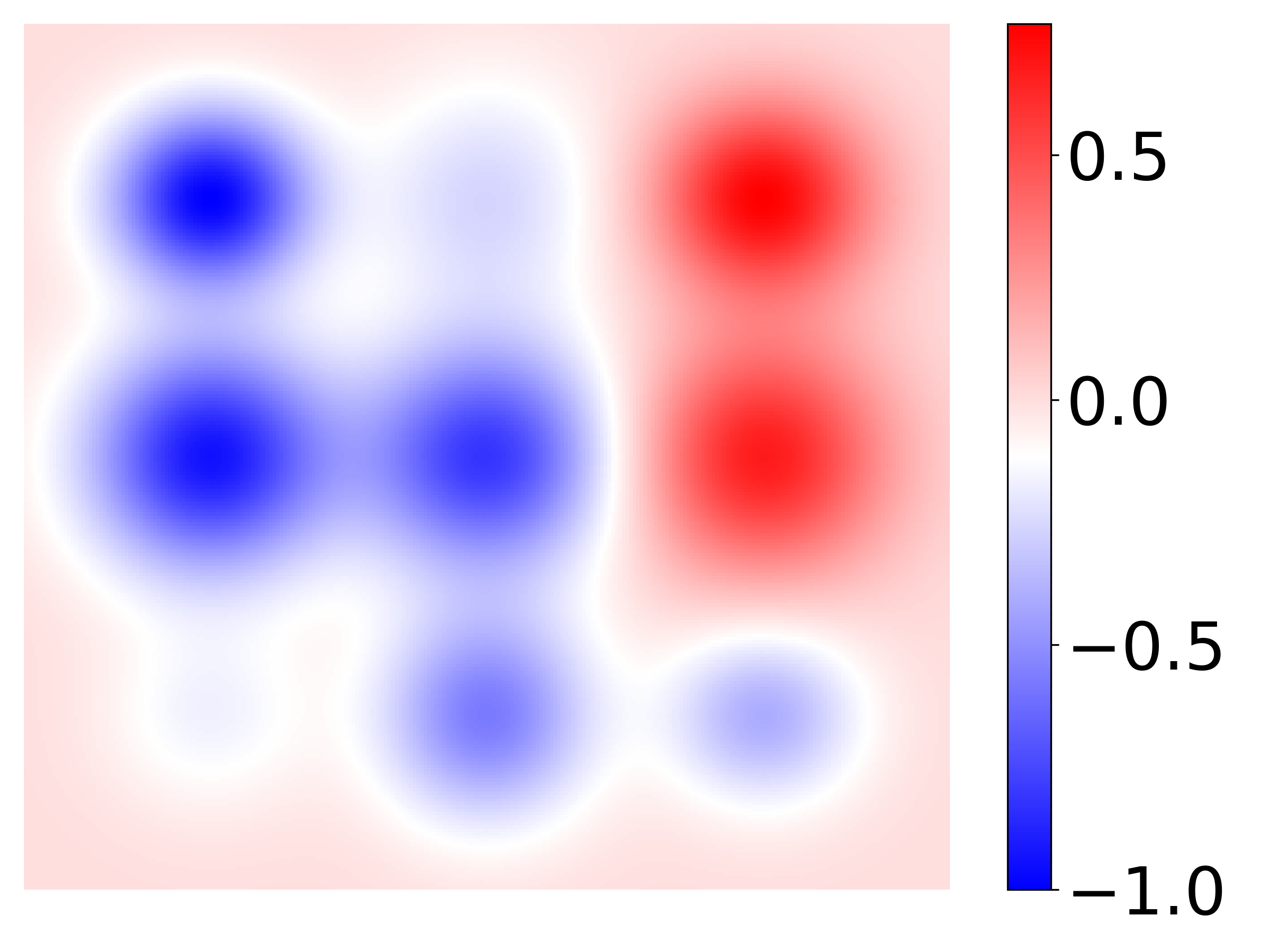}}

\subfigure[GRF]{
\includegraphics[width=0.3\textwidth]{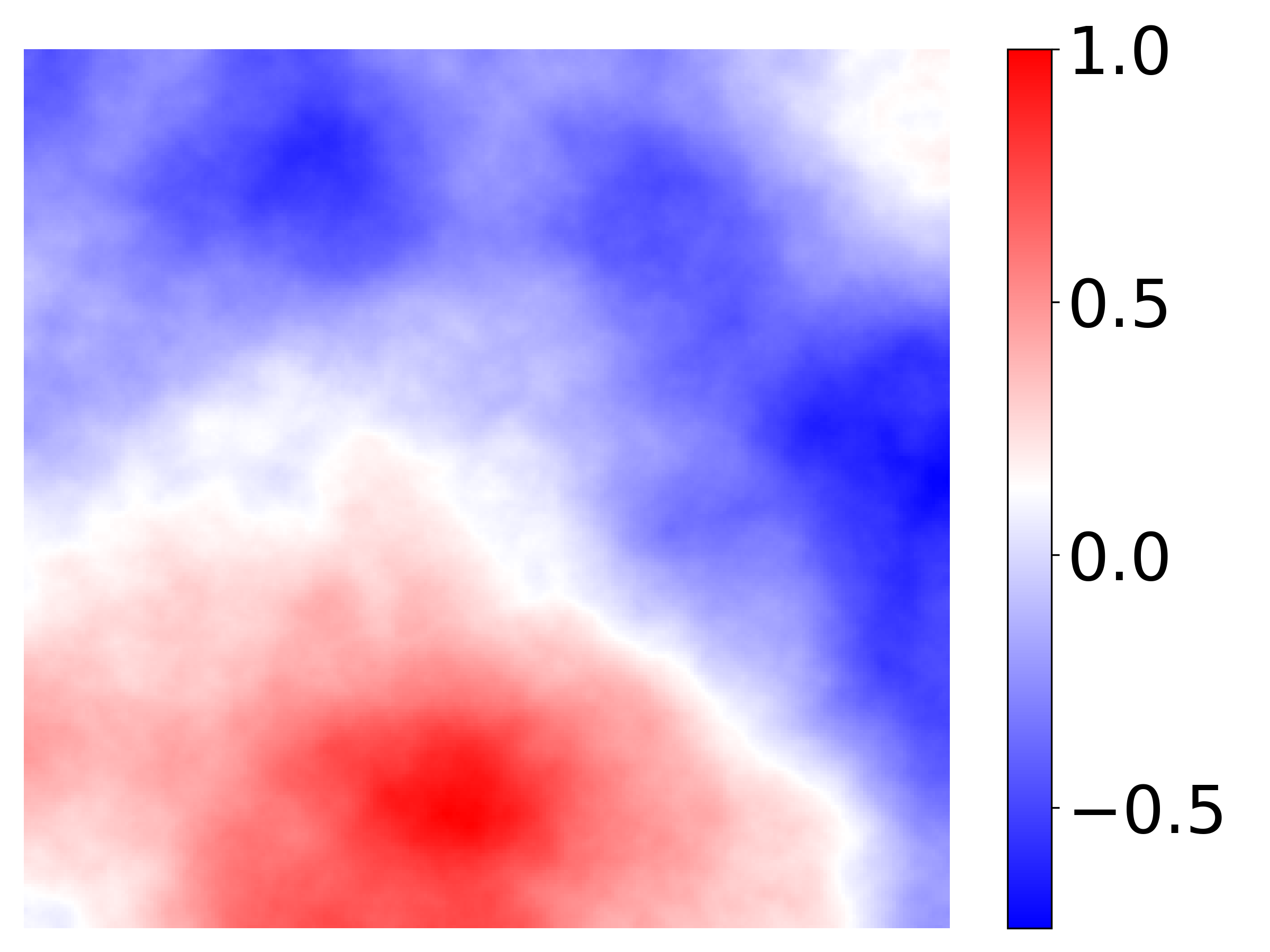}}
\hspace{0.04in}
\centering
\subfigure[Waves]{
\includegraphics[width=0.3\textwidth]{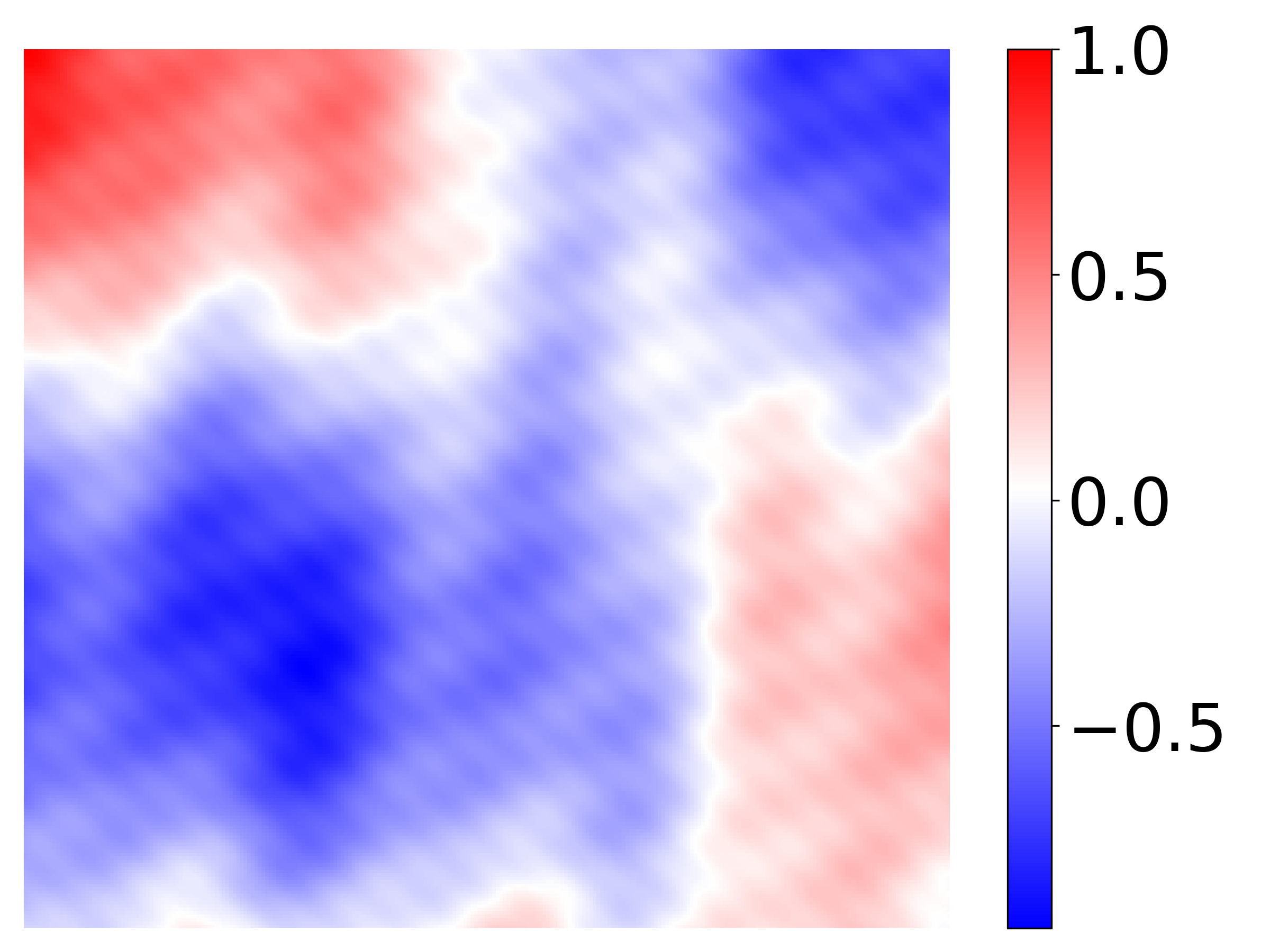}}
\hspace{0.04in}
\centering
\subfigure[Waves]{
\includegraphics[width=0.3\textwidth]{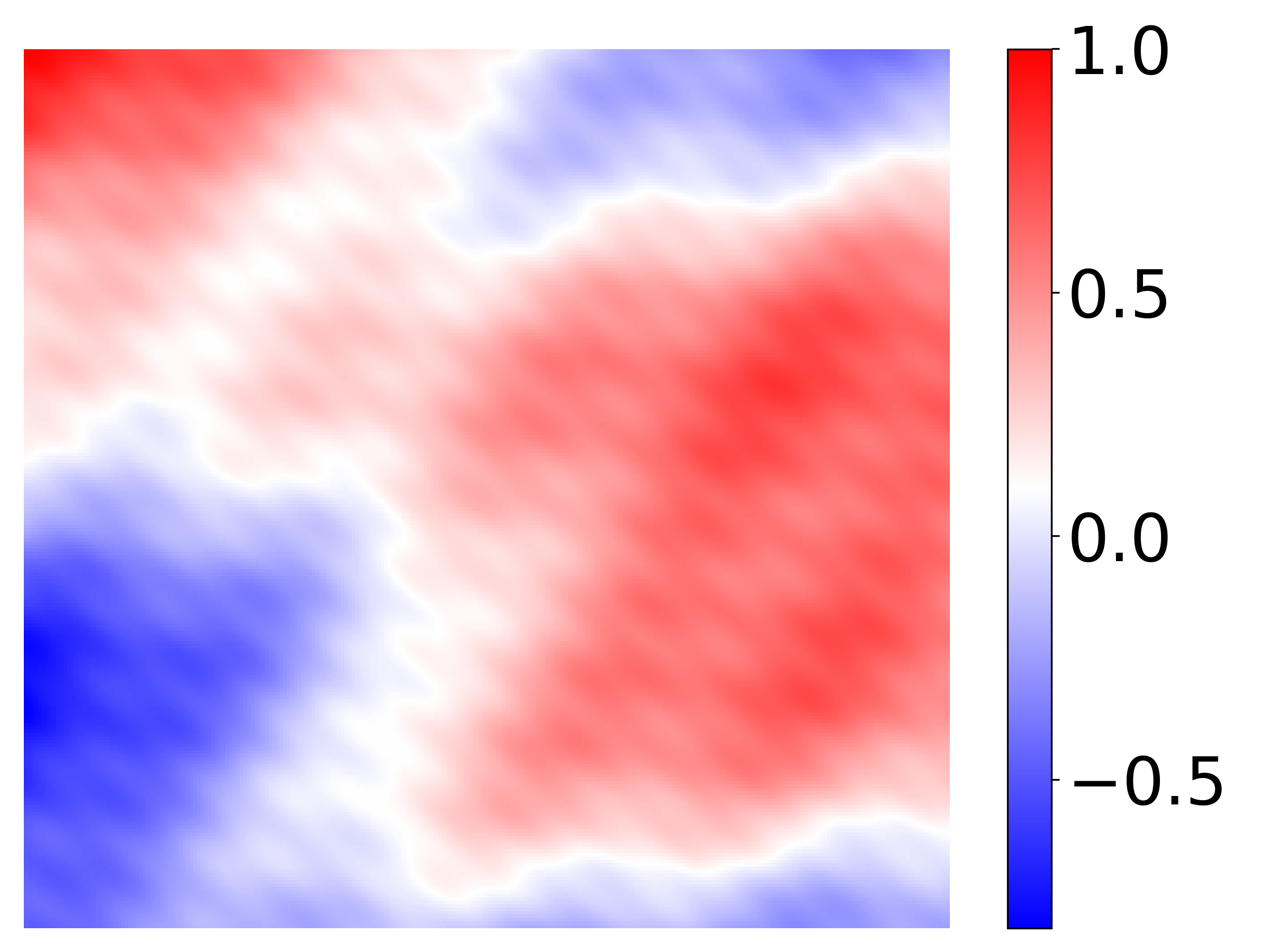}}
\caption{Examples of $f$. (a)-(c): Gaussian distribution with different rates of decay. (d): One example sampled from the Gaussian random field (GRF) distribution. (e)-(f): Two examples sampled from the wave distribution.}
\label{example_f}
\end{figure}

The three distributions for $f$ are listed below, along with examples sampled from each distribution shown in Fig. \ref{example_f}. We normalize each $f$ such that $\Vert f\Vert_{L^{\infty}(\it\Omega)}=1$.
\begin{itemize}
    \item \textbf{Gaussian(R)\supercite{zhang2023belnet}:} $f$ is taken as the sum of nine Gaussians. The centers of the Gaussians are fixed as $c_{i,j}=\left(\frac{3i-1}{10}, \frac{3j-1}{10}\right), i,j=1,2,3$. The decay rates of the Gaussians are uniformly sampled from $[R, 2R]$, where $R$ is an adjustable hyperparameter. Three examples with different decay rates are shown in Fig. \ref{example_f}(a)-(c). It can be seen that the larger the decat rates are, the more compactly supported $f$ is.
    \item \textbf{GRF\supercite{li2020fourier}:} $f$ is generated according to the Gaussian random field (GRF) $\mathcal{N}(0, (-\Delta+9I)^{-2})$ with zero Neumann boundary conditions on the Laplacian.
    \item \textbf{Wave:} The wave distribution is the weighted sum of six planar waves at different frequencies given as:
    \begin{equation}
        f(x,y)=\sum_{i=1}^6 \frac{1}{\mu_i}cos[\pi\mu_i(x\cos\theta_i+y\sin\theta_i)],
    \end{equation}
    where $\mu_i\sim U[2^{i-1}, 1.5\times 2^{i-1}]$, $\theta_i\sim U[0, 2\pi]$. Examples sampled from this distribution has explicit multi-scale property resulted from the superposition of waves at different frequencies.
\end{itemize}

After the coefficients $q$ and source term $f$ are generated, the corresponding exact solutions are then computed by finite difference method where the large linear system is solved by MUMPS\supercite{amestoy2001fully}.

\subsection{Experiment Results}
In this subsection, we first present the benchmark results for $k=20$. Following that, further experiments are carried out to further show the advantages of the NS-UNO in higher wavenumber scenario, data efficiency and computational cost.

\subsubsection{Benchmark Results}
By combining the distributions for coefficient $q$ and source term $f$, we present the benchmark results on the following six datasets as is shown in Table \ref{k20}. It can be seen that compared with FNO and UNO, both NS-FNO and NS-UNO have lower relative $L^2$-error on the test sets, indicating the effectiveness of the proposed Neumann series based framework. Besides, among the four models, NS-UNO achieves the highest accuracy over all datasets, especially on the datasets where $f$ has strong multi-scale characteristics, such as the GRF distribution. Compared with the state-of-the-art FNO, NS-UNO achieves at least 60\% reduction of relative error. Therefore, the proposed NS-UNO network architecture cam better capture the multi-scale features of the solution the Helmholtz equations.

\begin{table}
\centering
\caption{Benchmark relative $L^2$-error ($\times 10^{-2}$) for $k=20$ }
\vspace{0.1cm}
\begin{tabular}{|c|c|cccc|}
\hline
$q$  & $f$ & FNO  & UNO  & NS-FNO & NS-UNO \\
\hline
\multirow{3}{*}{T-shaped} & Gaussian(50) & 3.22 & 2.80 & 1.34 & \textbf{1.26} \\
& Gaussian(30) & 4.78 & 3.10 & 1.72 & \textbf{1.30} \\
& Gaussian(10) & 10.23 & 3.36 & 4.82 & \textbf{1.36} \\
\hline
T-shaped & GRF & 16.80 & 15.94 & 4.87 & \textbf{2.04}\\
\hline
Random circle & GRF & 14.43 & 11.27 & 5.48 & \textbf{1.54}\\
\hline
Smoothed circle & Waves & 7.68 & 6.15 & 3.08 & \textbf{1.43}\\
\hline
\end{tabular}
\label{k20}
\end{table}

Moereover, we observe that as the decay rate of the Gaussians in the Gaussian dataset for $f$ decreasing, the $L^2$-error of FNO-based methods significantly increases, while that of UNO-based methods only slightly increases. The reason is that for $f$ consists of Gaussians with smaller decay rates, the Gaussians decaying at different rates will overlap with each other, as can be seen from Fig. 6(a)-(c), adding more multi-scale features to $f$. Therefore, UNO-based methods which are more capable of handling multi-scale inputs outperform FNO-based methods when the decay rate decreases.

\begin{figure}
\centering
\subfigure[Ground truth, real part \qquad]{
\includegraphics[width=0.3\textwidth]{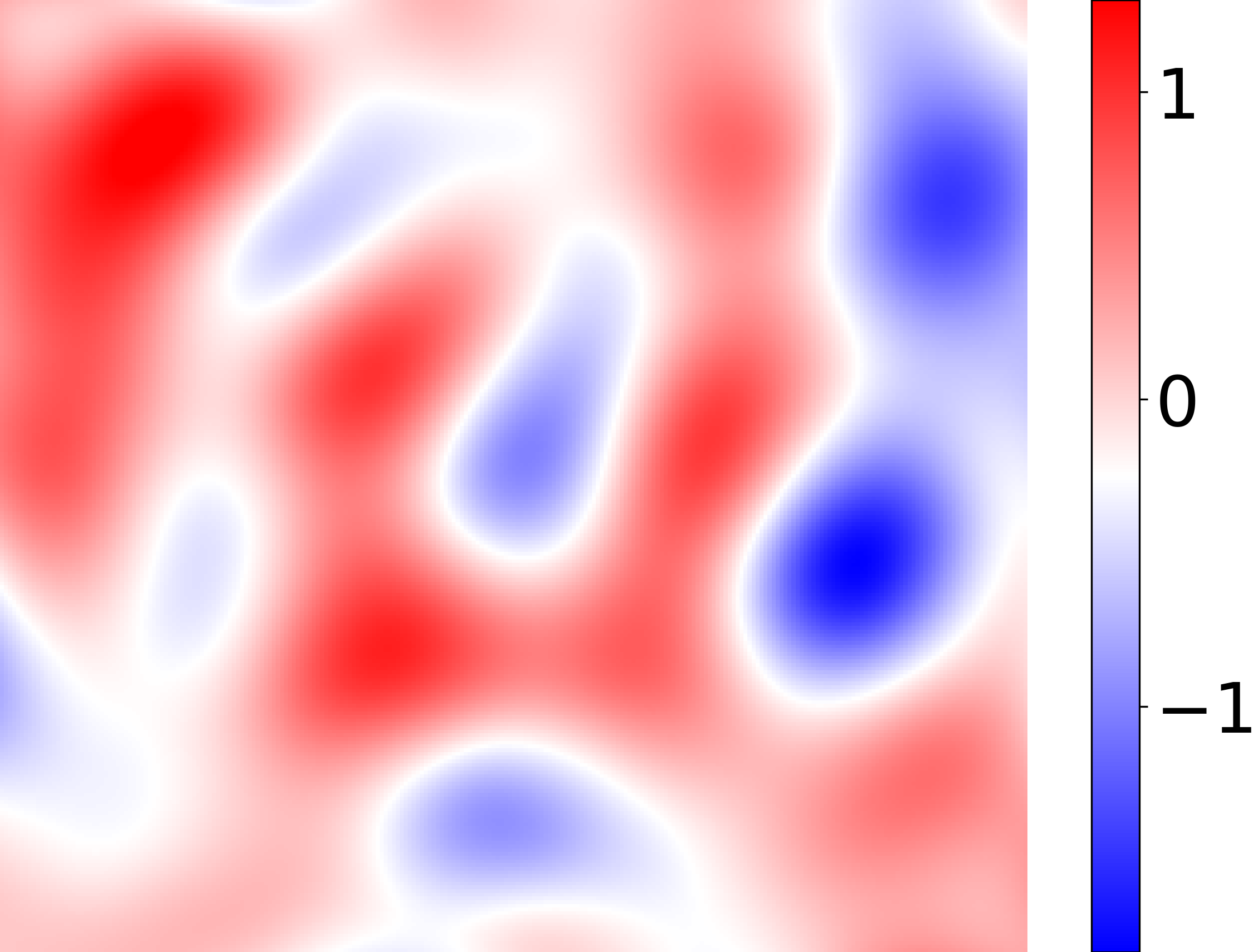}}
\hspace{0.04in}
\centering
\subfigure[Ground truth, imaginary part \qquad]{
\includegraphics[width=0.3\textwidth]{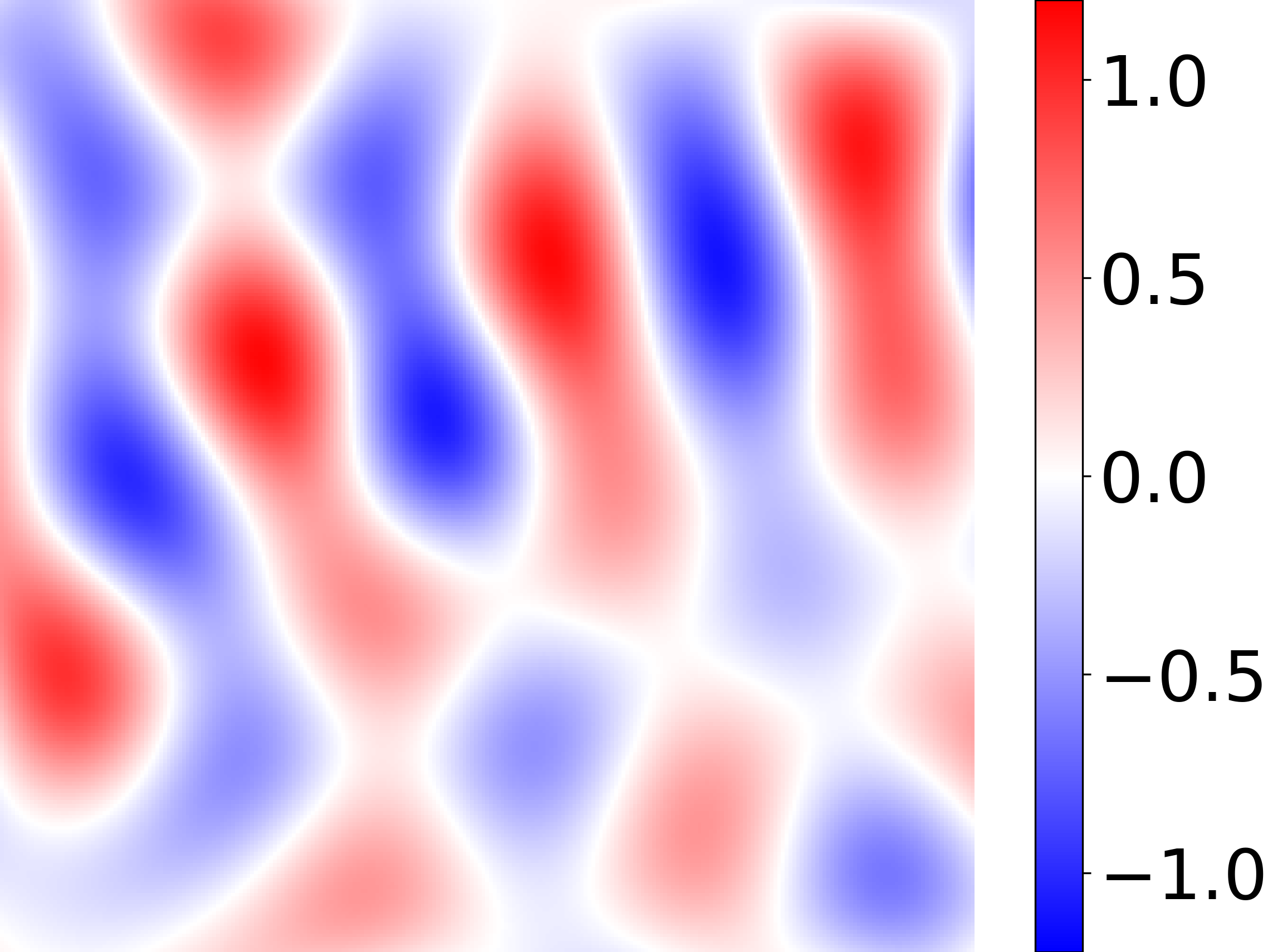}}
\hspace{0.04in}
\centering
\subfigure[$q$ (top left) and $f$\qquad]{
\includegraphics[width=0.3\textwidth]{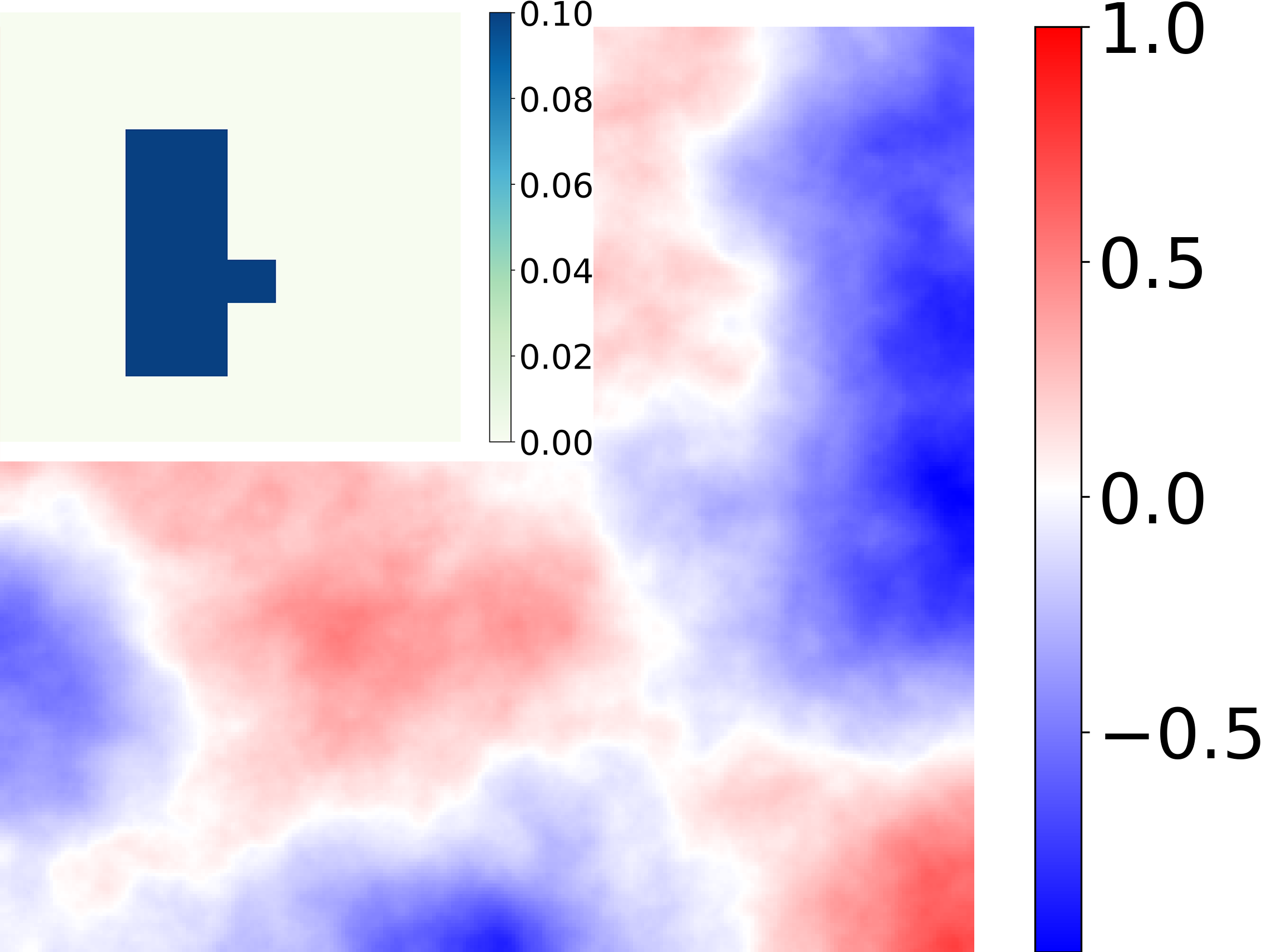}}

\subfigure[FNO, real part \quad]{
\includegraphics[width=0.3\textwidth]{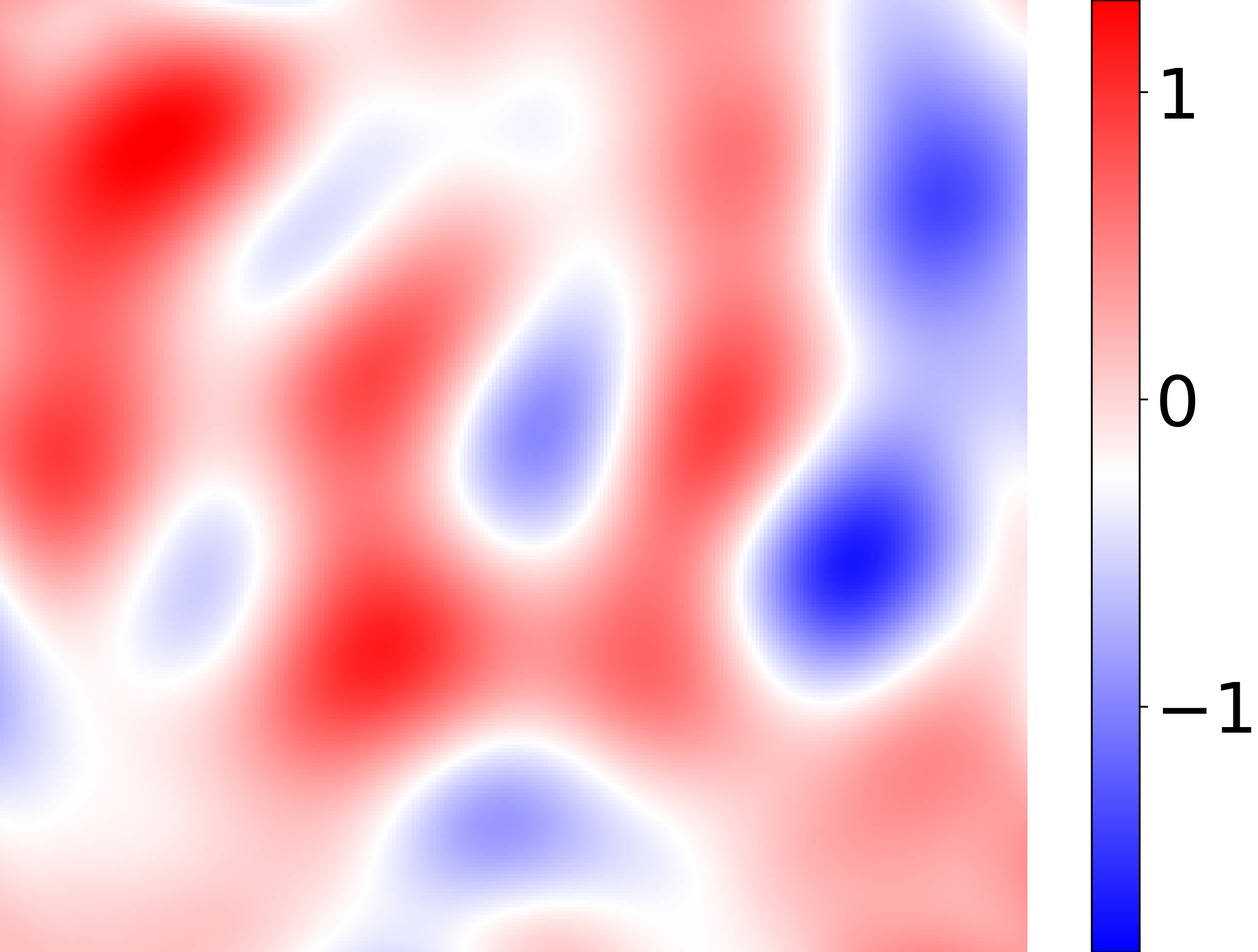}}
\hspace{0.04in}
\centering
\subfigure[FNO, imaginary part\qquad]{
\includegraphics[width=0.3\textwidth]{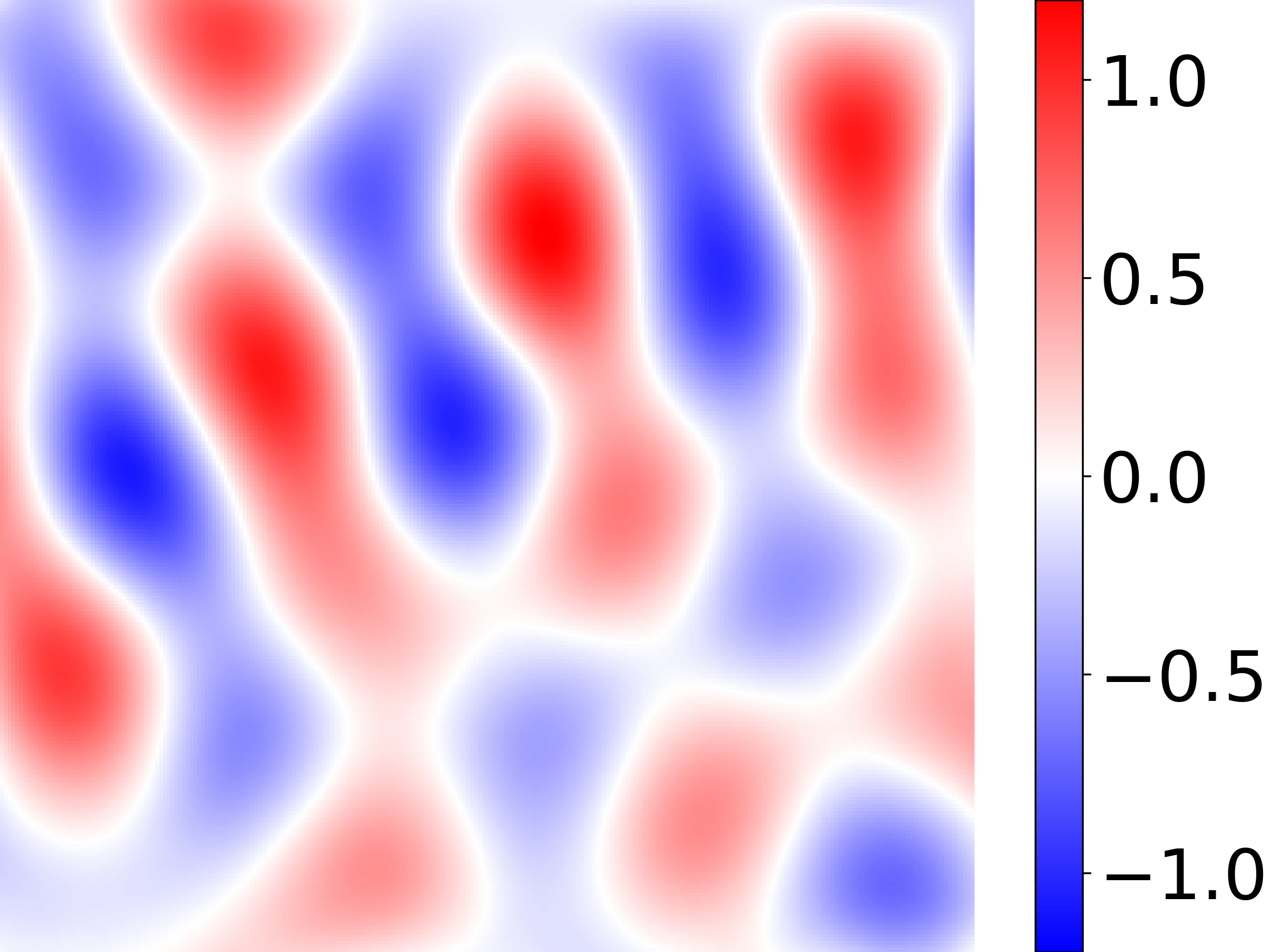}}
\hspace{0.04in}
\centering
\subfigure[FNO, error, 16.26\% \qquad]{
\includegraphics[width=0.3\textwidth]{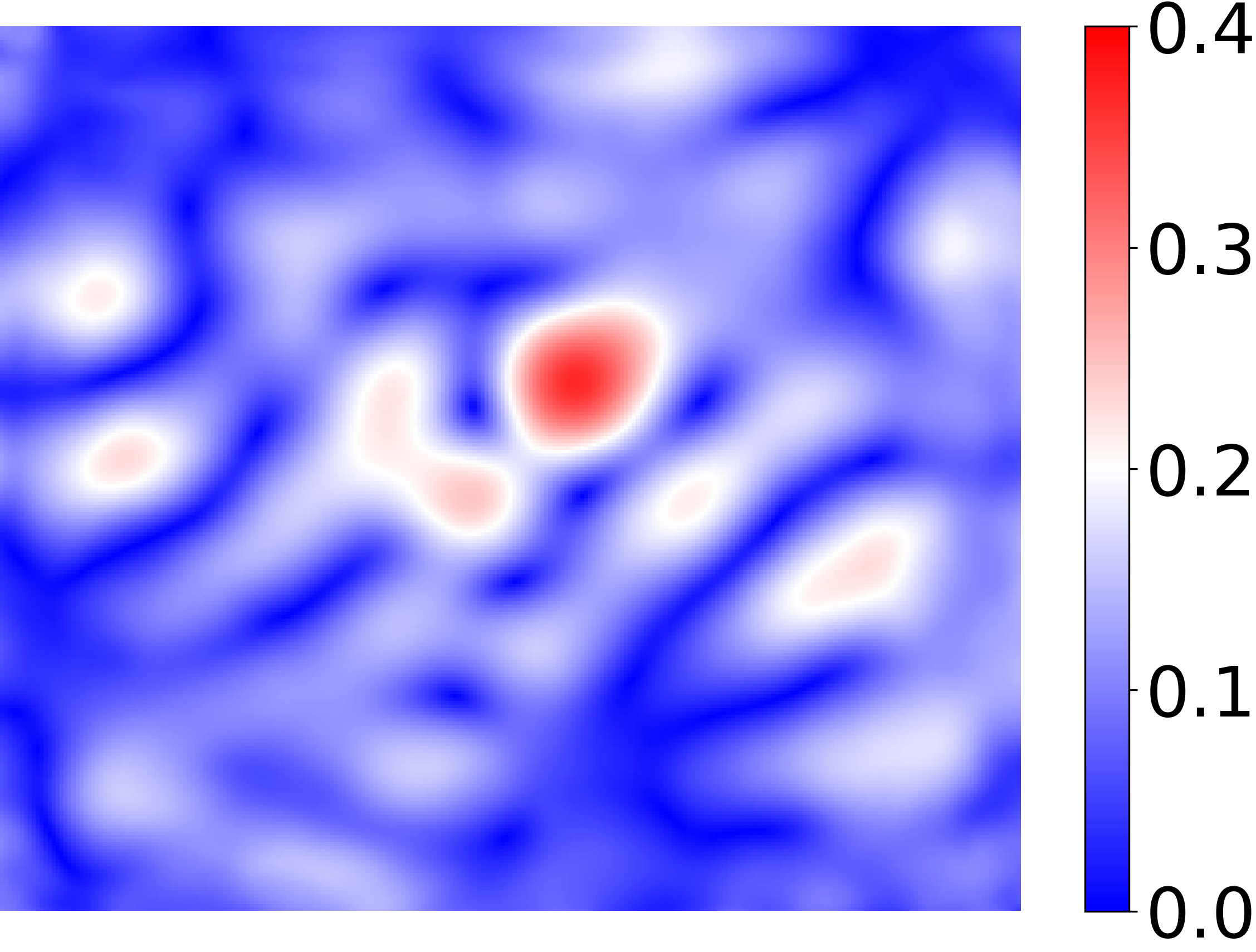}}

\subfigure[UNO, real part \quad]{
\includegraphics[width=0.3\textwidth]{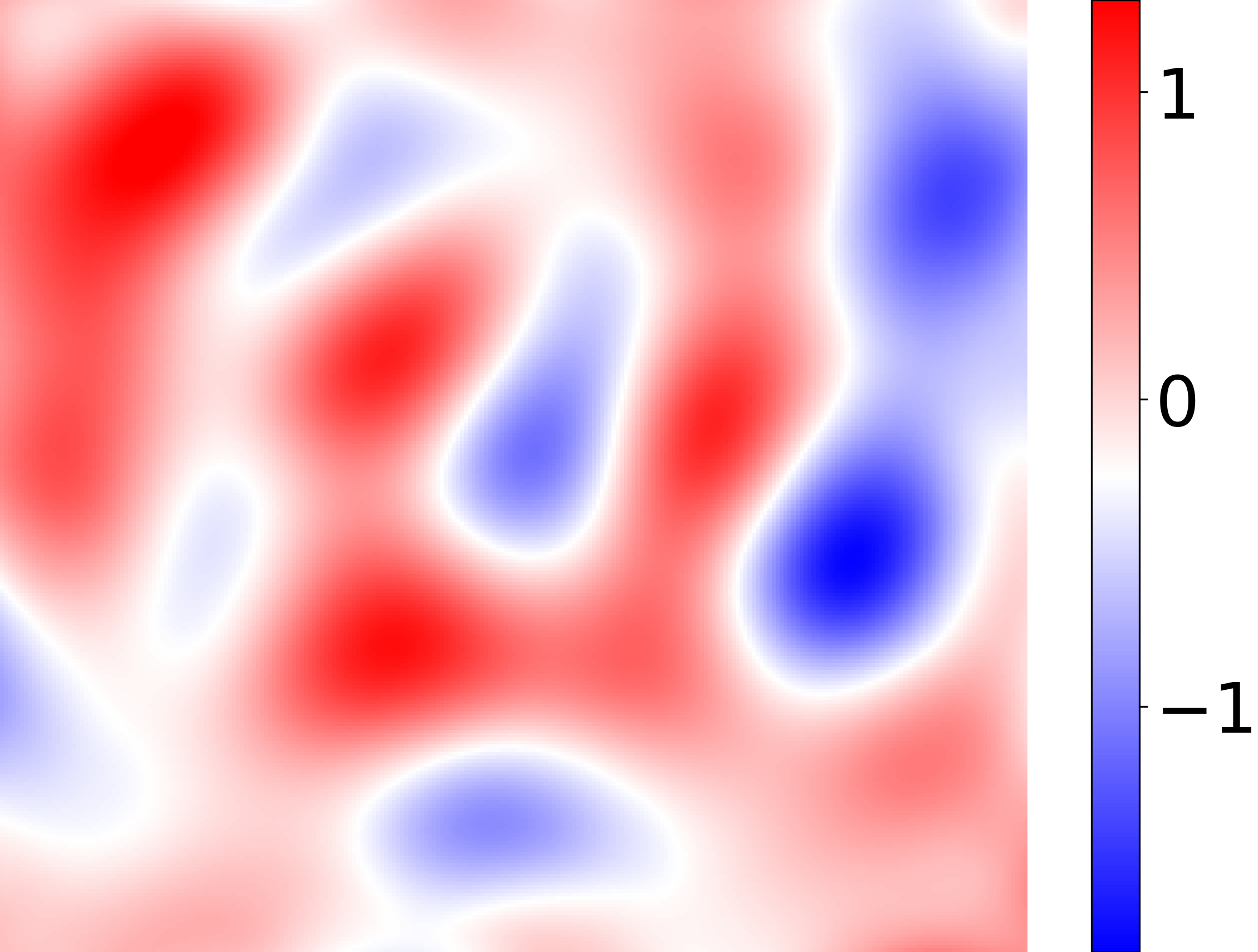}}
\hspace{0.04in}
\centering
\subfigure[UNO, imaginary part\qquad]{
\includegraphics[width=0.3\textwidth]{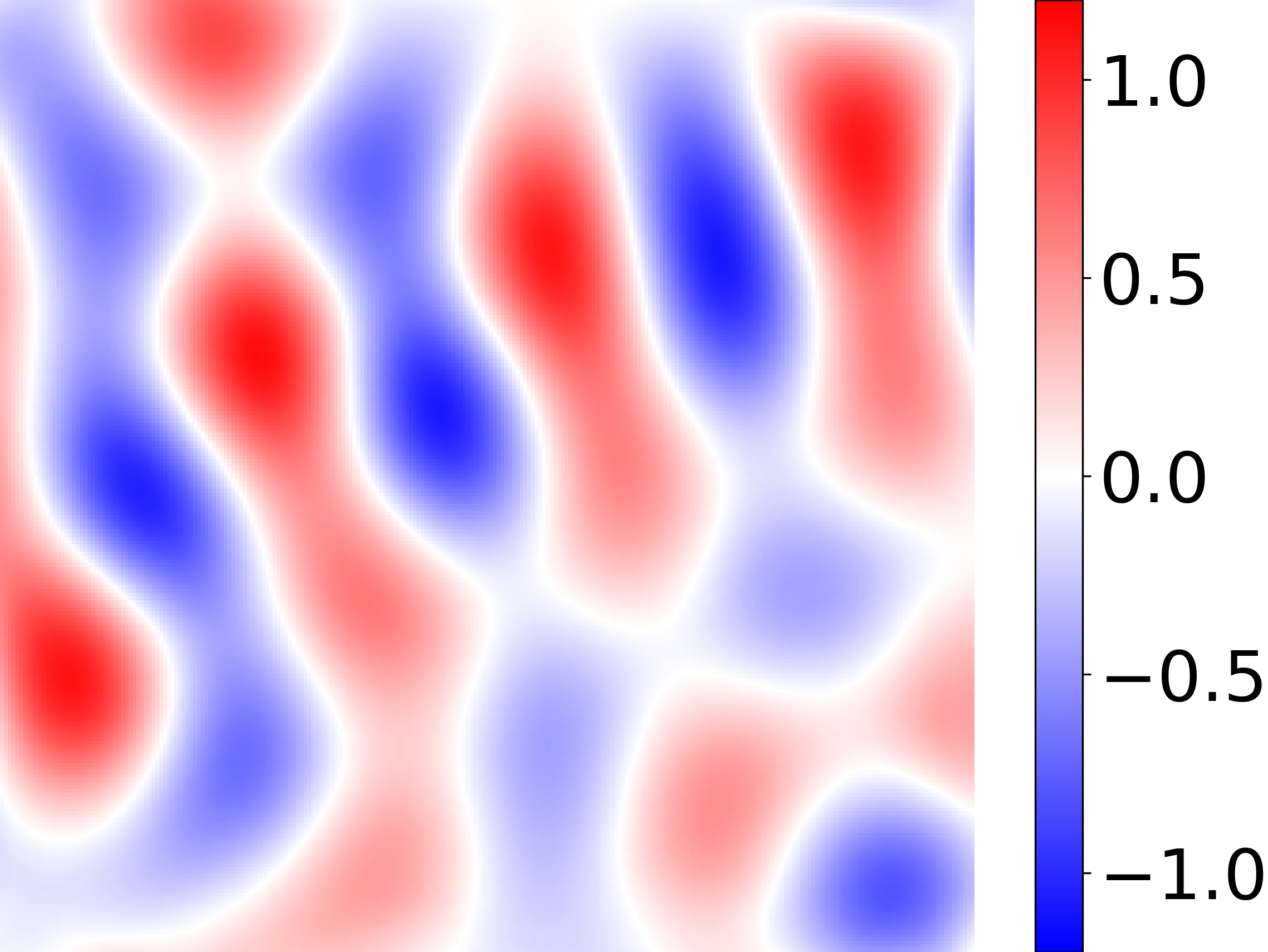}}
\hspace{0.04in}
\centering
\subfigure[UNO, error, 15.45\% \qquad]{
\includegraphics[width=0.3\textwidth]{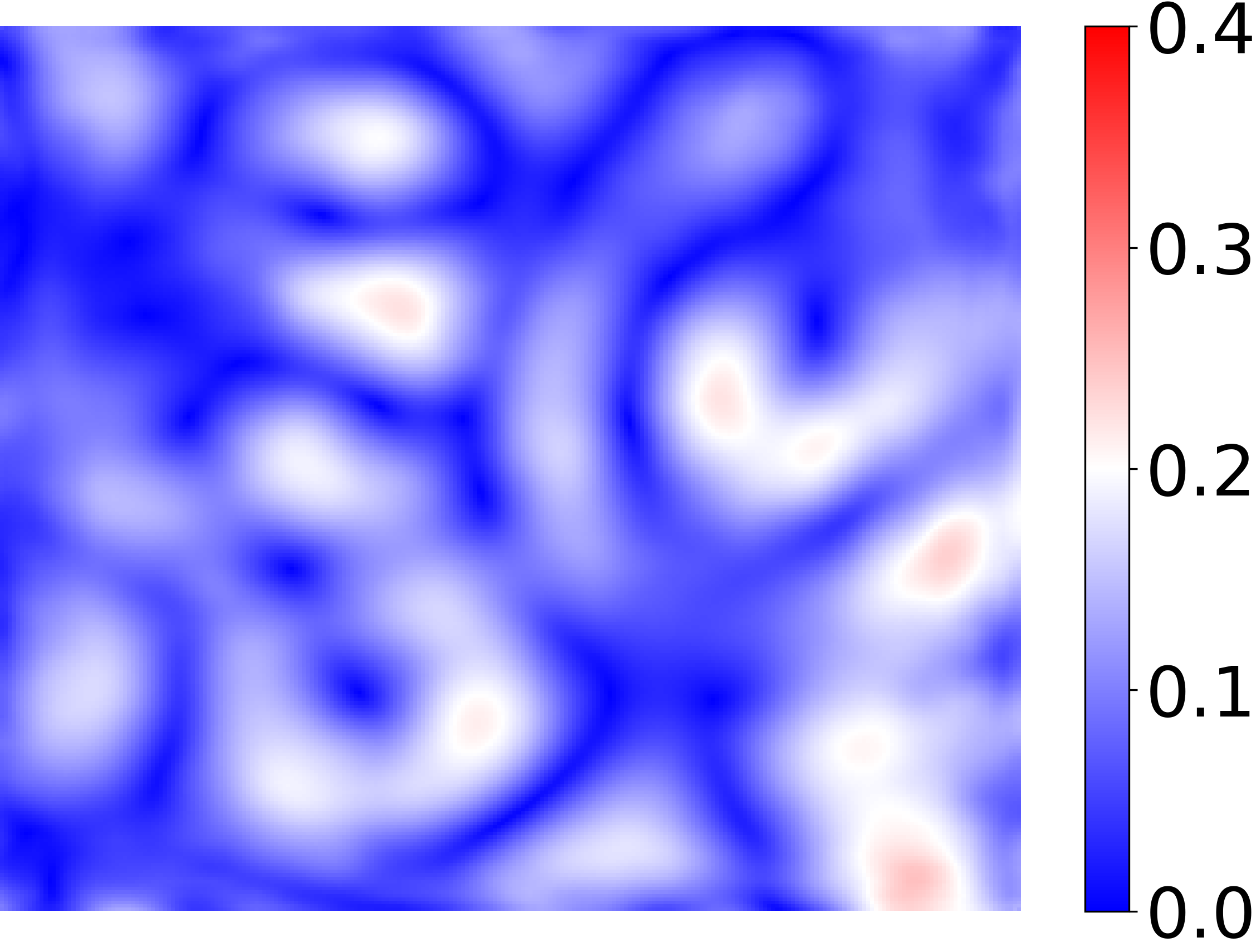}}

\subfigure[NS-FNO, real part \quad]{
\includegraphics[width=0.3\textwidth]{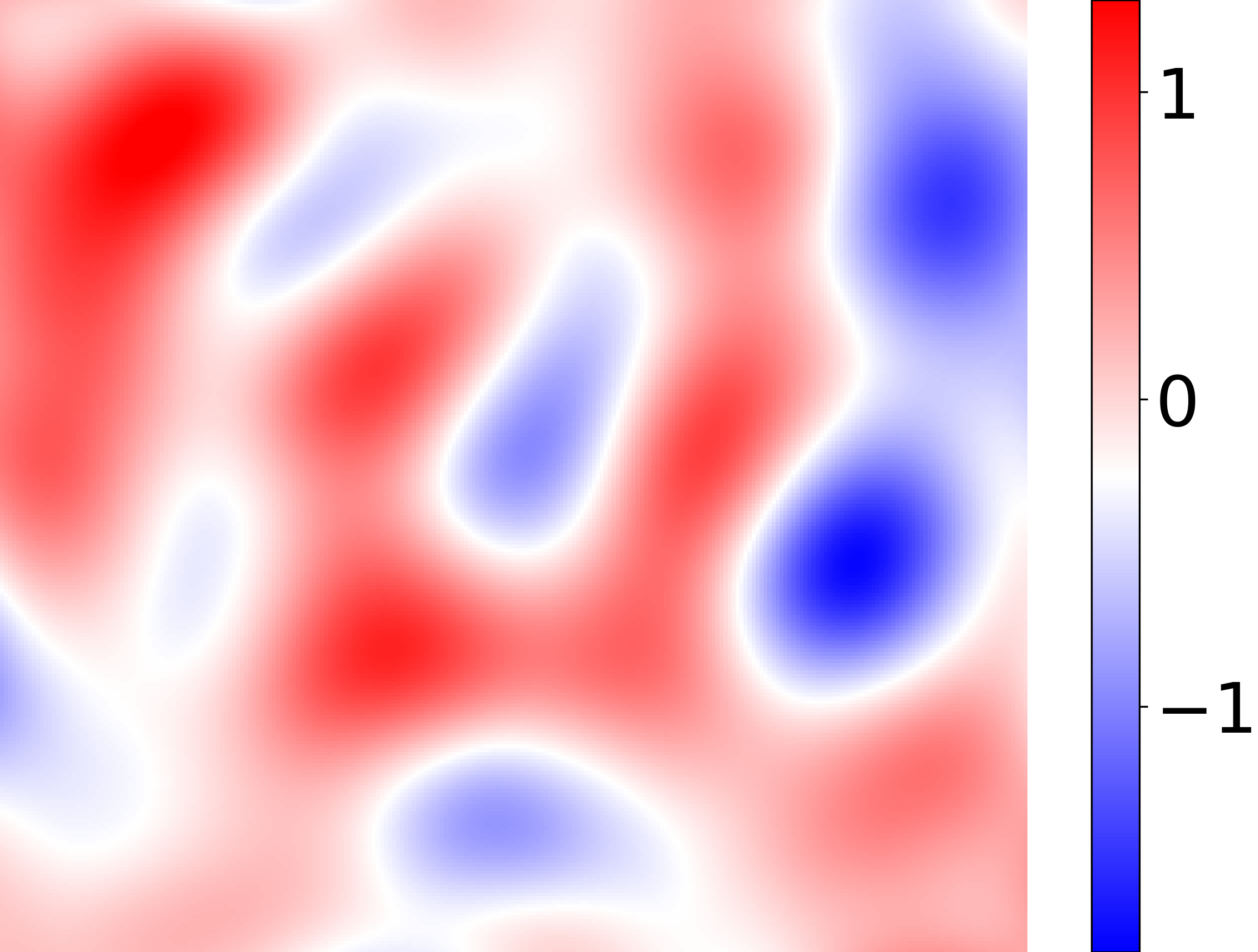}}
\hspace{0.04in}
\centering
\subfigure[NS-FNO, imaginary part\qquad]{
\includegraphics[width=0.3\textwidth]{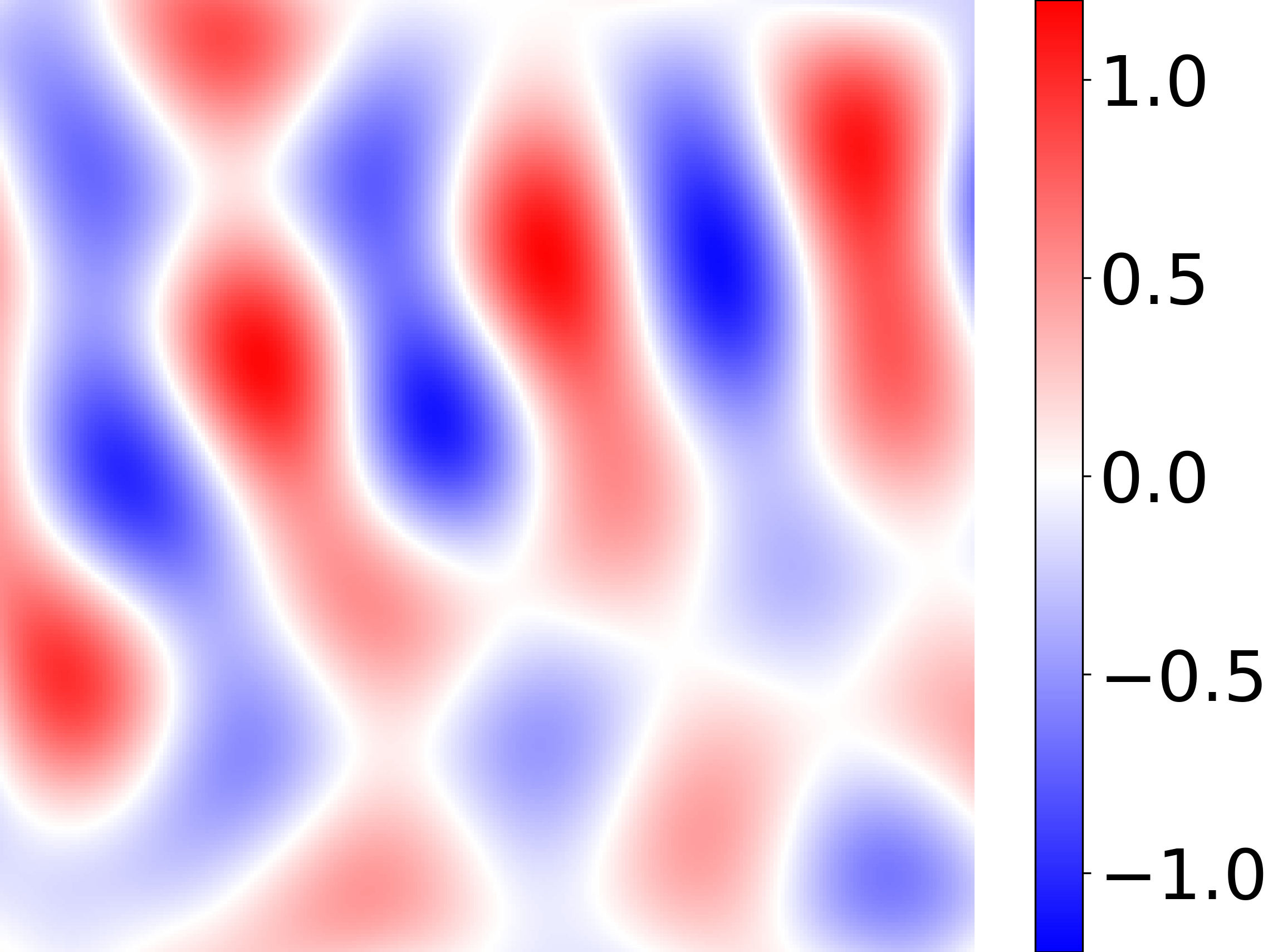}}
\hspace{0.04in}
\centering
\subfigure[NS-FNO, error, 5.11\% \qquad]{
\includegraphics[width=0.3\textwidth]{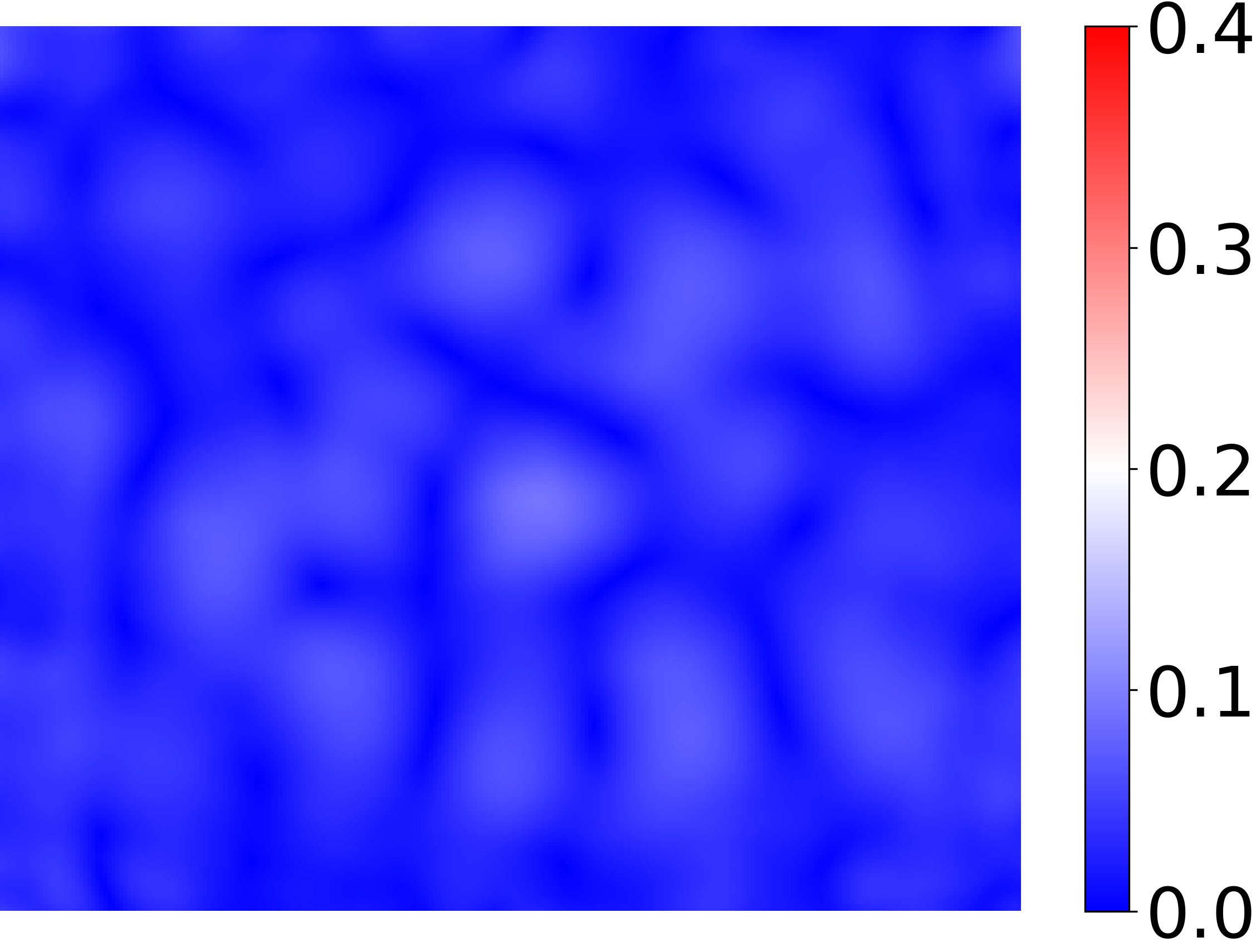}}

\subfigure[NS-UNO, real part \quad]{
\includegraphics[width=0.3\textwidth]{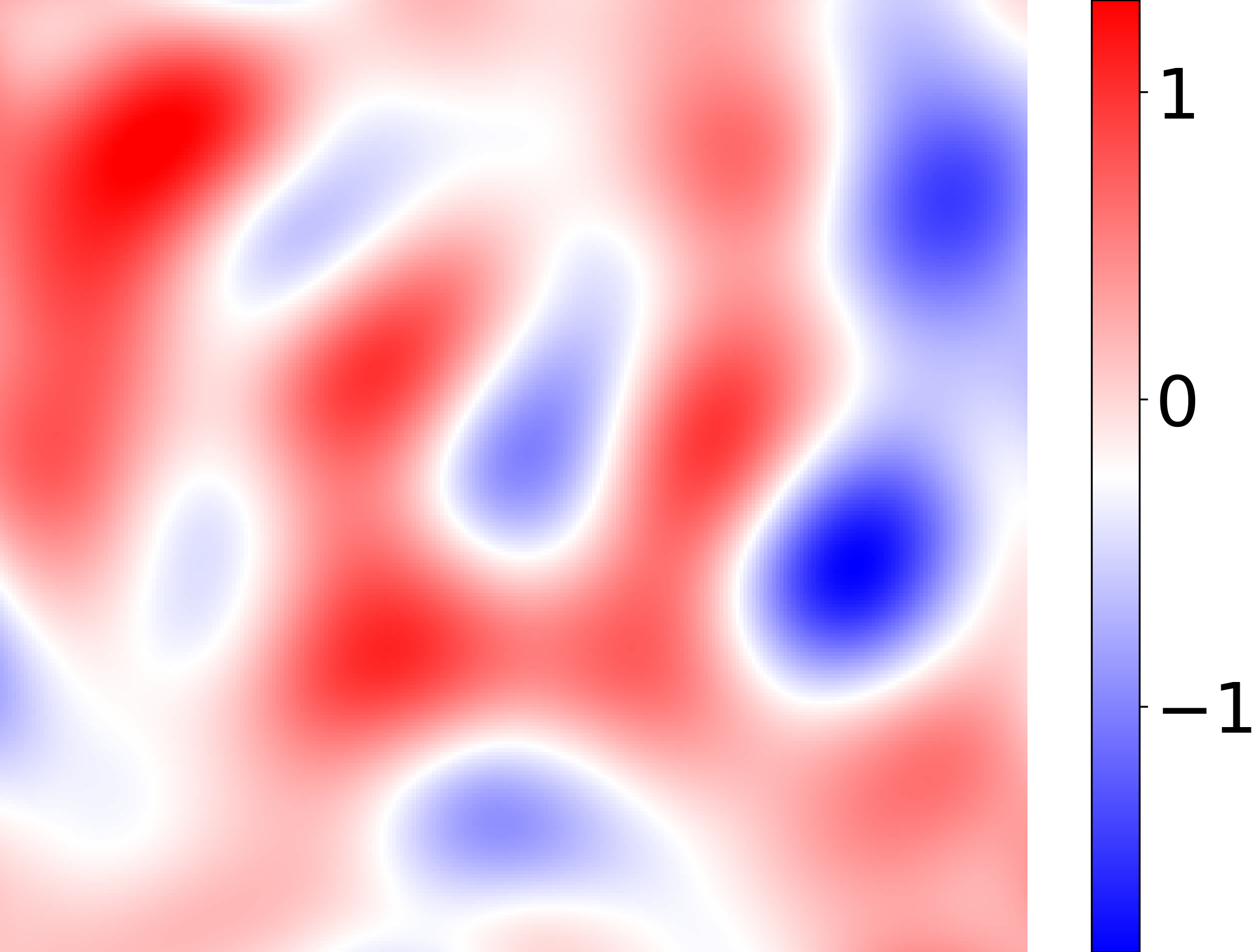}}
\hspace{0.04in}
\centering
\subfigure[NS-UNO, imaginary part\qquad]{
\includegraphics[width=0.3\textwidth]{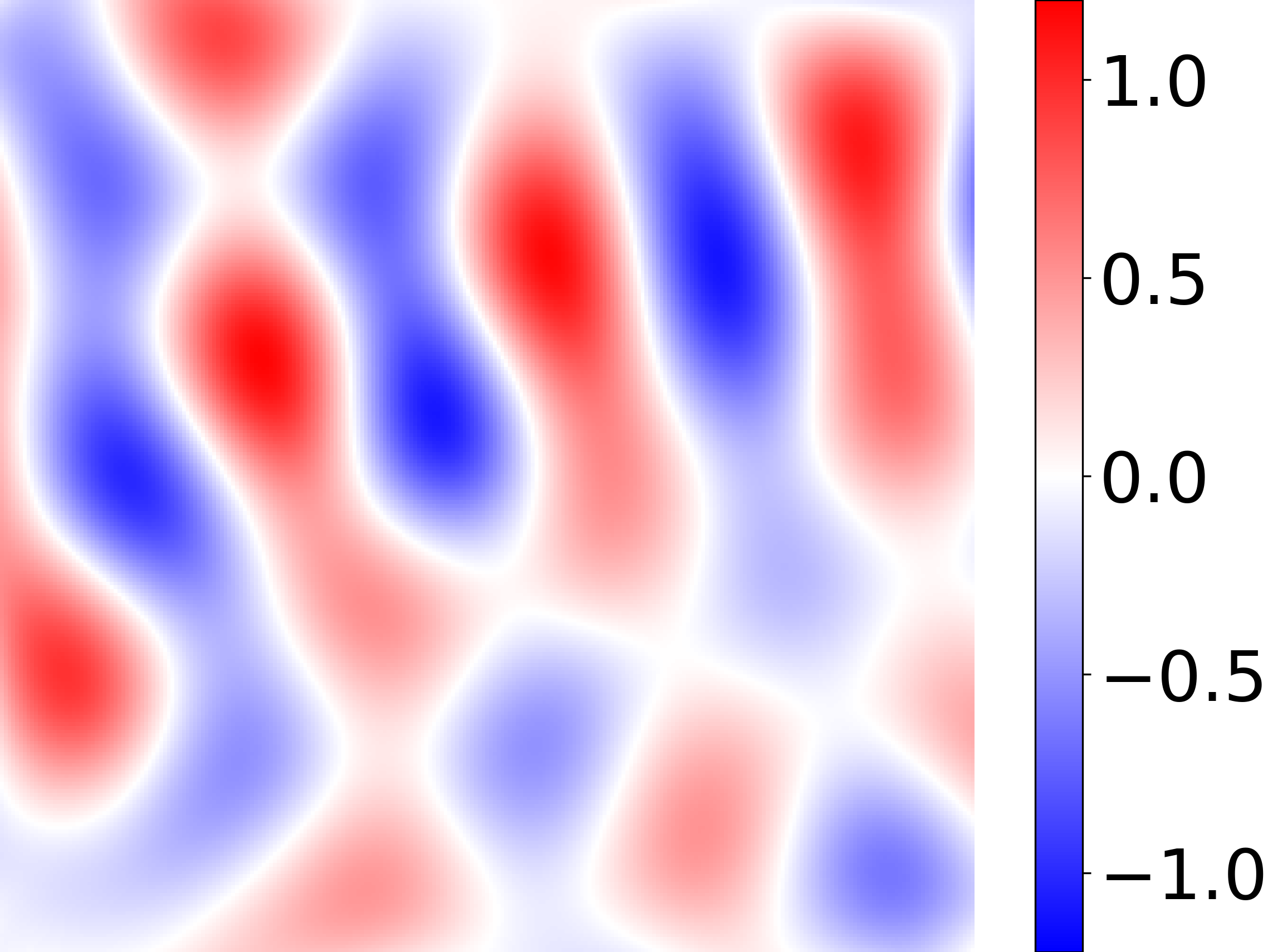}}
\hspace{0.04in}
\centering
\subfigure[NS-UNO, error, 1.94\% \qquad]{
\includegraphics[width=0.3\textwidth]{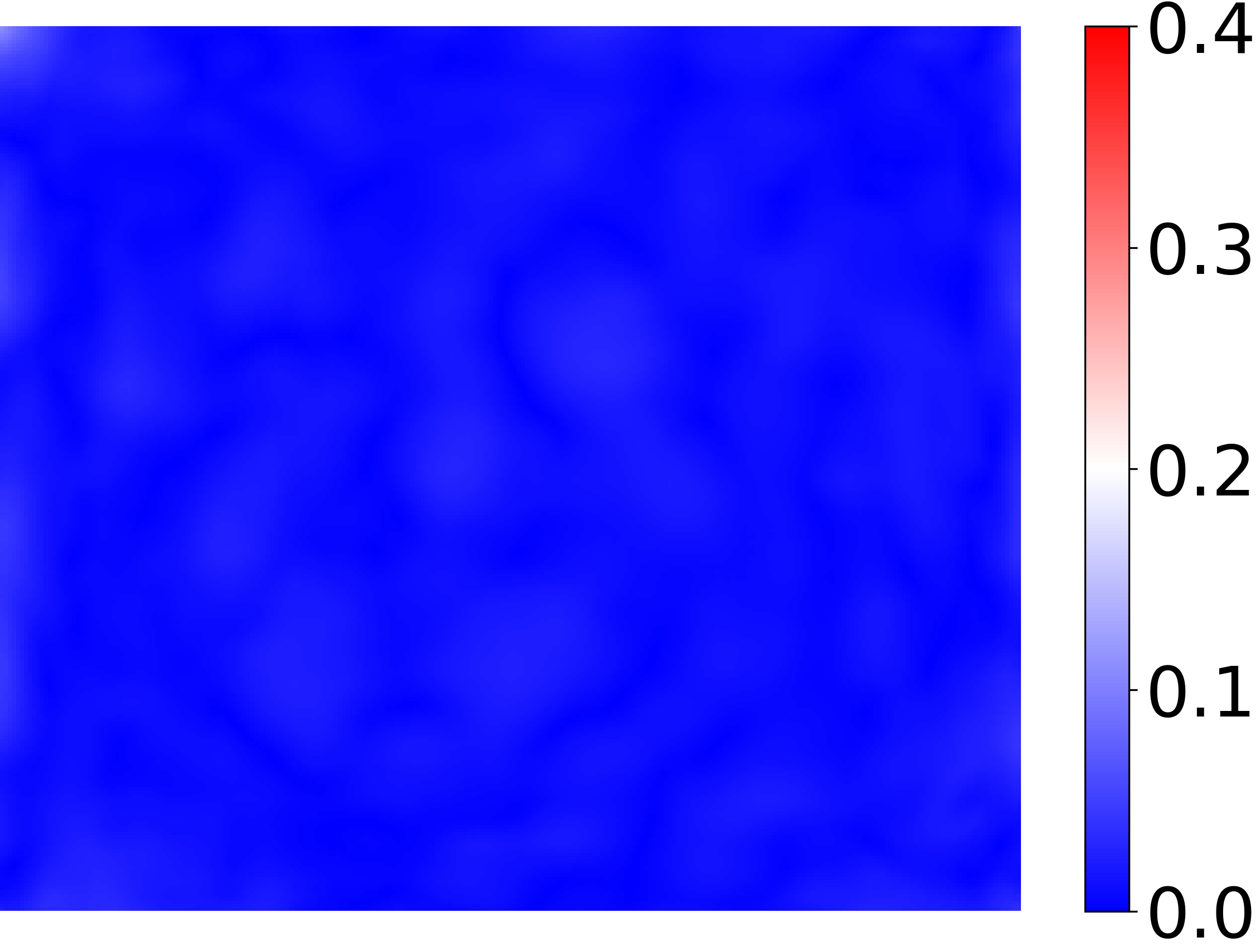}}

\caption{Example of exact solution, numerical solutions and absolute error for dataset with T-shaped $q$ and GRF $f$ when $k=20$}
\label{eg1}
\end{figure}

\begin{figure}
\centering
\subfigure[Ground truth, real part \qquad]{
\includegraphics[width=0.3\textwidth]{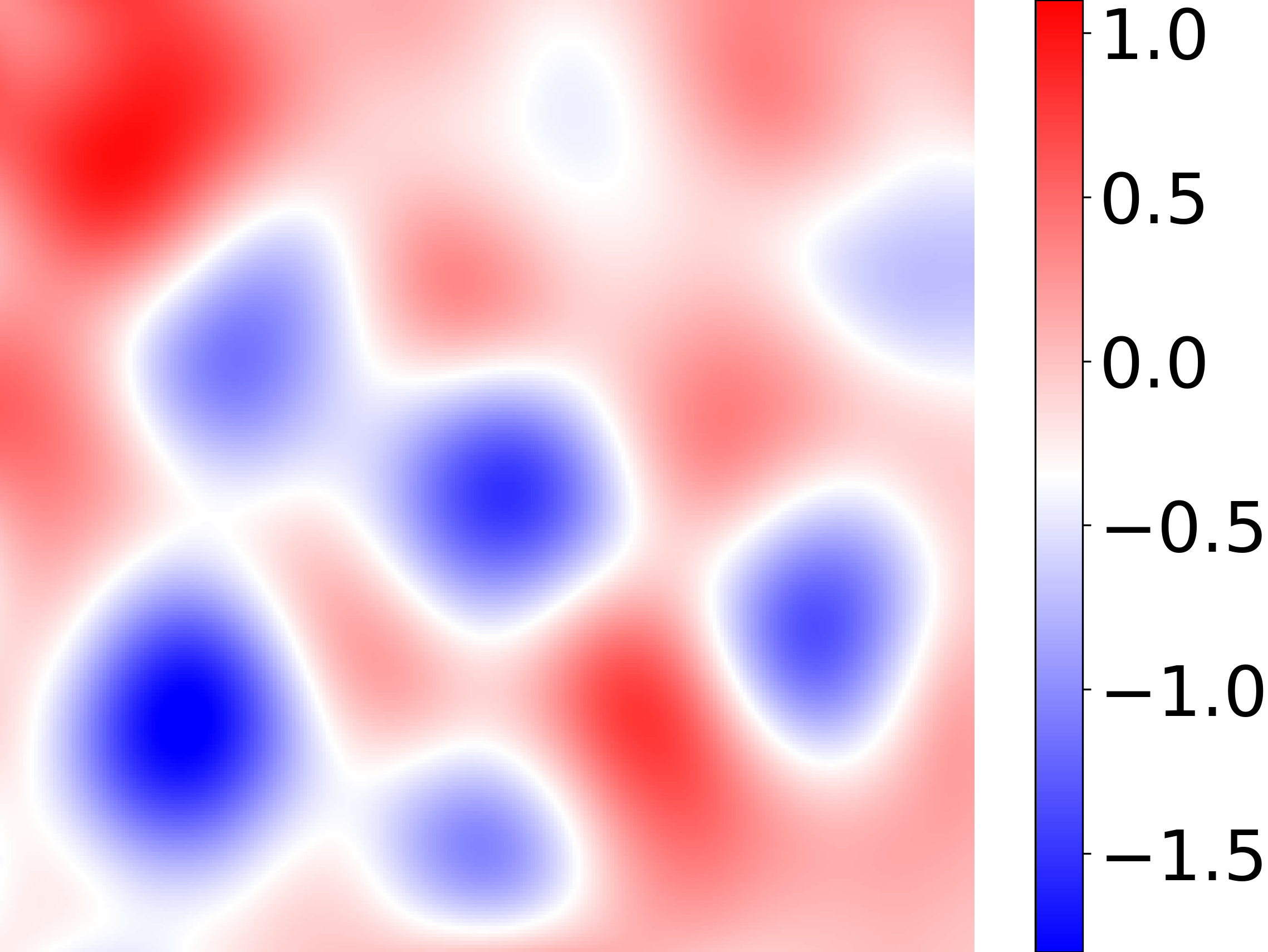}}
\hspace{0.04in}
\centering
\subfigure[Ground truth, imaginary part \qquad]{
\includegraphics[width=0.3\textwidth]{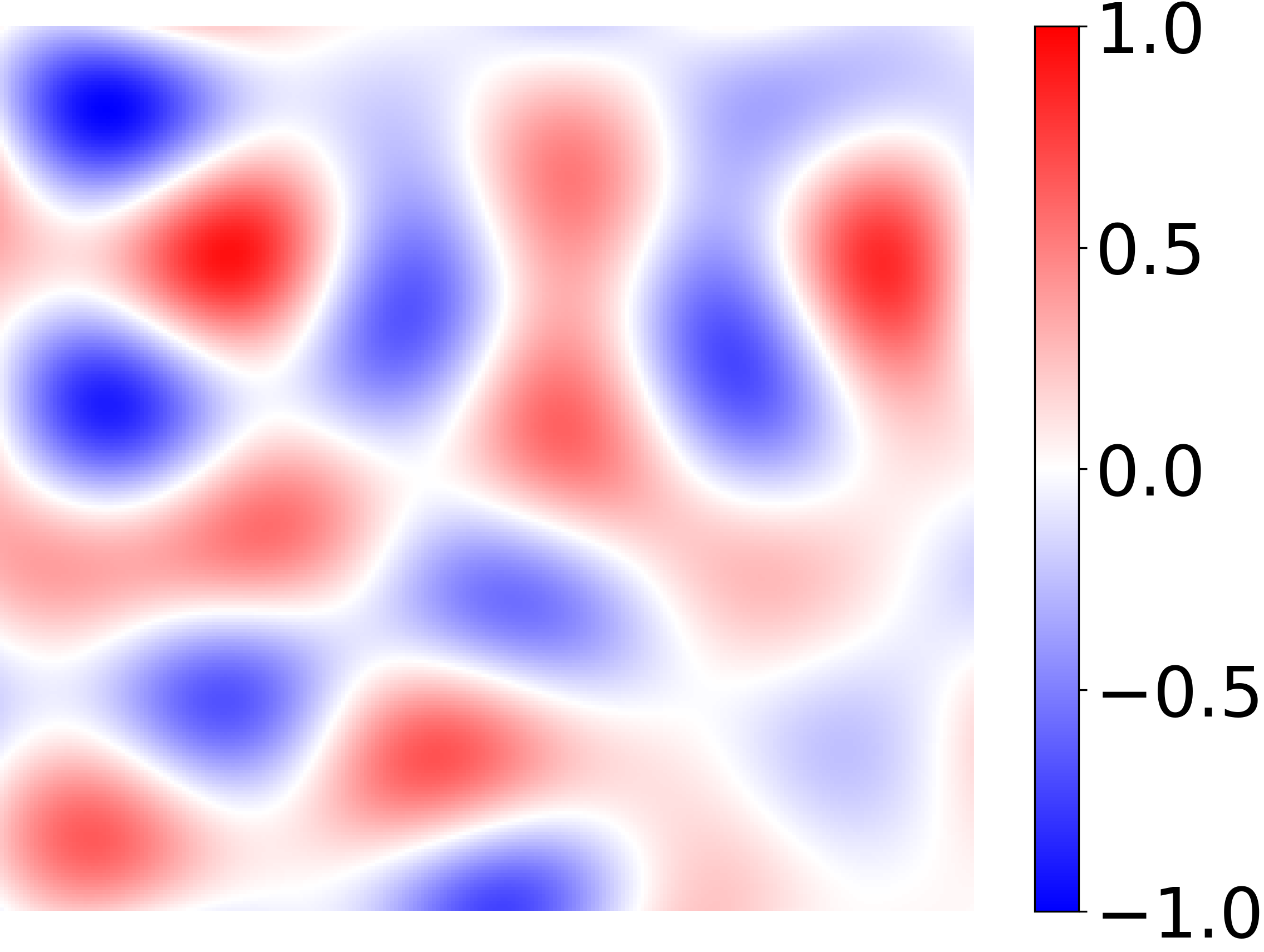}}
\hspace{0.04in}
\centering
\subfigure[$q$ (top left) and $f$\qquad]{
\includegraphics[width=0.3\textwidth]{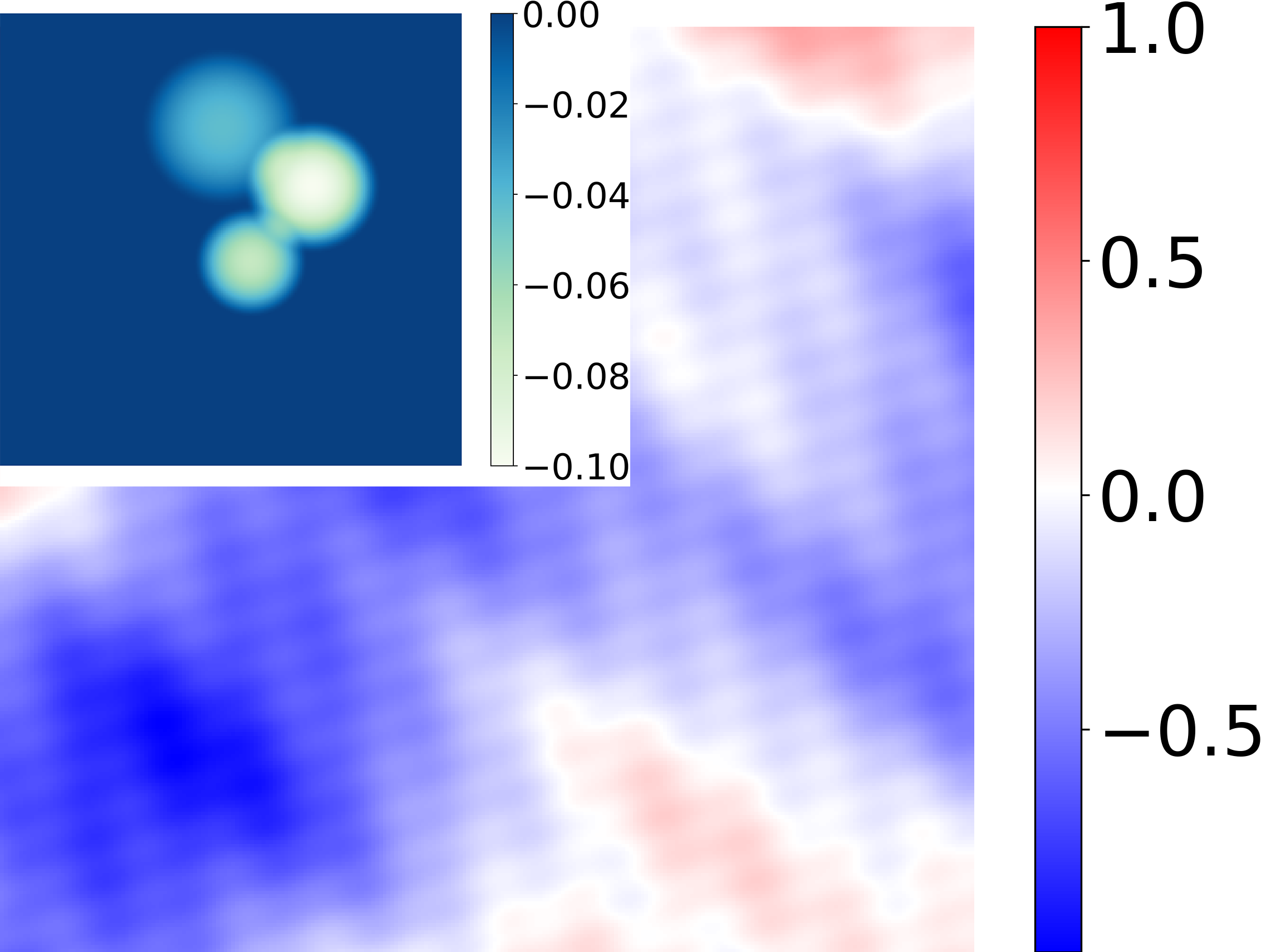}}

\subfigure[FNO, real part \quad]{
\includegraphics[width=0.3\textwidth]{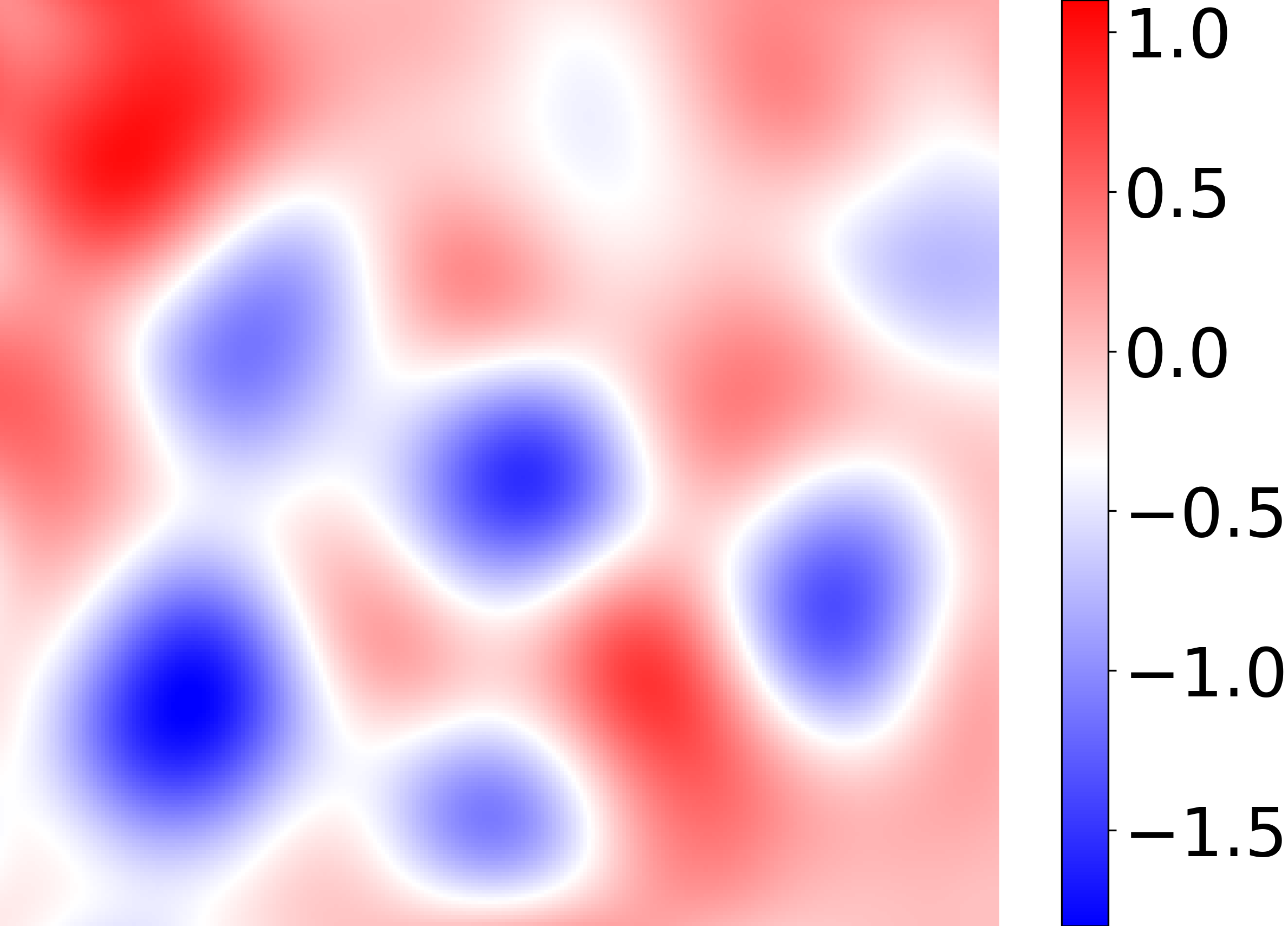}}
\hspace{0.04in}
\centering
\subfigure[FNO, imaginary part\qquad]{
\includegraphics[width=0.3\textwidth]{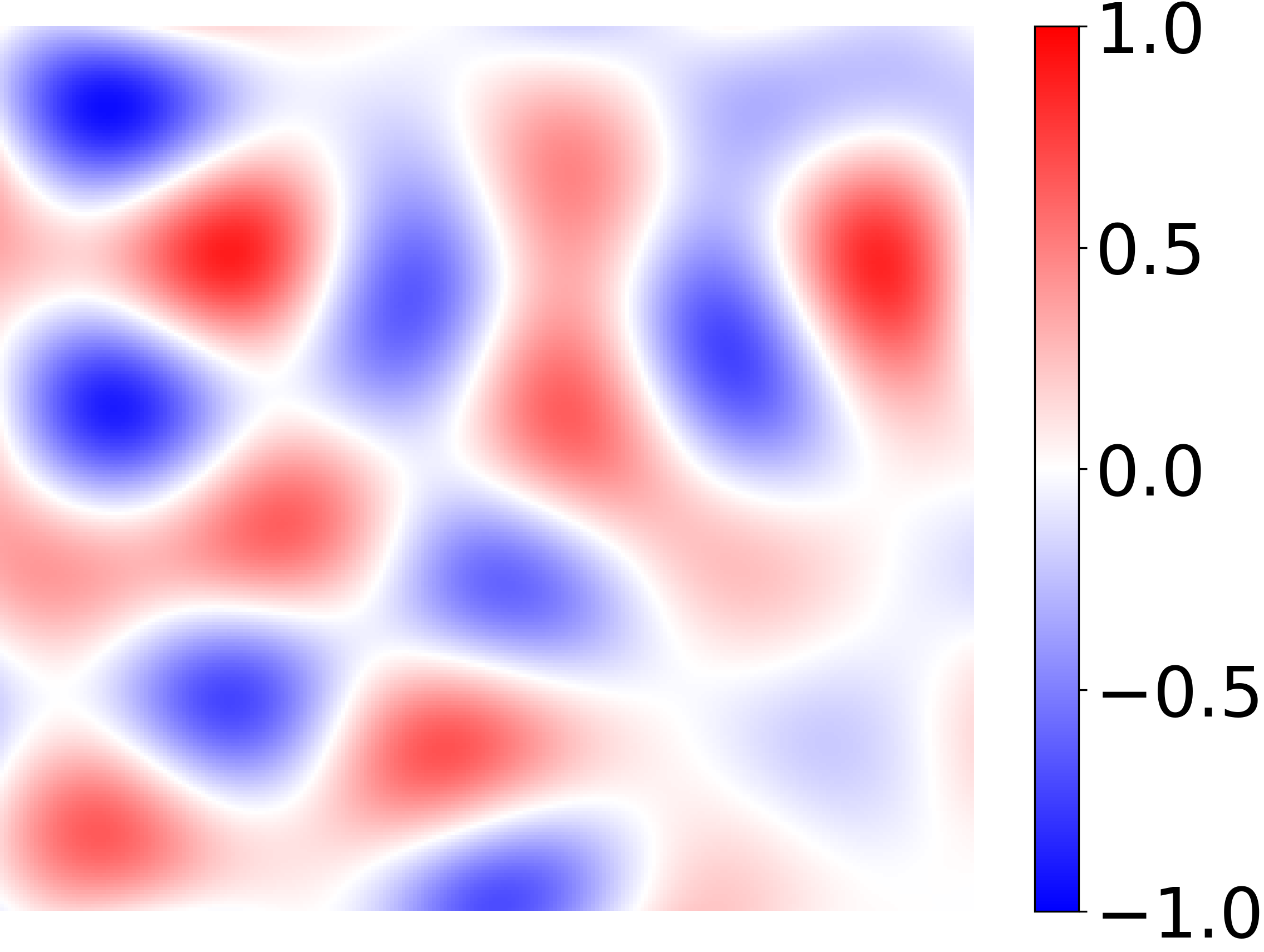}}
\hspace{0.04in}
\centering
\subfigure[FNO, error, 9.57\% \qquad]{
\includegraphics[width=0.3\textwidth]{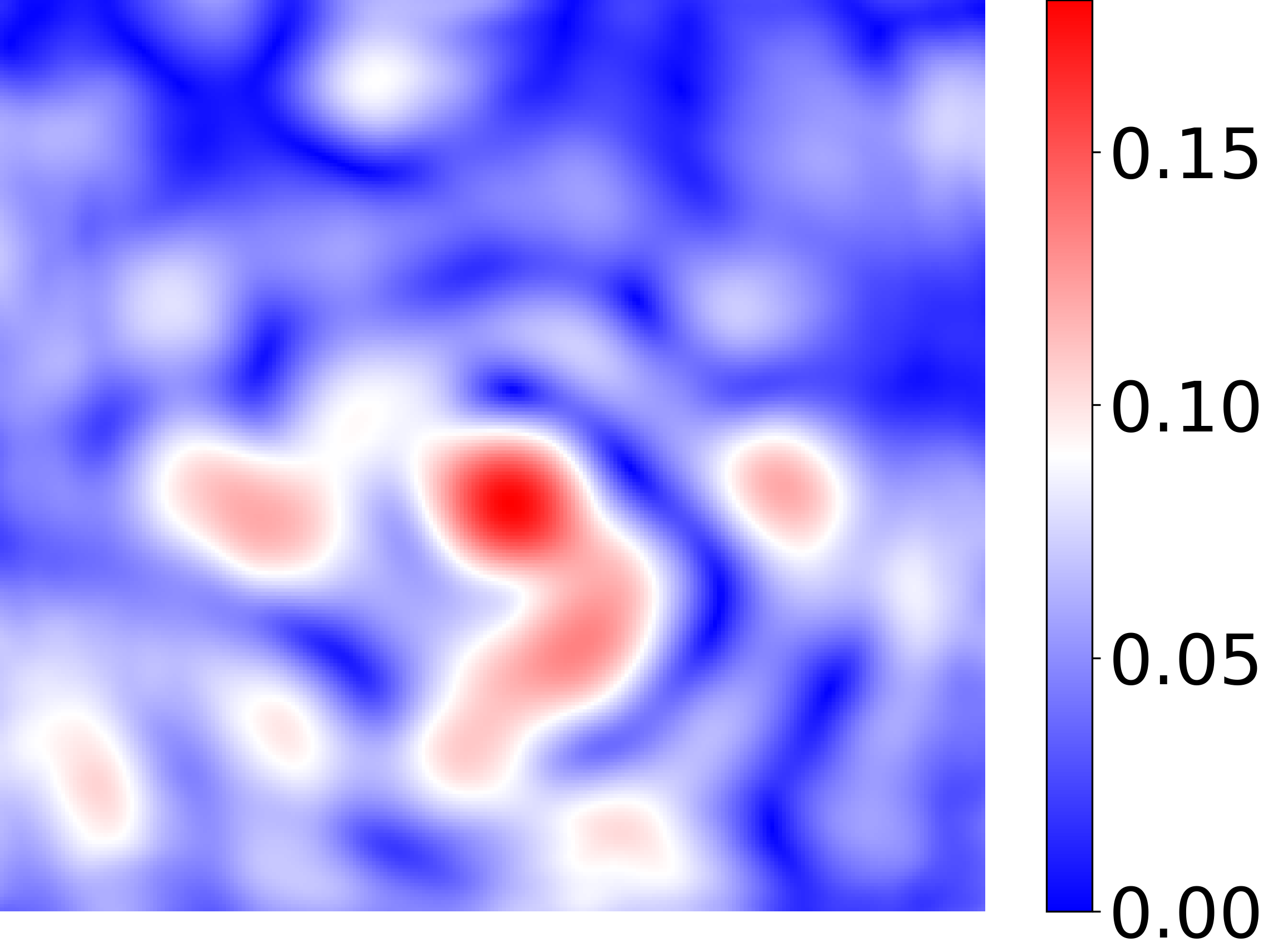}}

\subfigure[UNO, real part \quad]{
\includegraphics[width=0.3\textwidth]{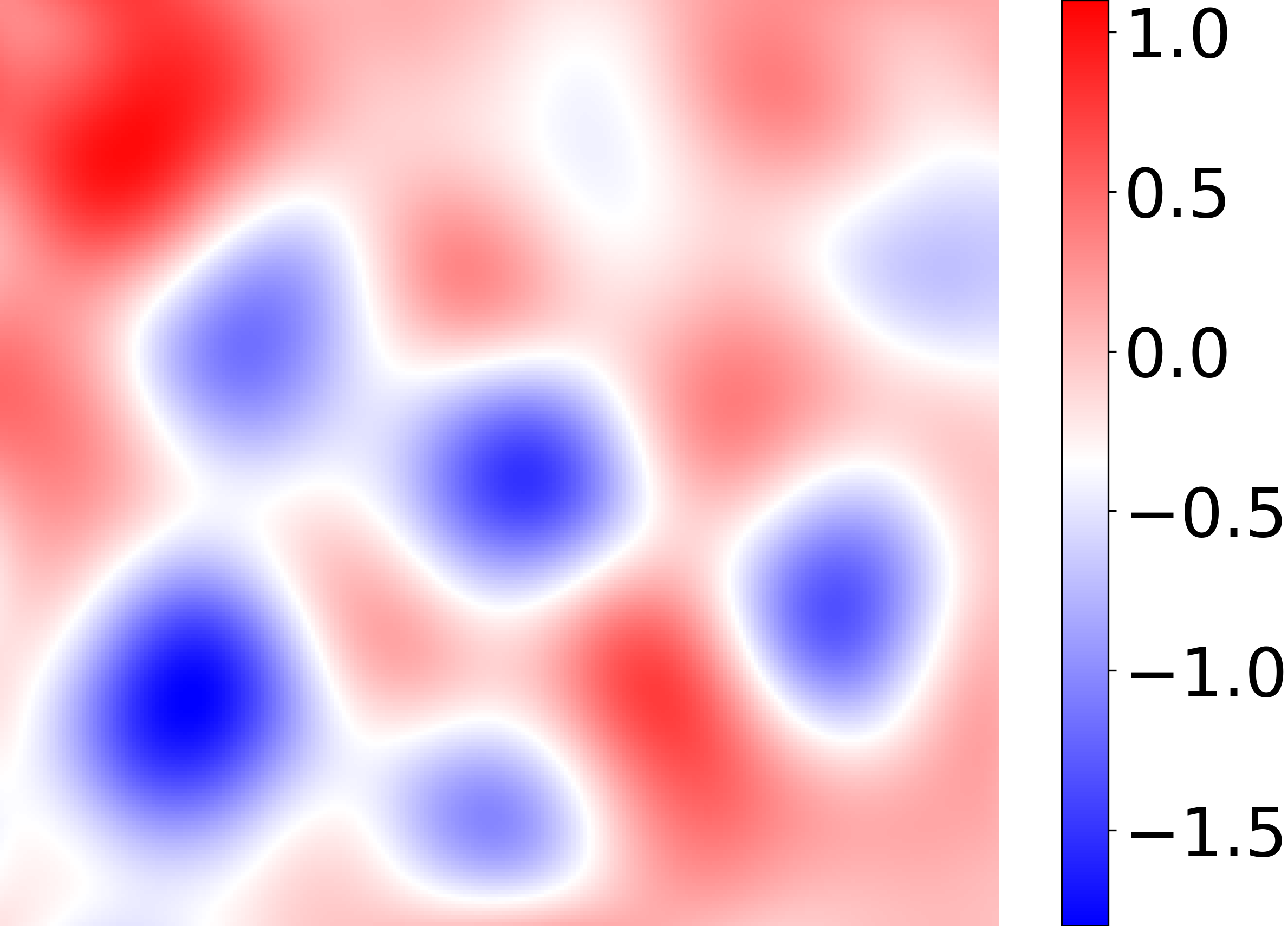}}
\hspace{0.04in}
\centering
\subfigure[UNO, imaginary part\qquad]{
\includegraphics[width=0.3\textwidth]{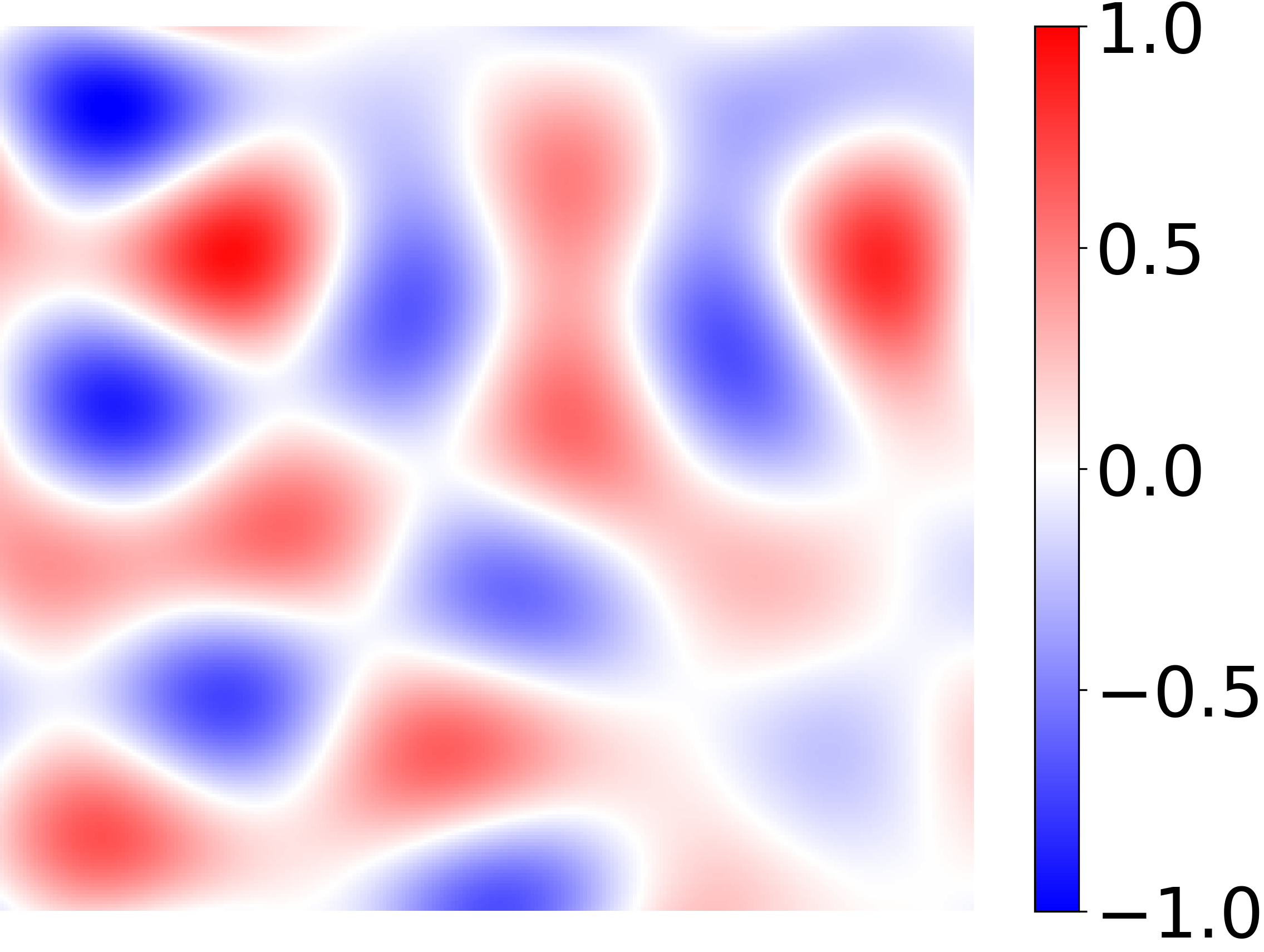}}
\hspace{0.04in}
\centering
\subfigure[UNO, error, 7.70\% \qquad]{
\includegraphics[width=0.3\textwidth]{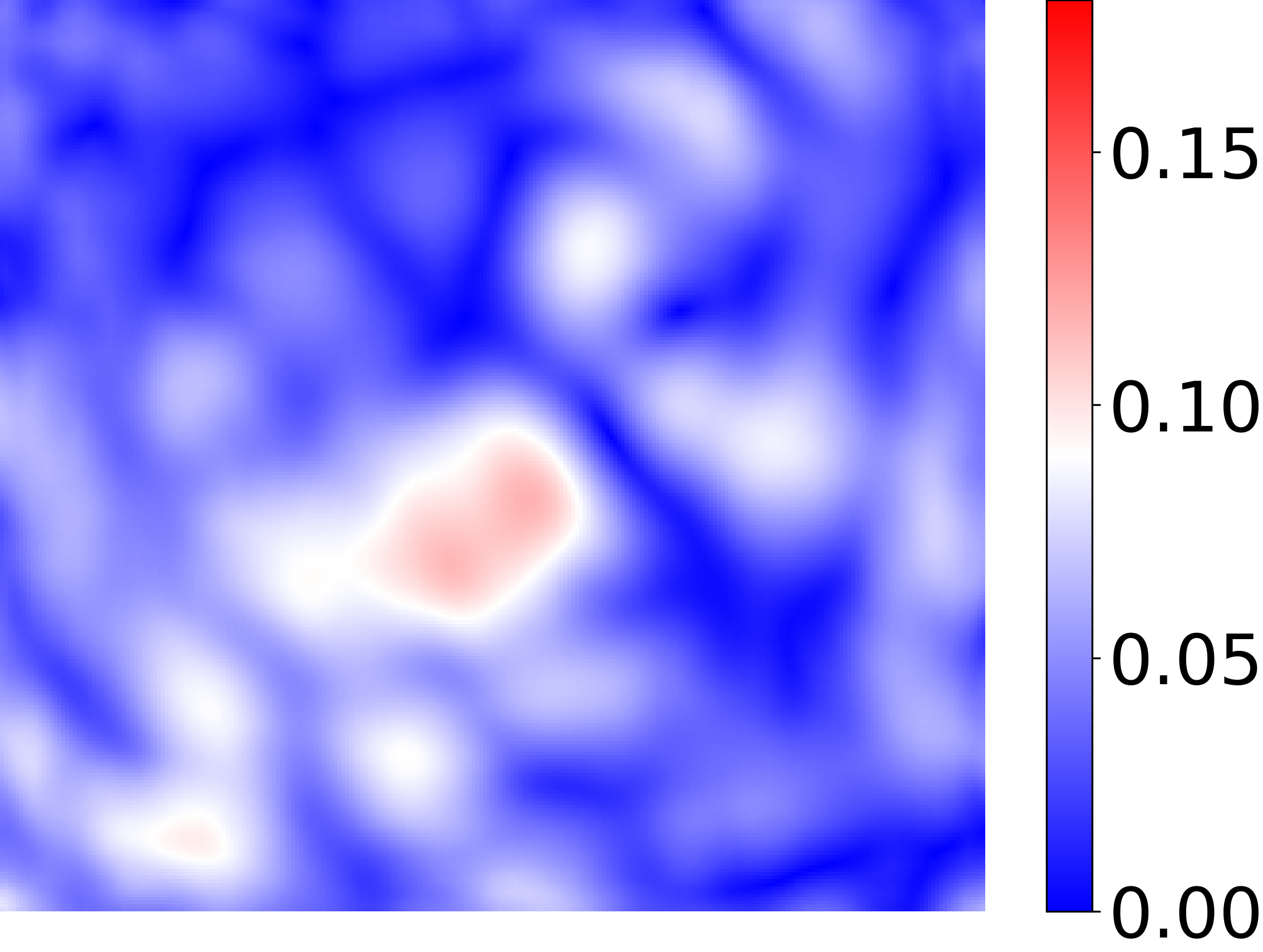}}

\subfigure[NS-FNO, real part \quad]{
\includegraphics[width=0.3\textwidth]{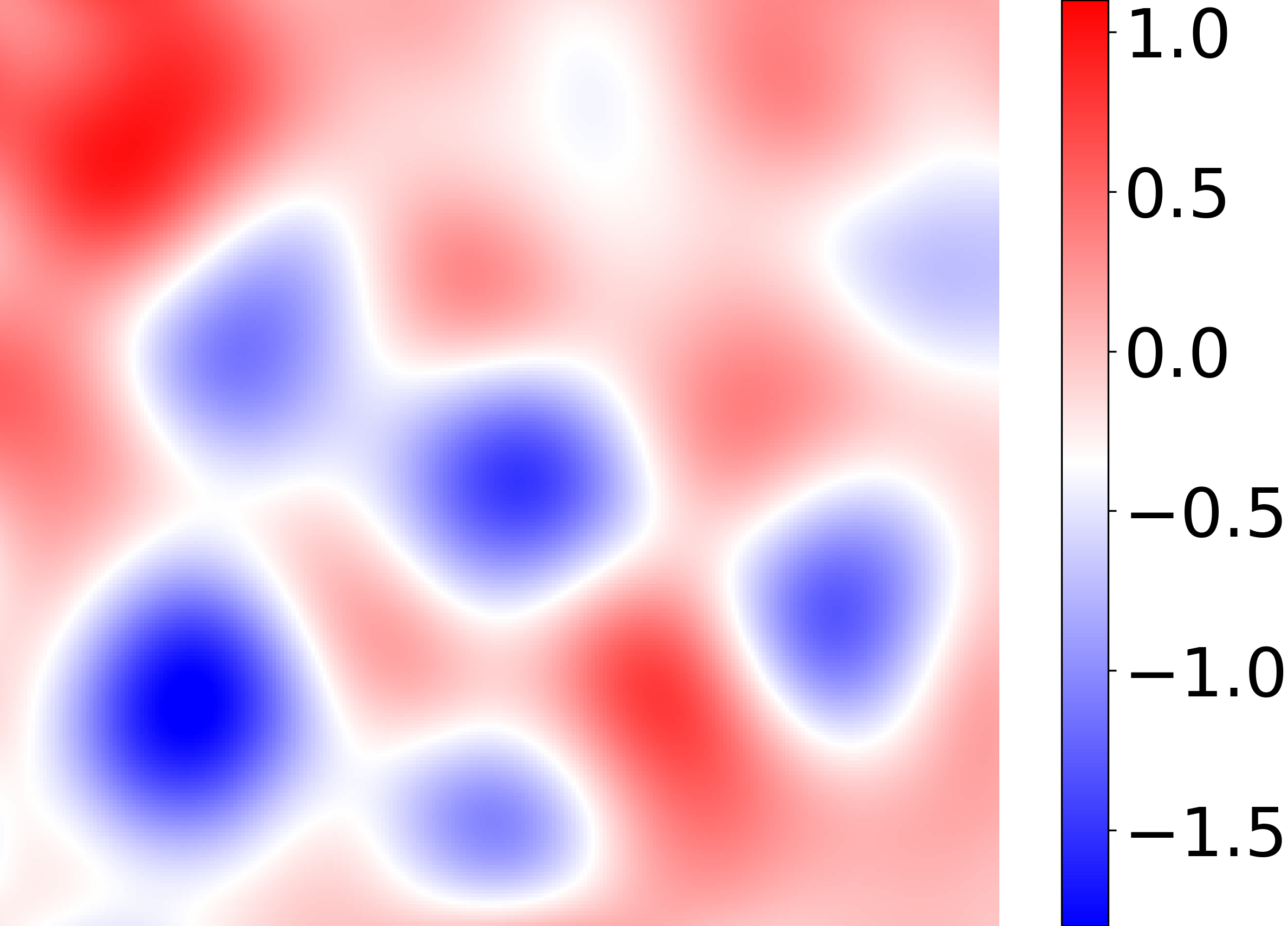}}
\hspace{0.04in}
\centering
\subfigure[NS-FNO, imaginary part\qquad]{
\includegraphics[width=0.3\textwidth]{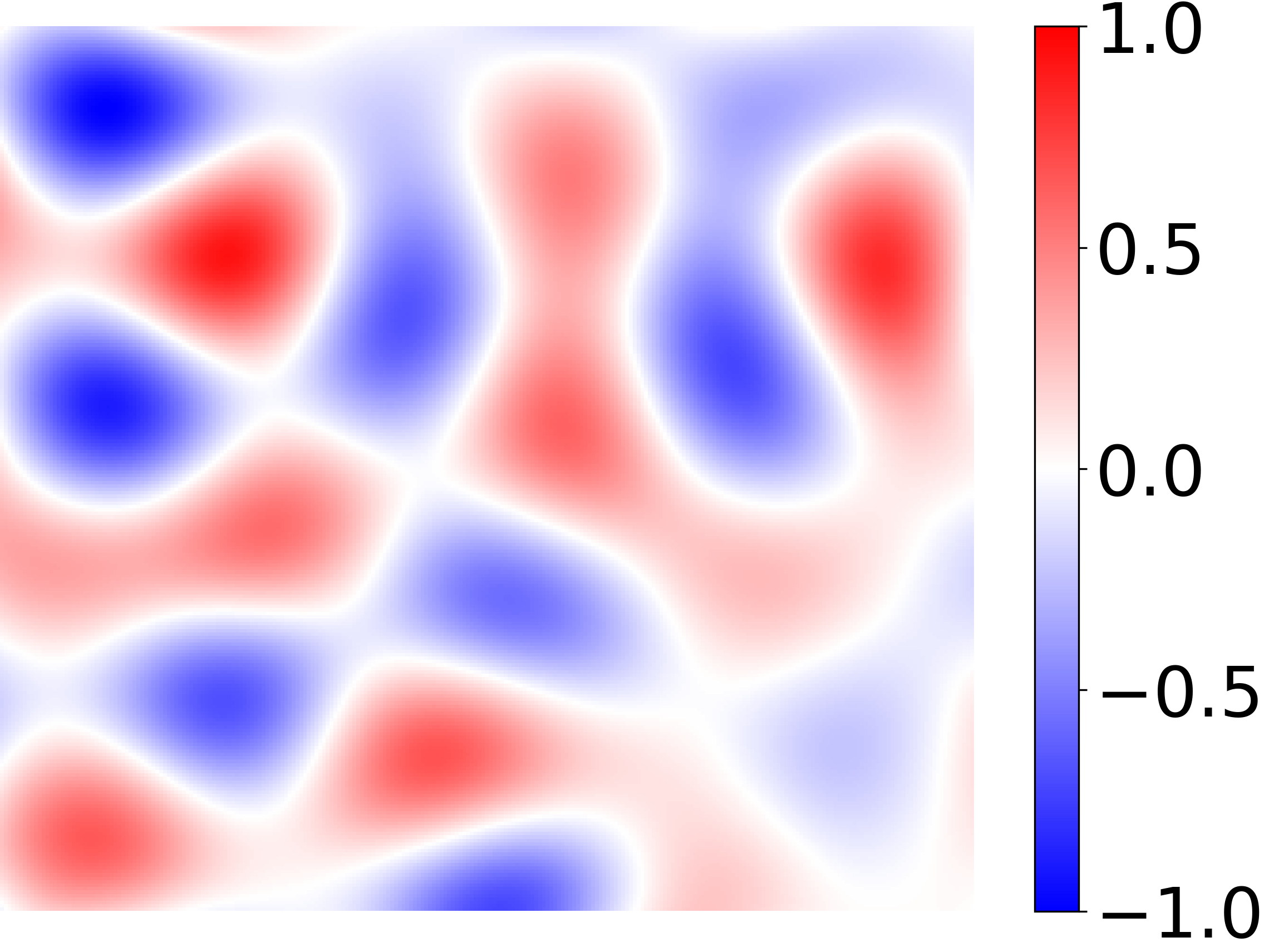}}
\hspace{0.04in}
\centering
\subfigure[NS-FNO, error, 2.70\% \qquad]{
\includegraphics[width=0.3\textwidth]{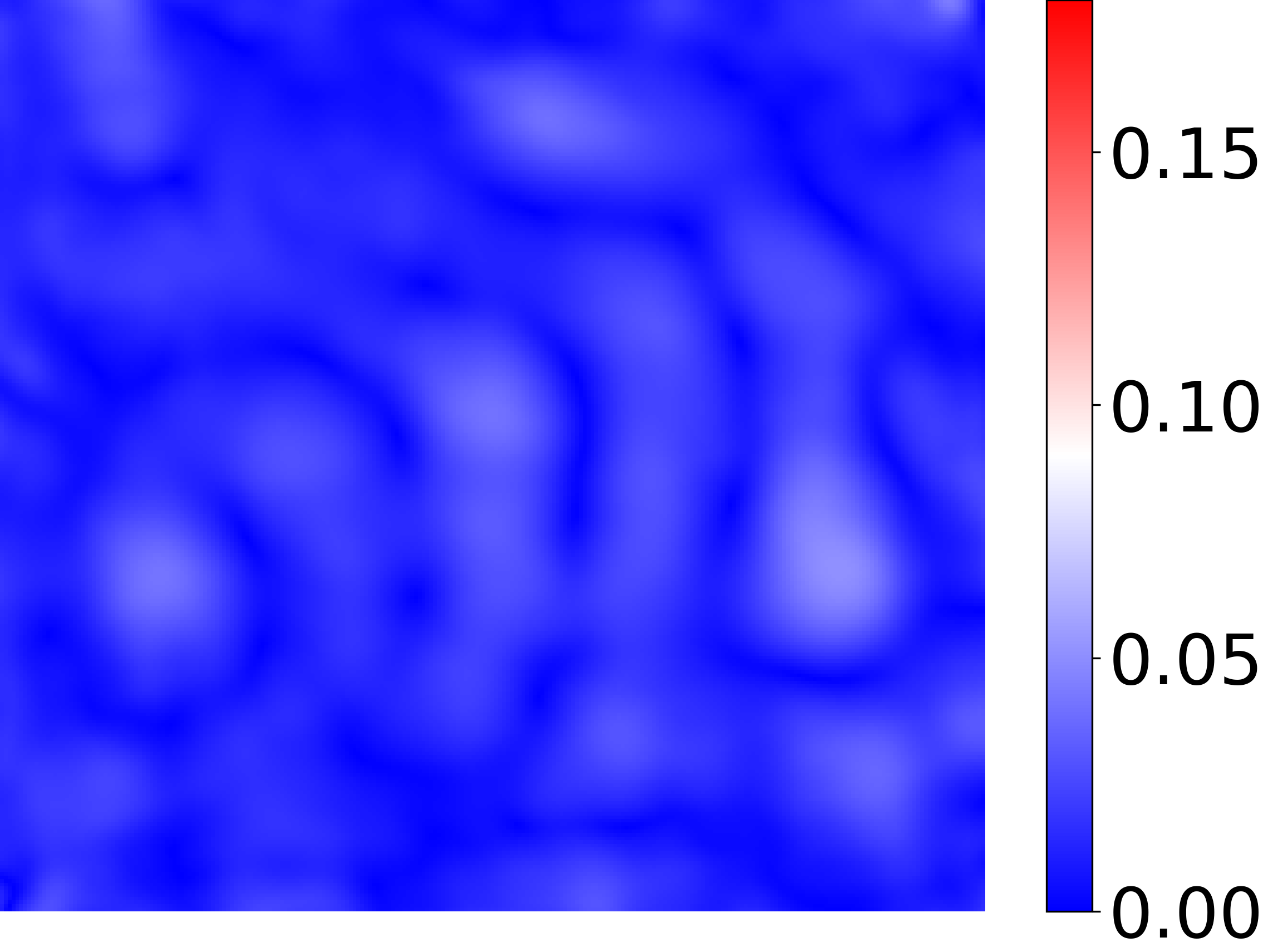}}

\subfigure[NS-UNO, real part \quad]{
\includegraphics[width=0.3\textwidth]{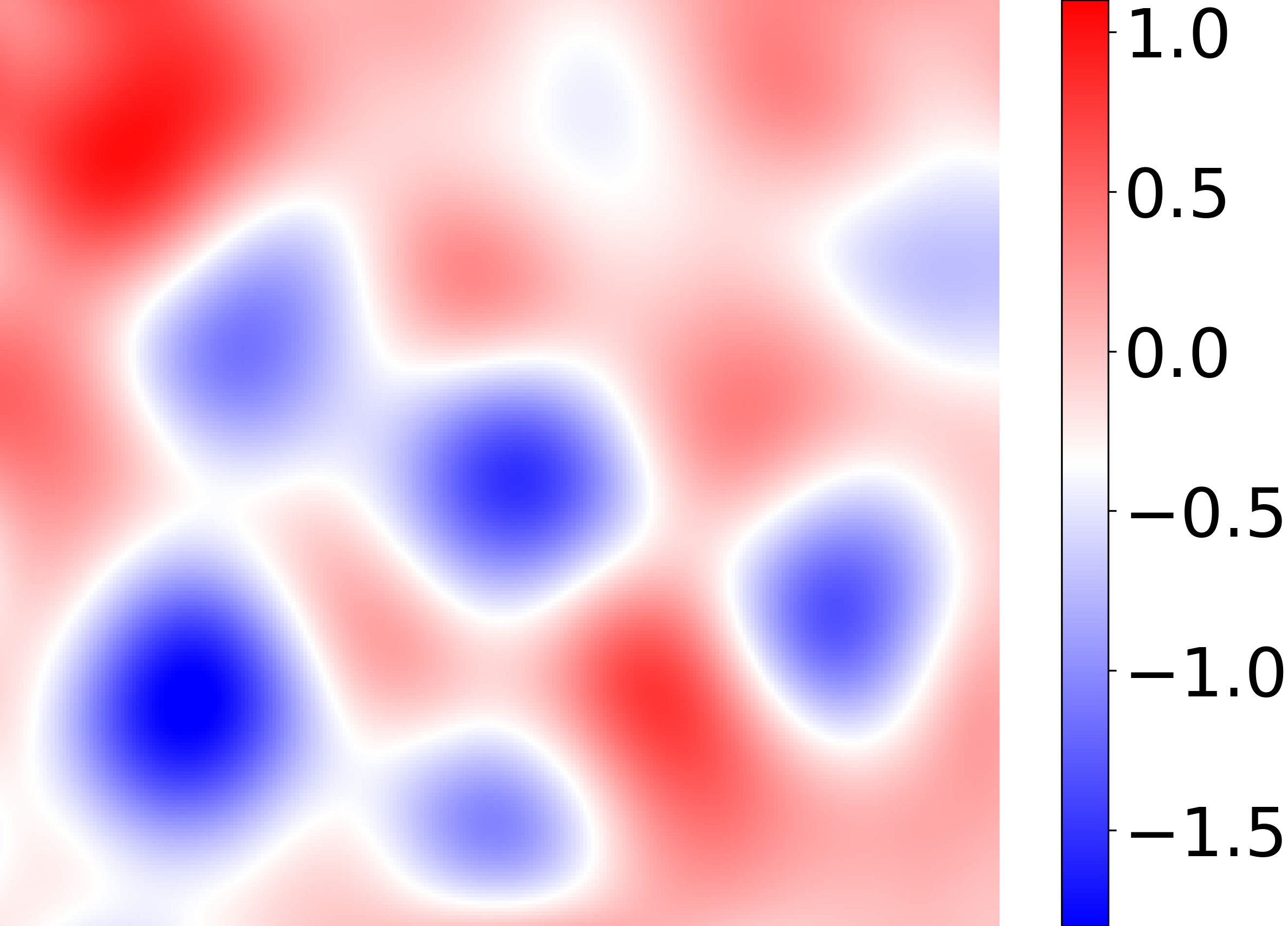}}
\hspace{0.04in}
\centering
\subfigure[NS-UNO, imaginary part\qquad]{
\includegraphics[width=0.3\textwidth]{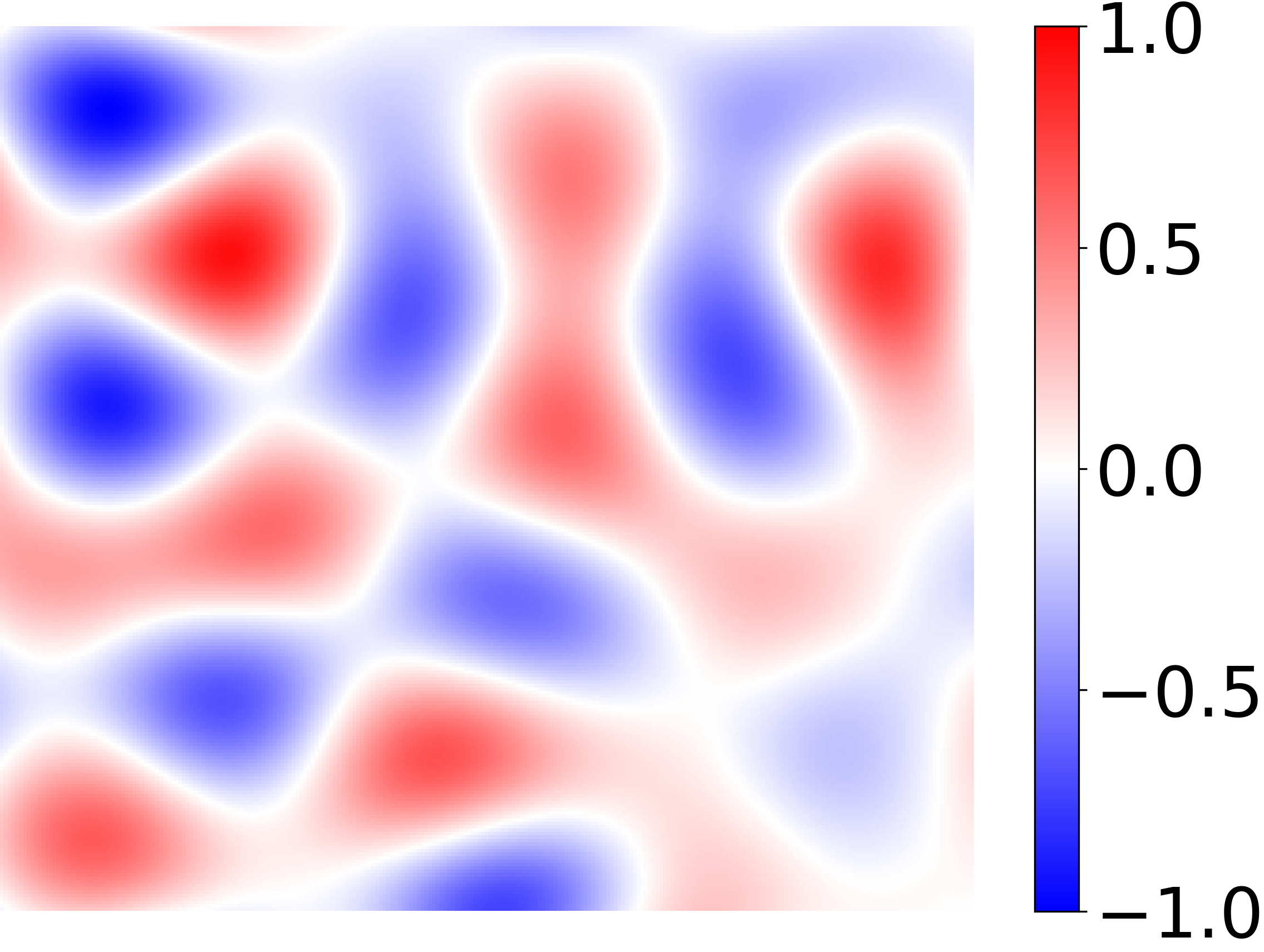}}
\hspace{0.04in}
\centering
\subfigure[NS-UNO, error, 1.32\% \qquad]{
\includegraphics[width=0.3\textwidth]{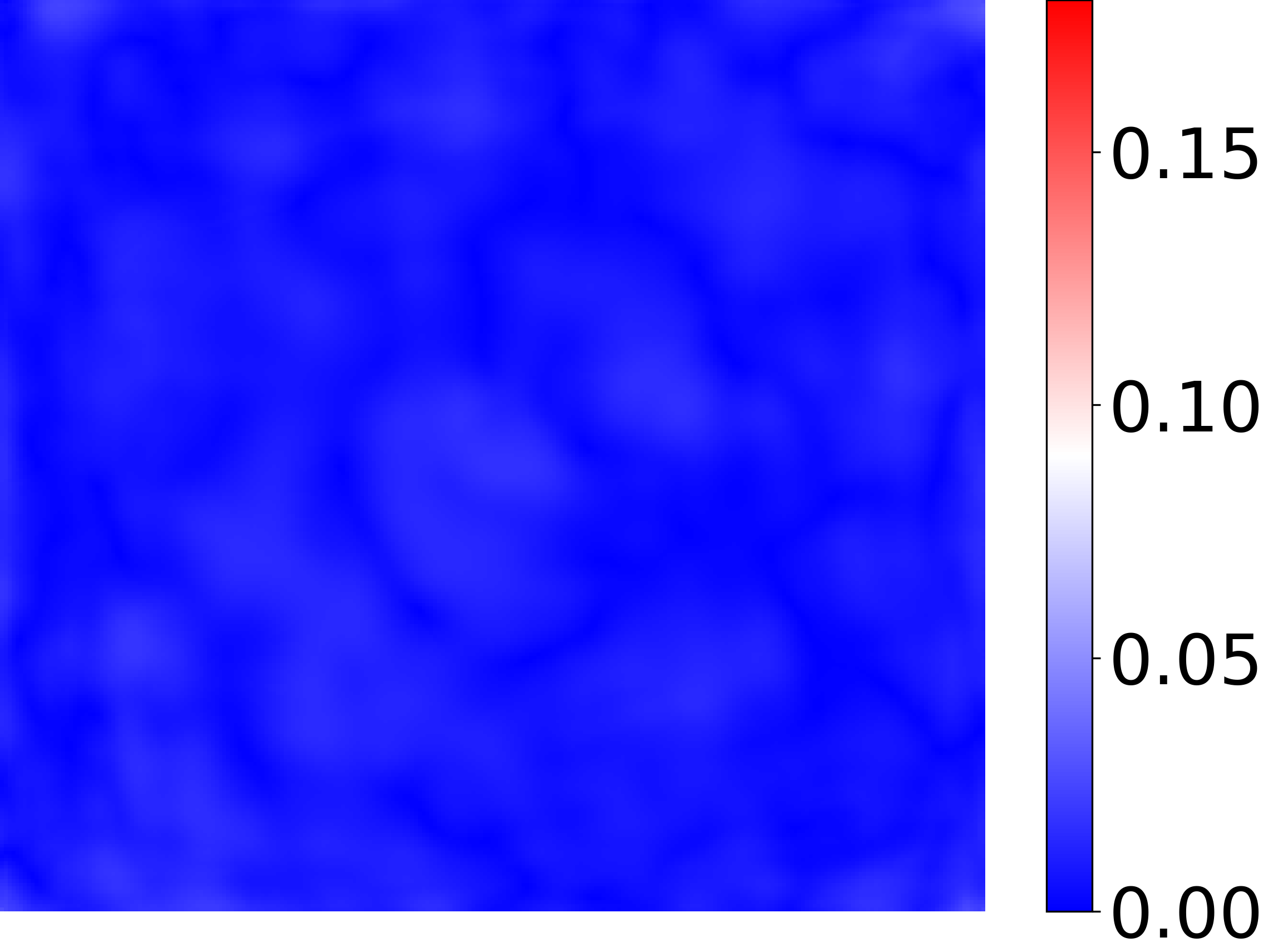}}

\caption{Example of exact solution, numerical solutions and absolute error for dataset with smooth circle $q$ and waves $f$ when $k=20$}
\label{eg2}
\end{figure}

In Fig. \ref{eg1} and Fig. \ref{eg2}, We showcase the exact solution and the numerical solution obtained by the four benchmark models for the dataset with T-shaped $q$ and $f$ generated by GRF as well as the dataset with smooth circles $q$ and $f$ generated by the wave distribution, respectively. The first row presents the real and imaginary parts of the exact solution and the corresponding parameters $(q, f)$, respectively. The real and imaginary parts and the absolute error of FNO, UNO, NS-FNO and NS-UNO are listed in the second to the fifth row in Fig. \ref{eg1} and Fig. \ref{eg2}, respectively. It can be seen that in both datasets NS-UNO has the best accuracy. Compared with the other three models, the error is uniformly small, i.e., the numerical solution only largely deviates from the exact solution at few points, further indicating the effectiveness of the proposed NS-UNO.

\subsubsection{Higher Wavenumber scenario}

\begin{figure}
\centering
\setcounter{subfigure}{0}
\subfigure[$q$ T-shaped, $f$ Gaussian(30)]{
\includegraphics[width=0.45\textwidth]{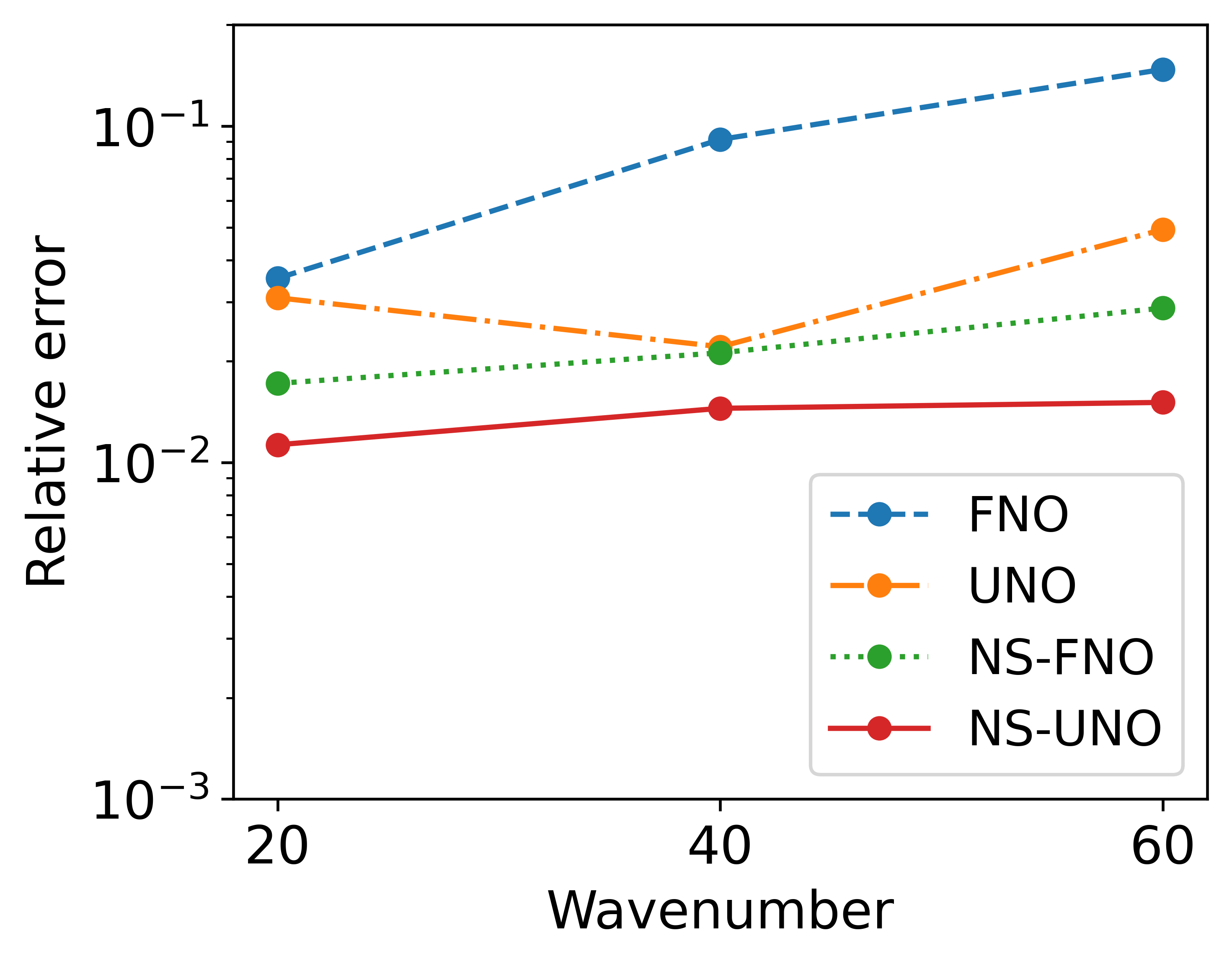}}
\hspace{0.04in}
\centering
\subfigure[$q$ T-shaped, $f$ GRF]{
\includegraphics[width=0.45\textwidth]{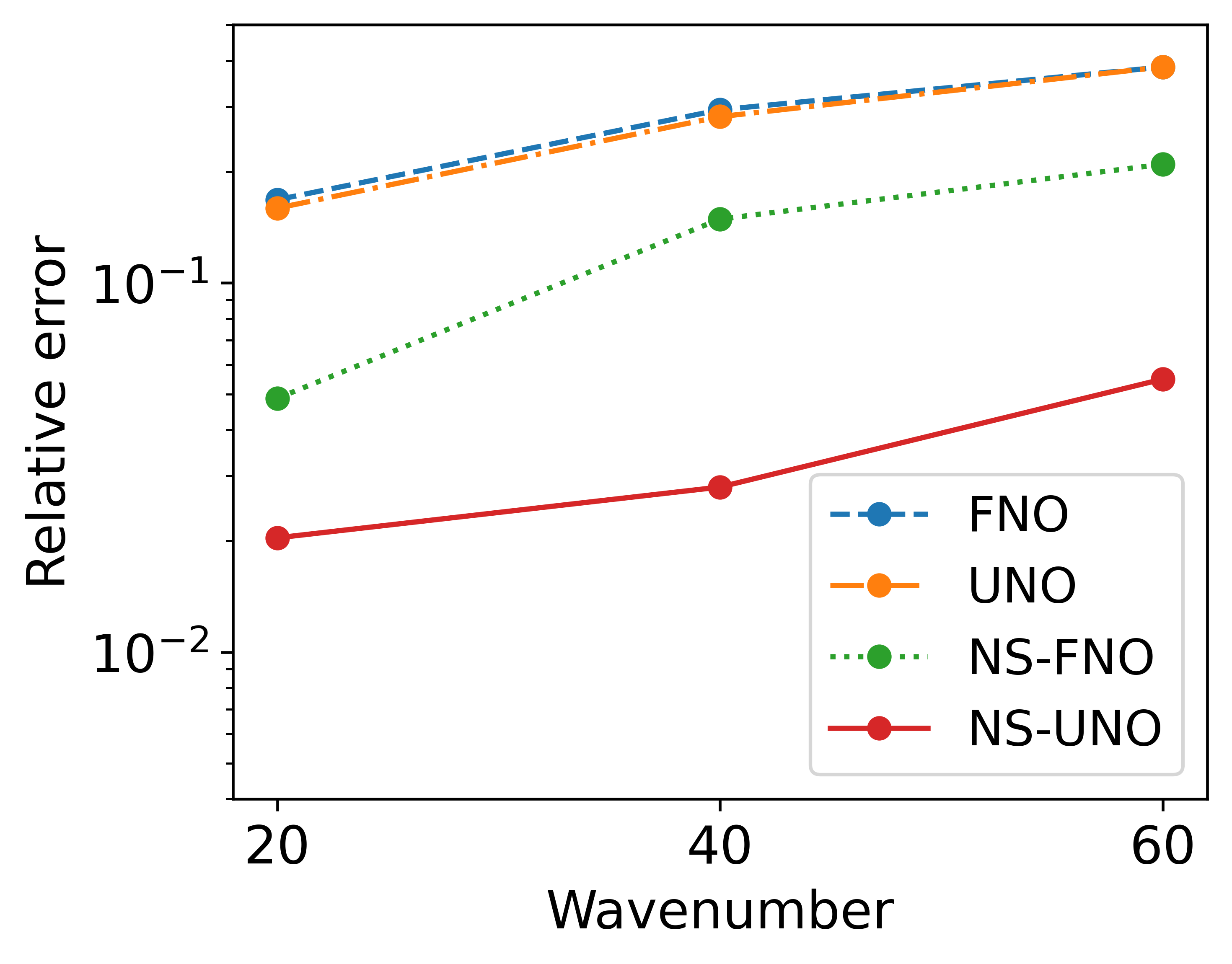}}

\subfigure[$q$ random circles, $f$ GRF]{
\includegraphics[width=0.45\textwidth]{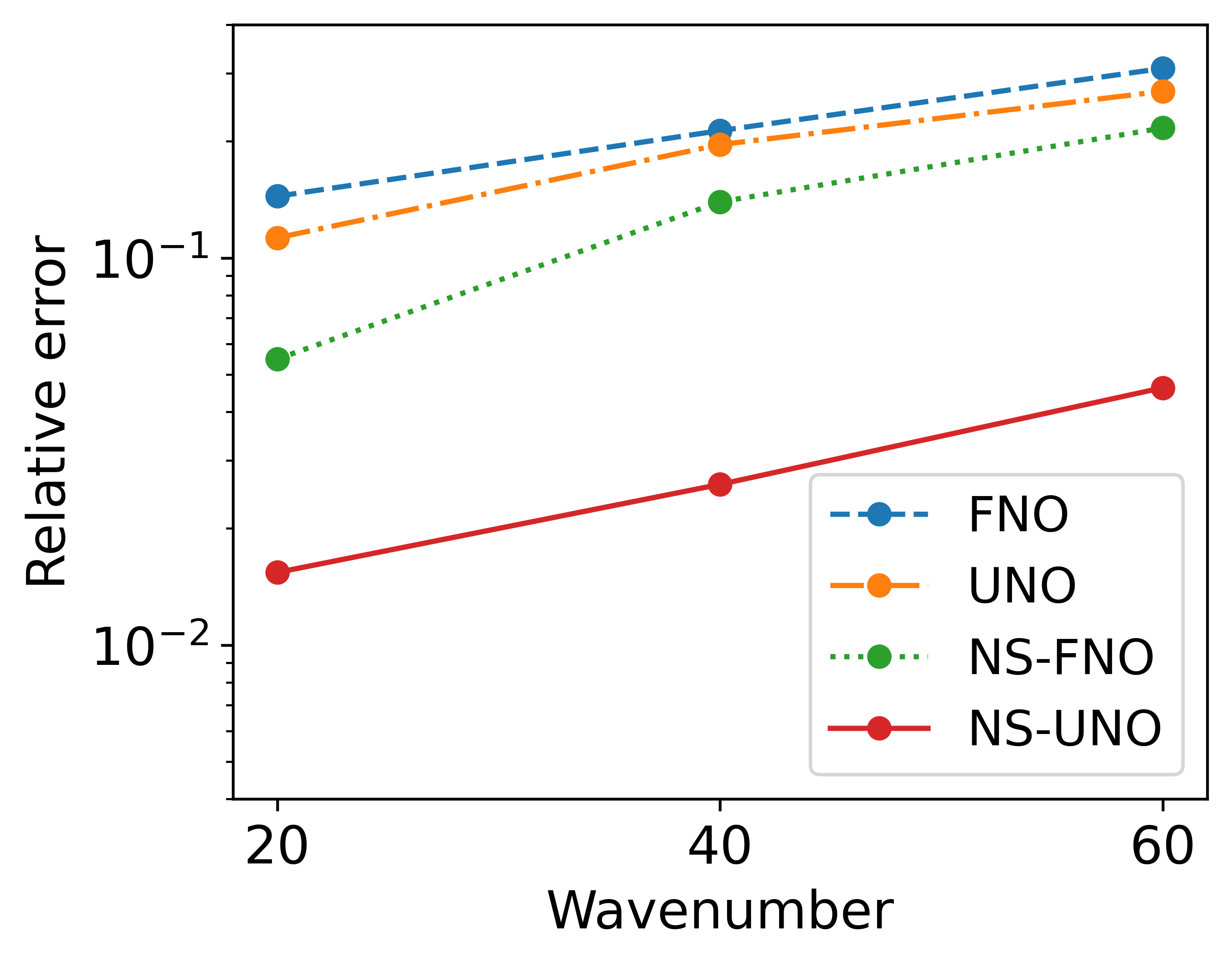}}
\hspace{0.04in}
\centering
\subfigure[$q$ smoothed circles, $f$ wave]{
\includegraphics[width=0.45\textwidth]{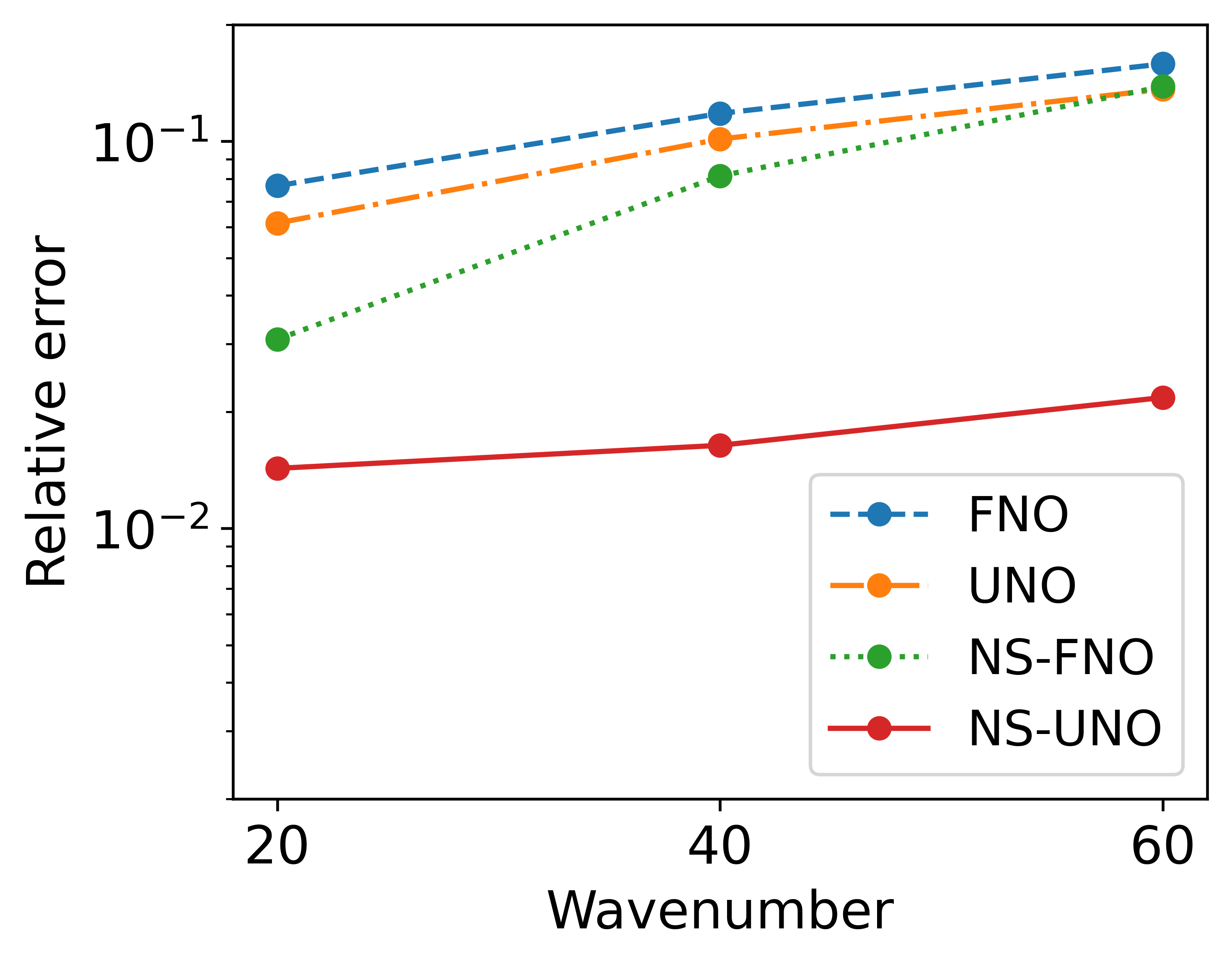}}
\caption{Relative $L^2$-error of FNO, UNO, NS-FNO and NS-UNO on four datasets specified in the subcaptions for wavenumber $k=20, 40, 60$.}
\label{fig:k}
\end{figure}

To further test the performance of the models on multi-scale problems, we further present the results for $k=40, 60$ in the Fig. \ref{fig:k}. It can be seen that among the four models, NS-UNO has the lowest relative error, followed by NS-FNO, indicating the effectiveness of the proposed NSNO. Besides, for the datasets with $f$ sampled from Gaussian random field or superposition of planar waves, which have significant multi-scale properties, NS-UNO shows superior advantage over the other models with nearly one order of magnitude lower relative error, which further reveals the potential of UNO architecture in dealing with multi-scale problems.

\begin{figure}
\centering
\subfigure[Ground truth, real part \qquad]{
\includegraphics[width=0.3\textwidth]{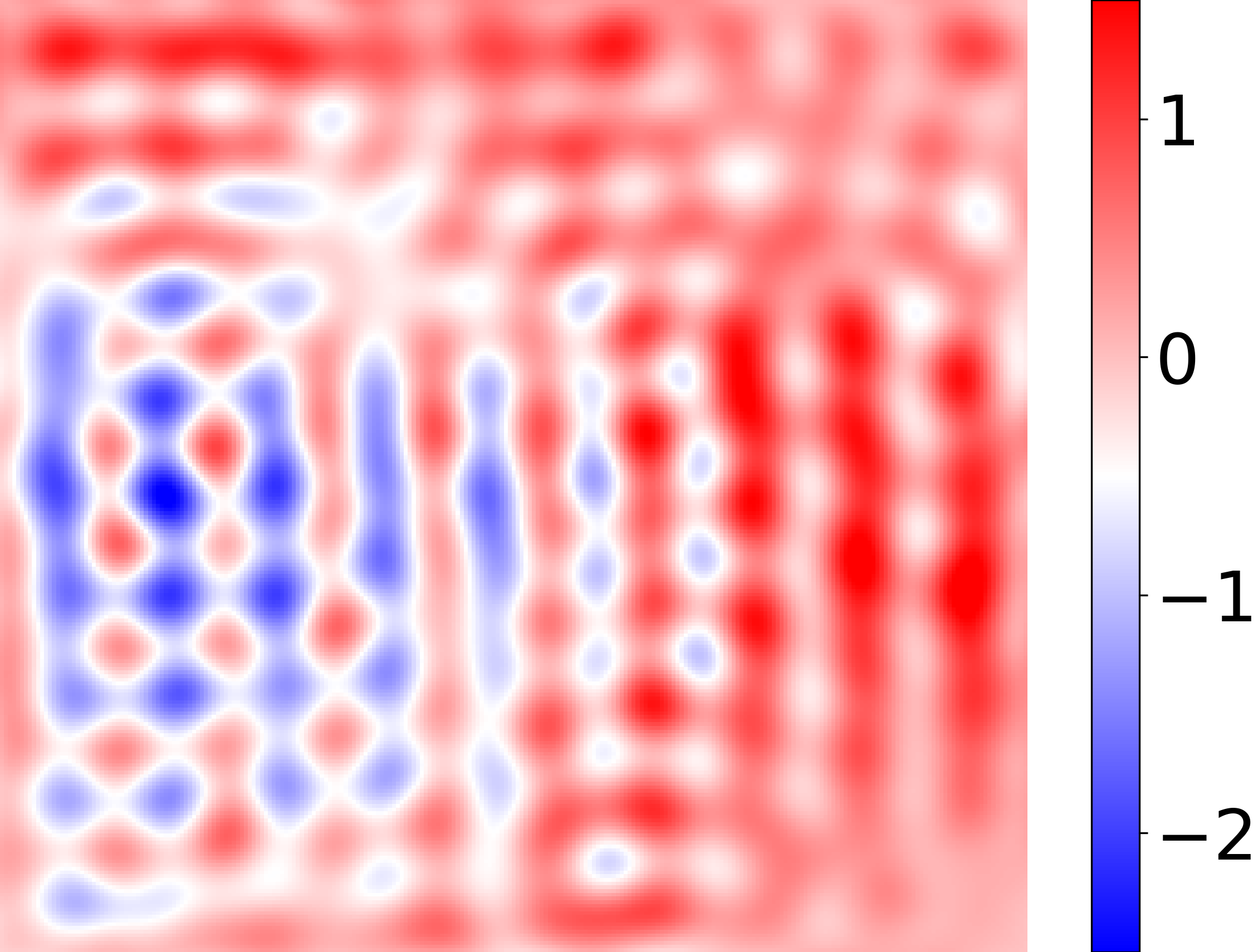}}
\hspace{0.04in}
\centering
\subfigure[Ground truth, imaginary part \qquad]{
\includegraphics[width=0.3\textwidth]{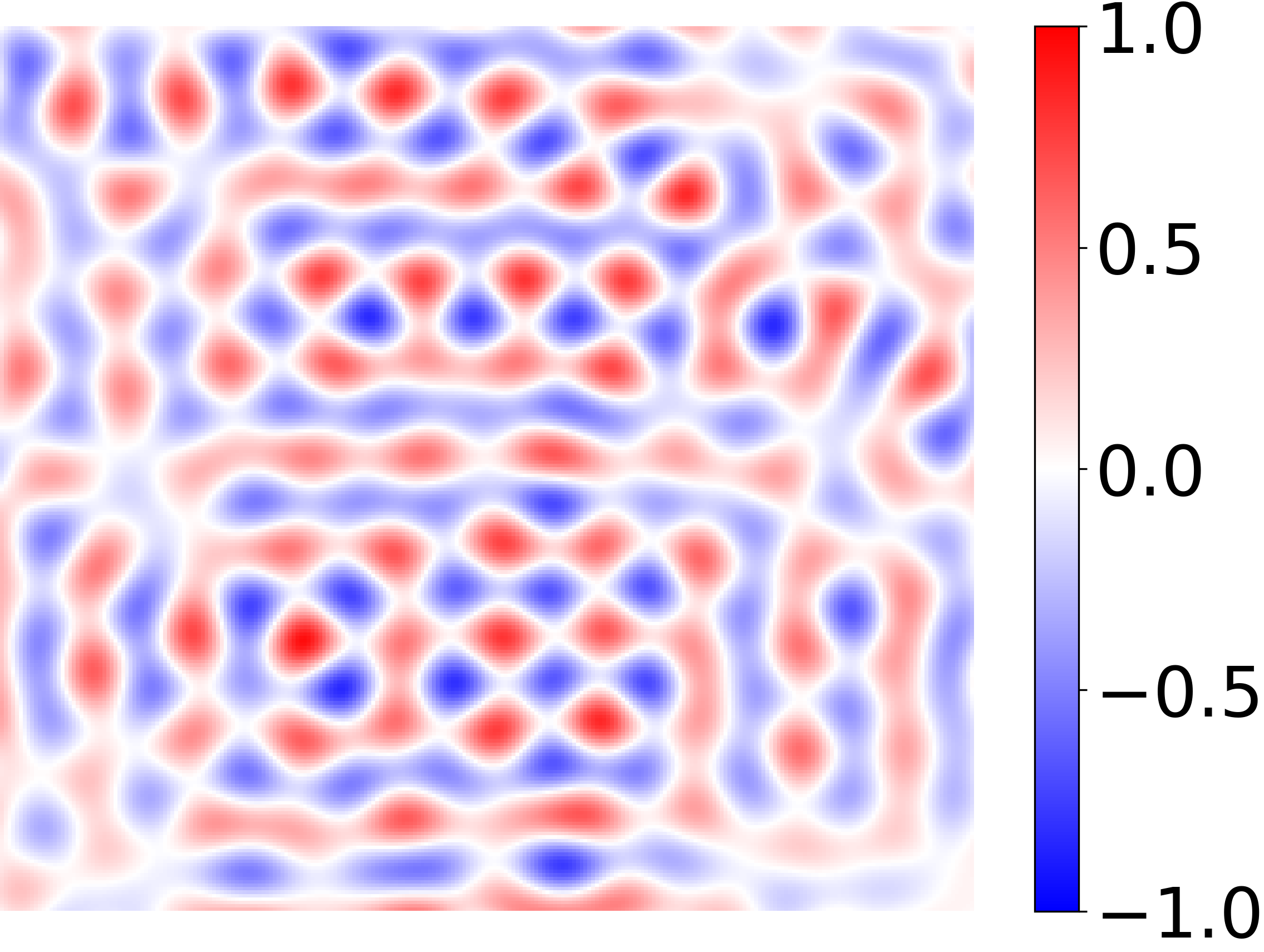}}
\hspace{0.04in}
\centering
\subfigure[$q$ (top left) and $f$\qquad]{
\includegraphics[width=0.3\textwidth]{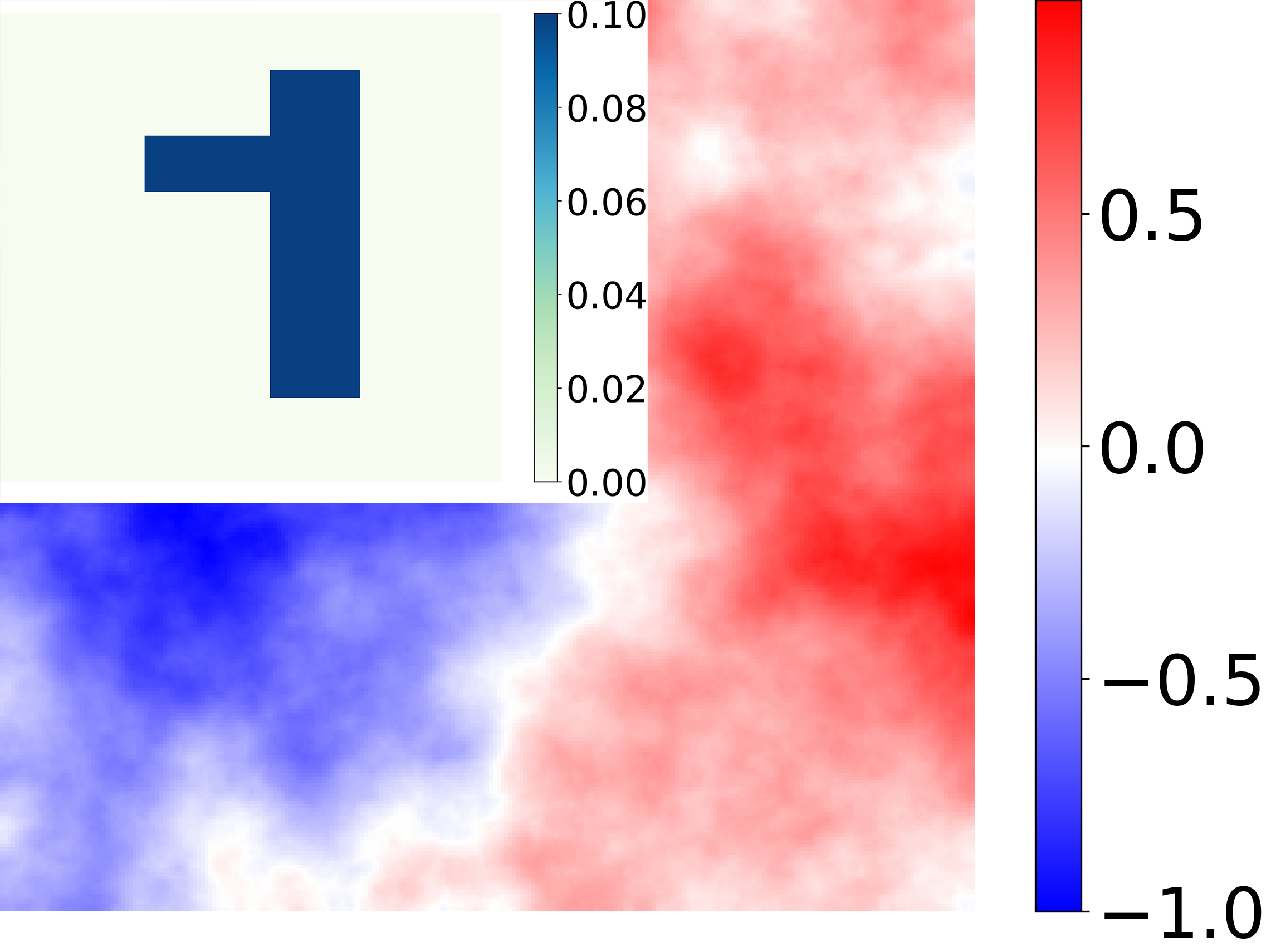}}

\subfigure[FNO, real part \quad]{
\includegraphics[width=0.3\textwidth]{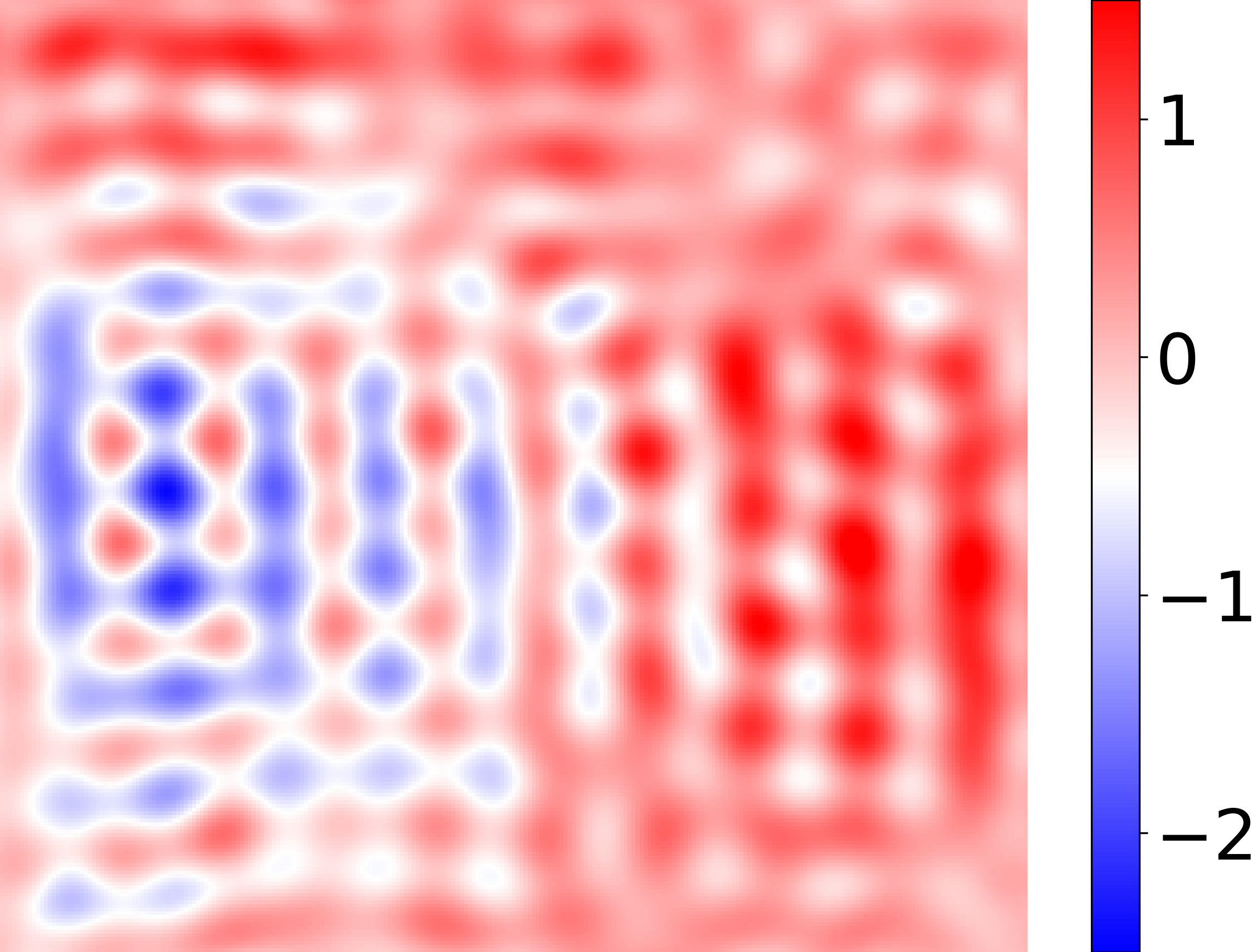}}
\hspace{0.04in}
\centering
\subfigure[FNO, imaginary part\qquad]{
\includegraphics[width=0.3\textwidth]{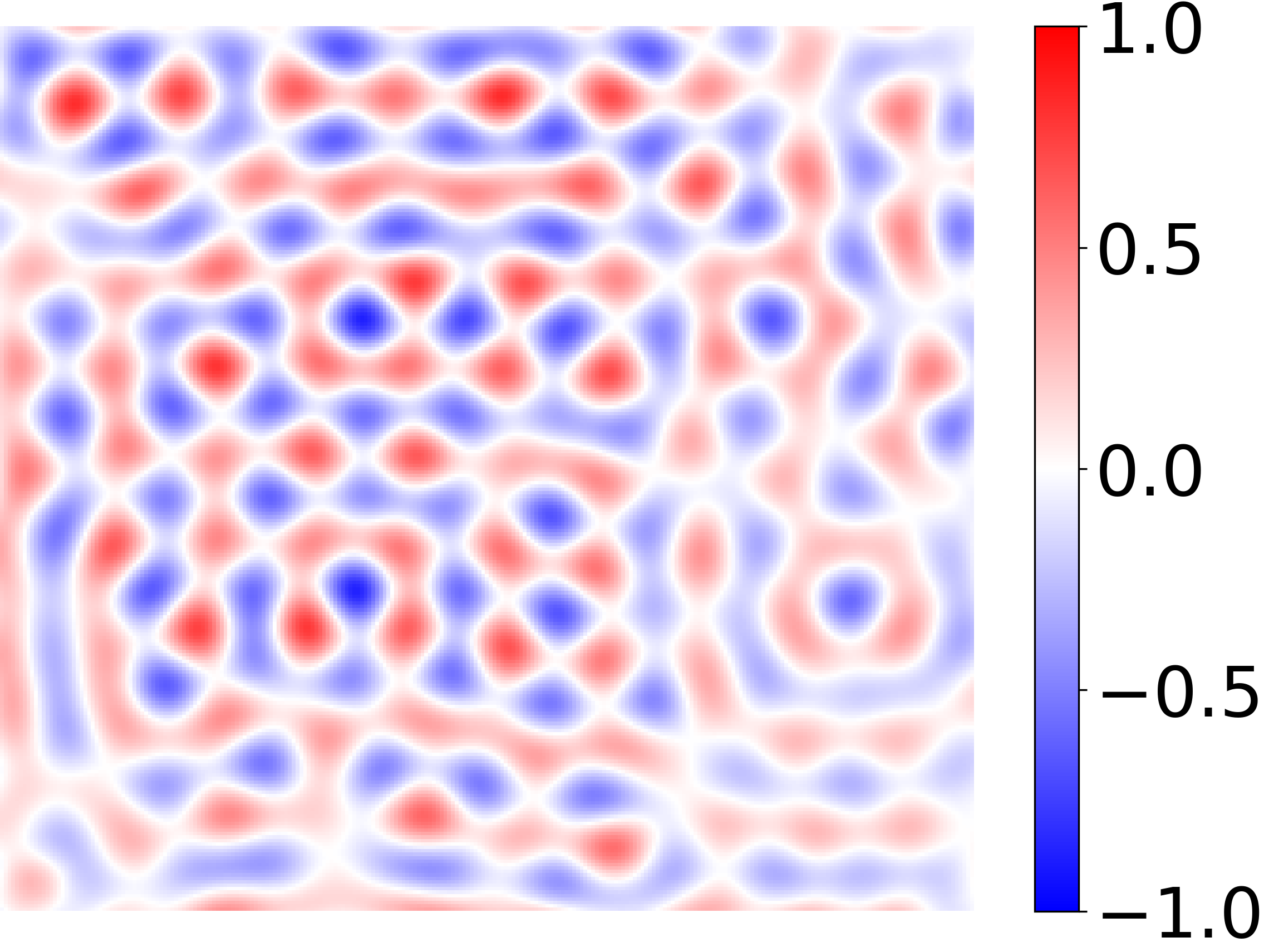}}
\hspace{0.04in}
\centering
\subfigure[FNO, error, 42.10\% \qquad]{
\includegraphics[width=0.3\textwidth]{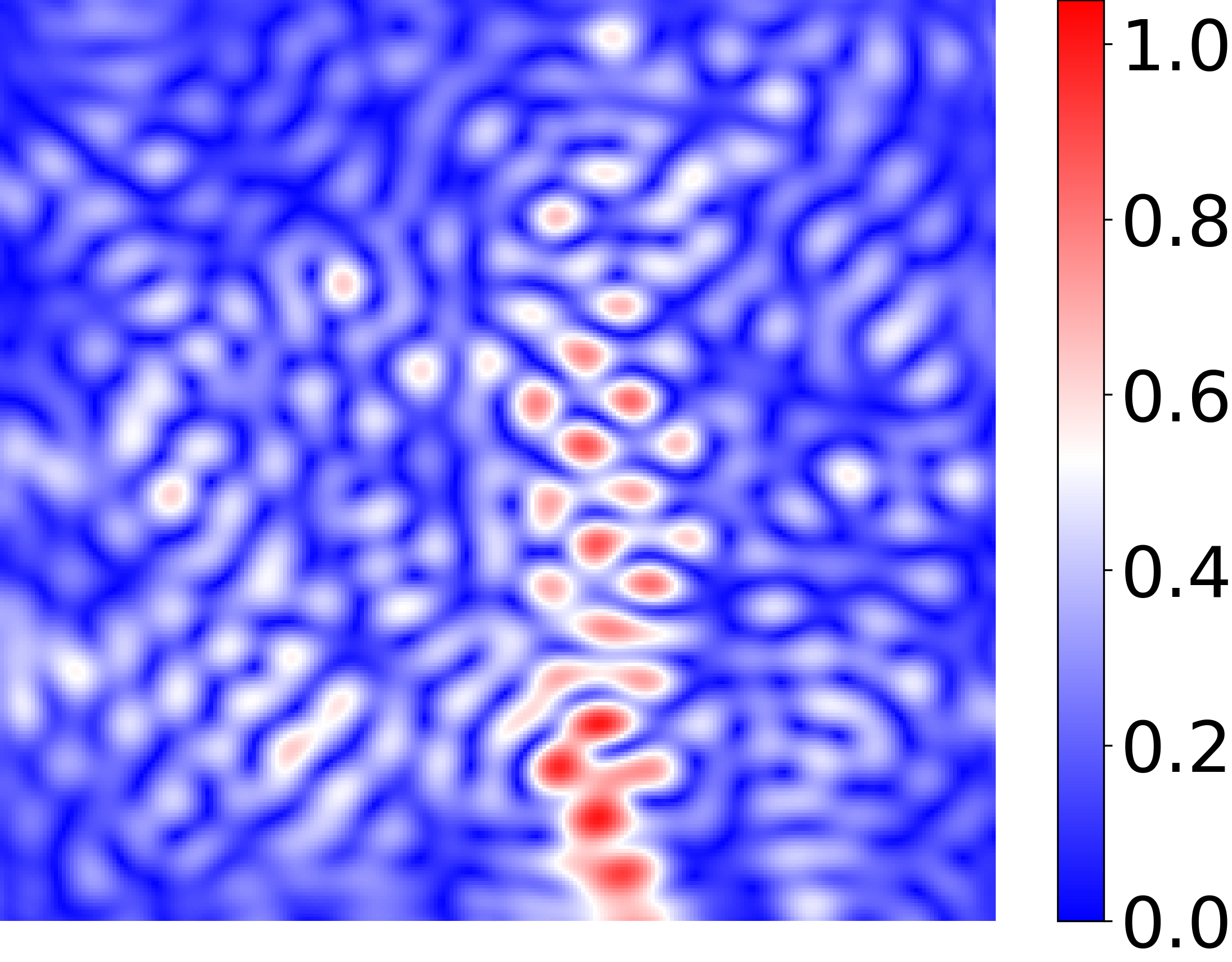}}

\subfigure[UNO, real part \quad]{
\includegraphics[width=0.3\textwidth]{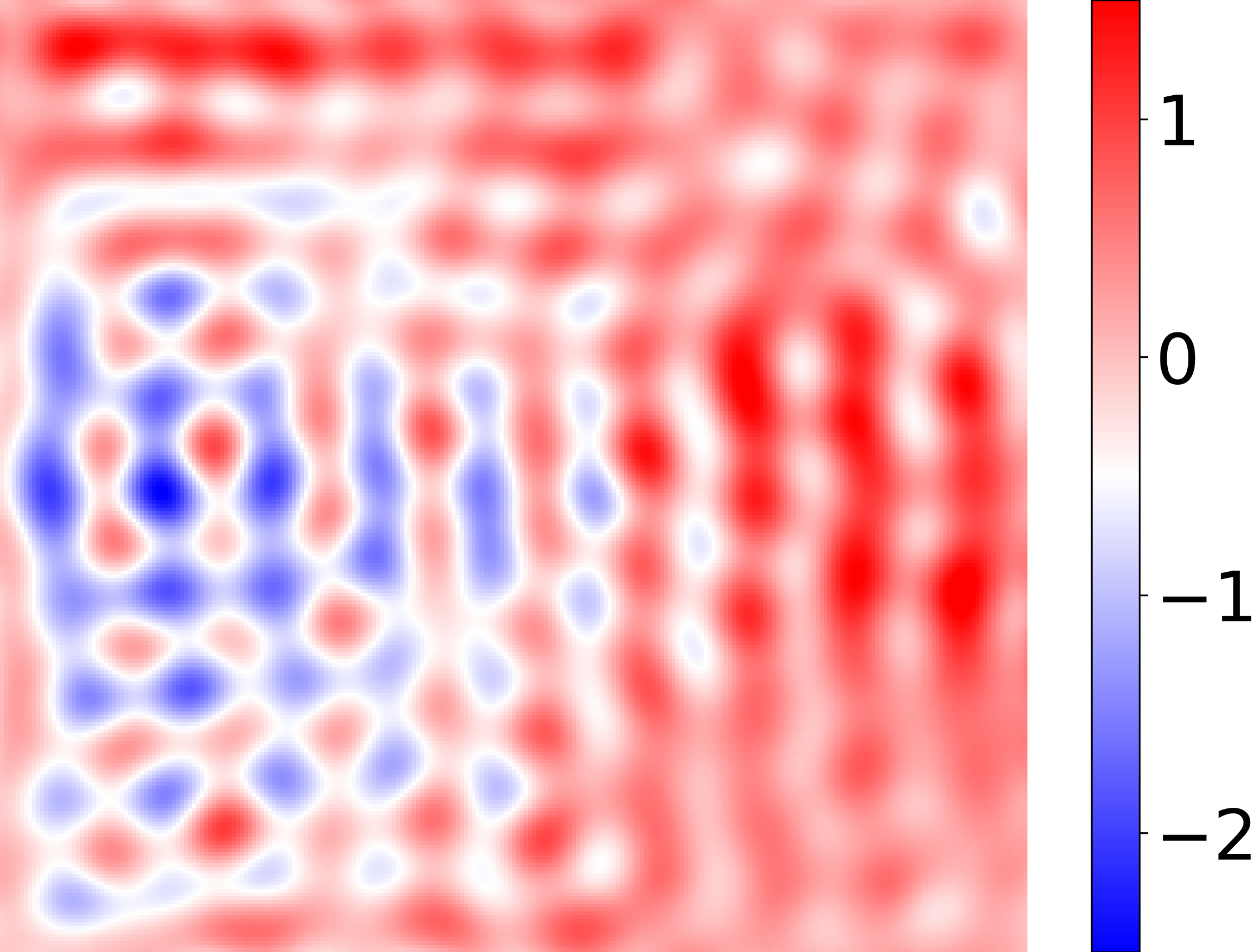}}
\hspace{0.04in}
\centering
\subfigure[UNO, imaginary part\qquad]{
\includegraphics[width=0.3\textwidth]{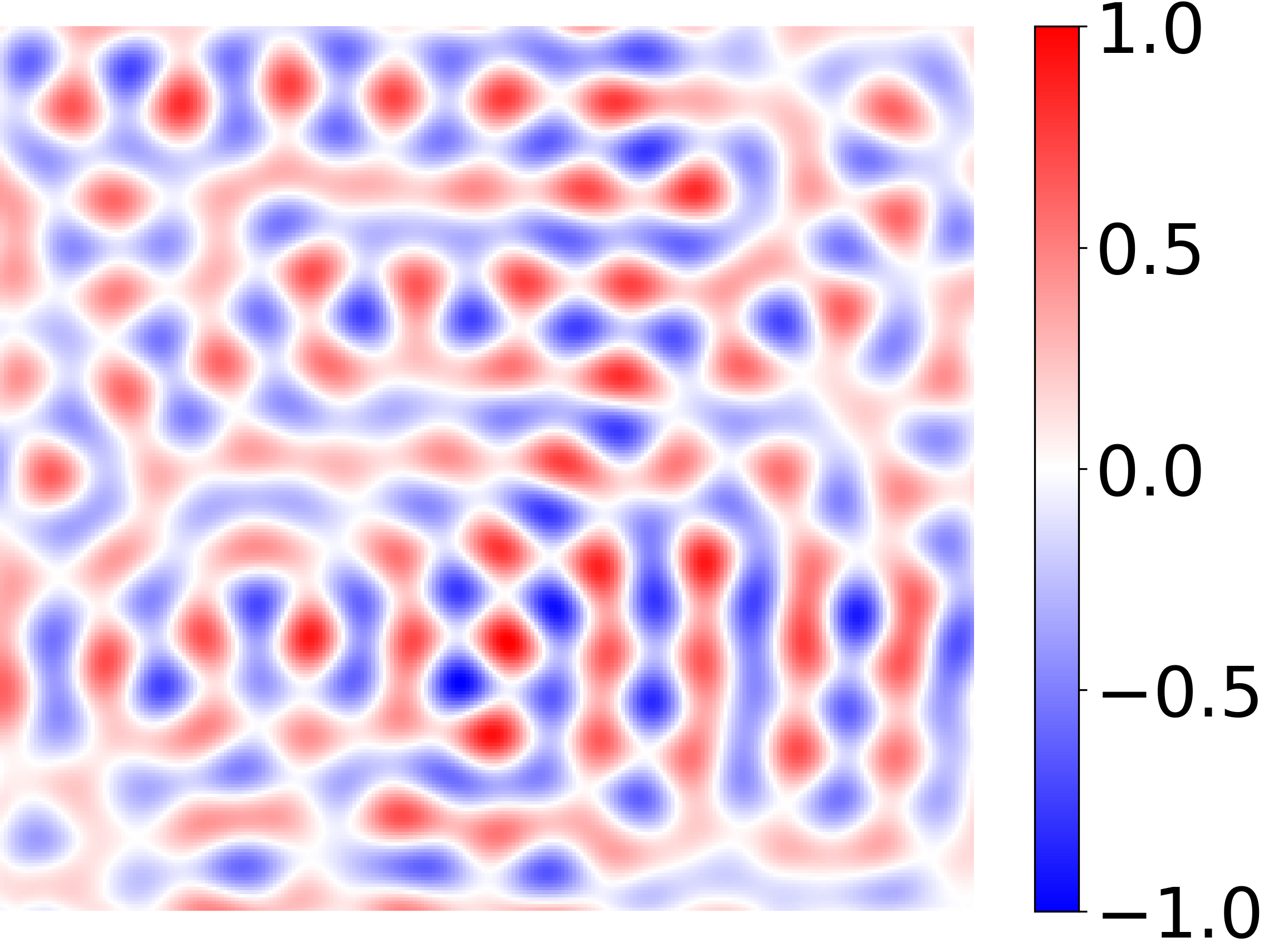}}
\hspace{0.04in}
\centering
\subfigure[UNO, error, 36.42\% \qquad]{
\includegraphics[width=0.3\textwidth]{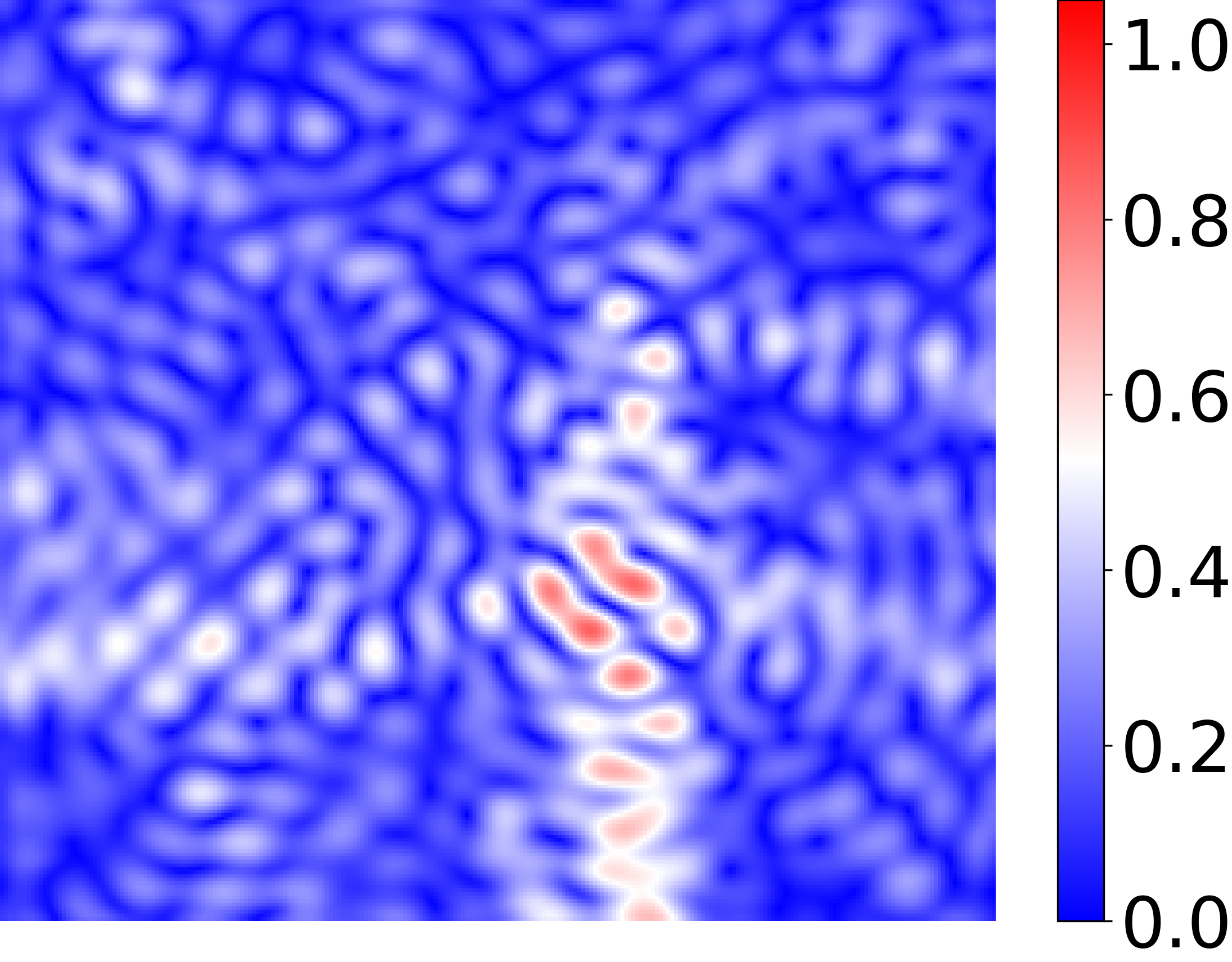}}

\subfigure[NS-FNO, real part \quad]{
\includegraphics[width=0.3\textwidth]{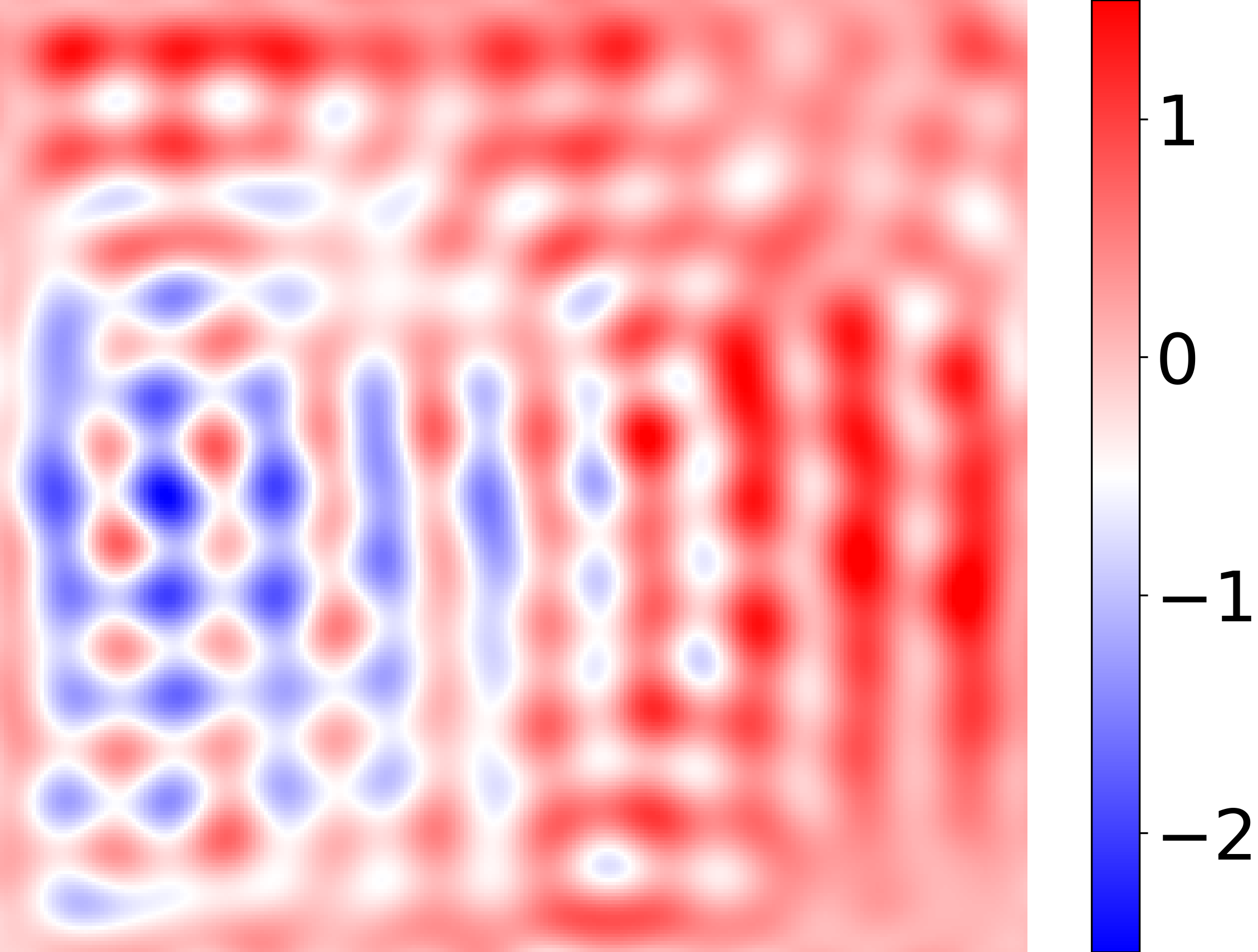}}
\hspace{0.04in}
\centering
\subfigure[NS-FNO, imaginary part\qquad]{
\includegraphics[width=0.3\textwidth]{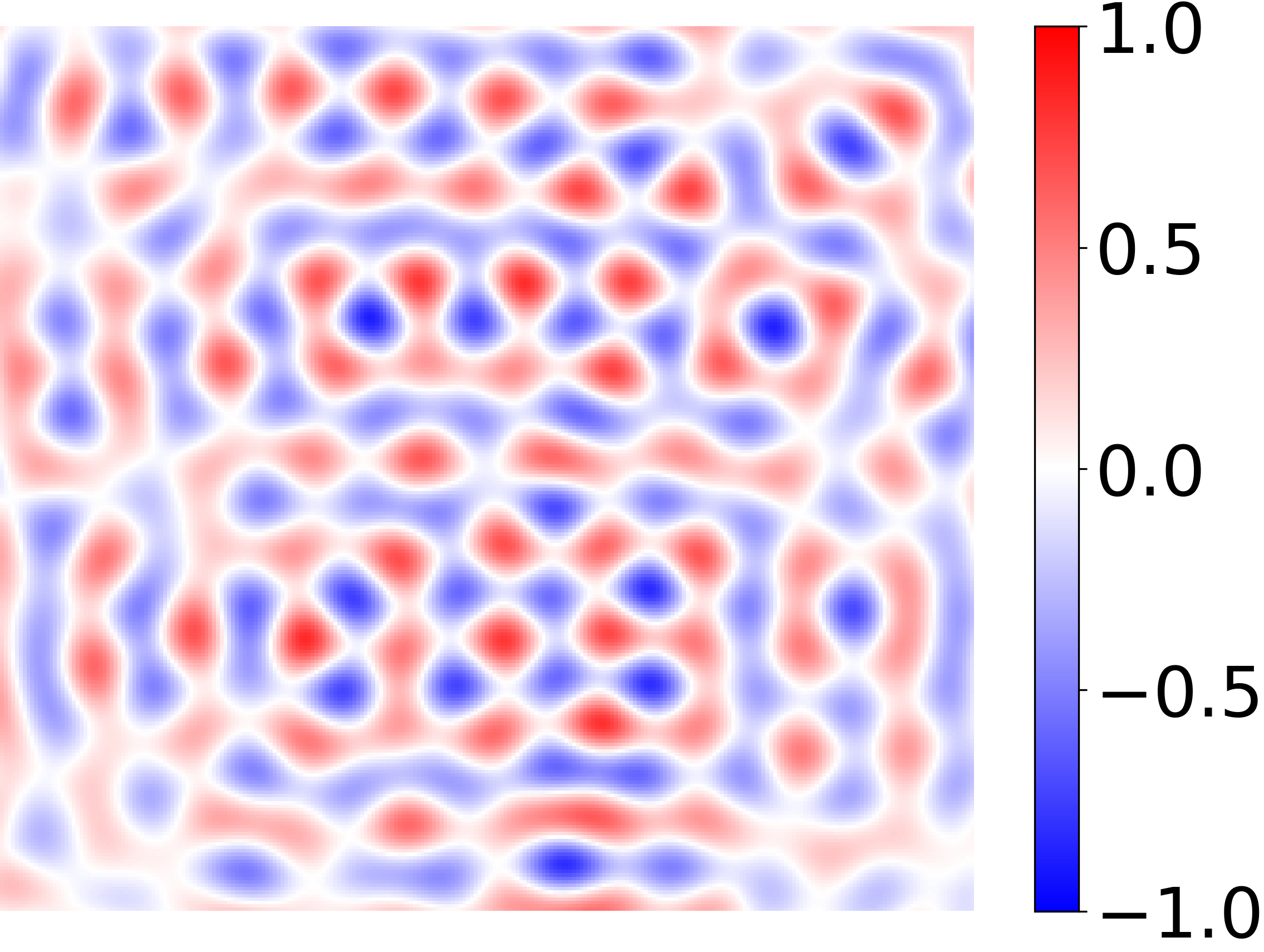}}
\hspace{0.04in}
\centering
\subfigure[NS-FNO, error, 17.82\% \qquad]{
\includegraphics[width=0.3\textwidth]{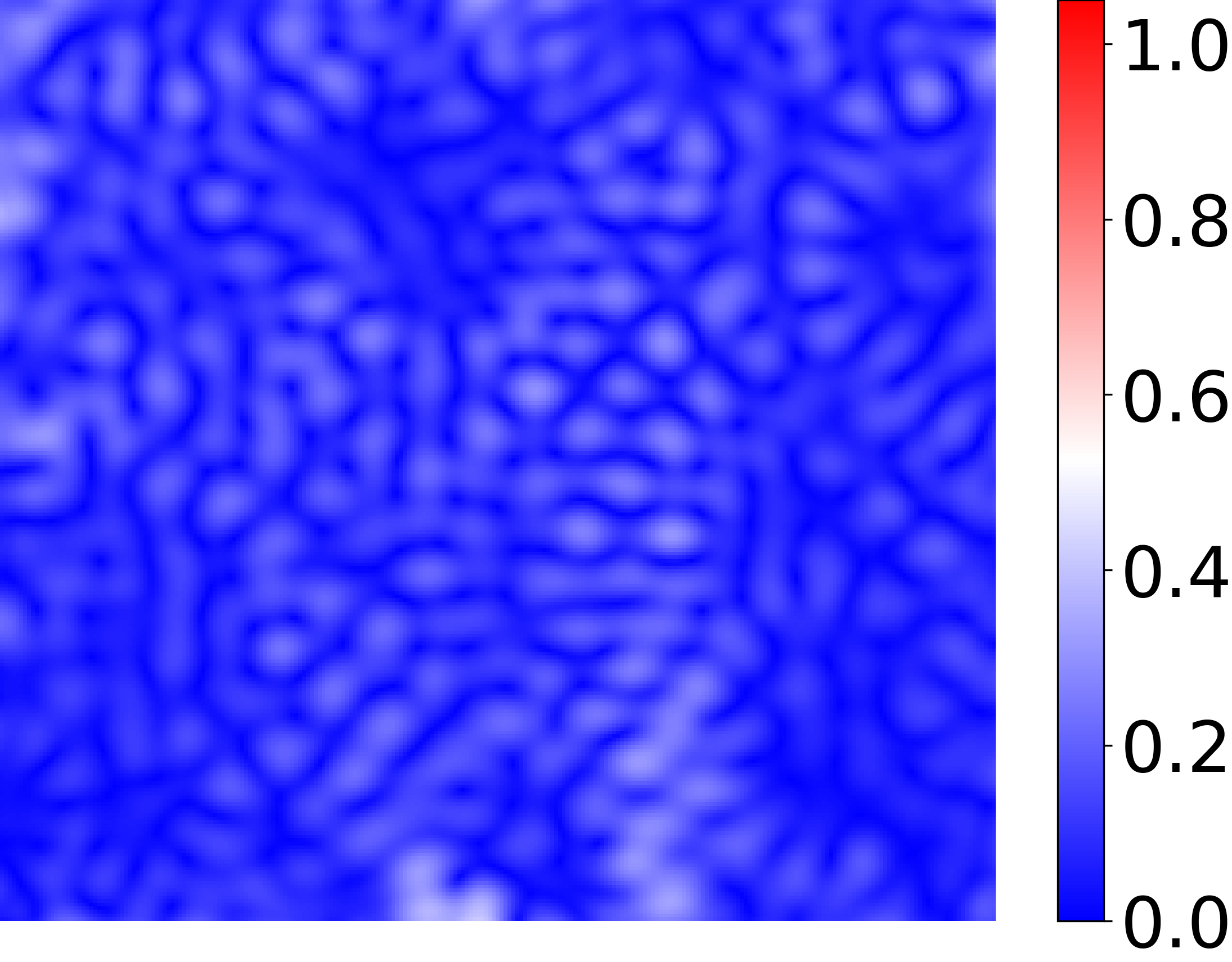}}

\subfigure[NS-UNO, real part \quad]{
\includegraphics[width=0.3\textwidth]{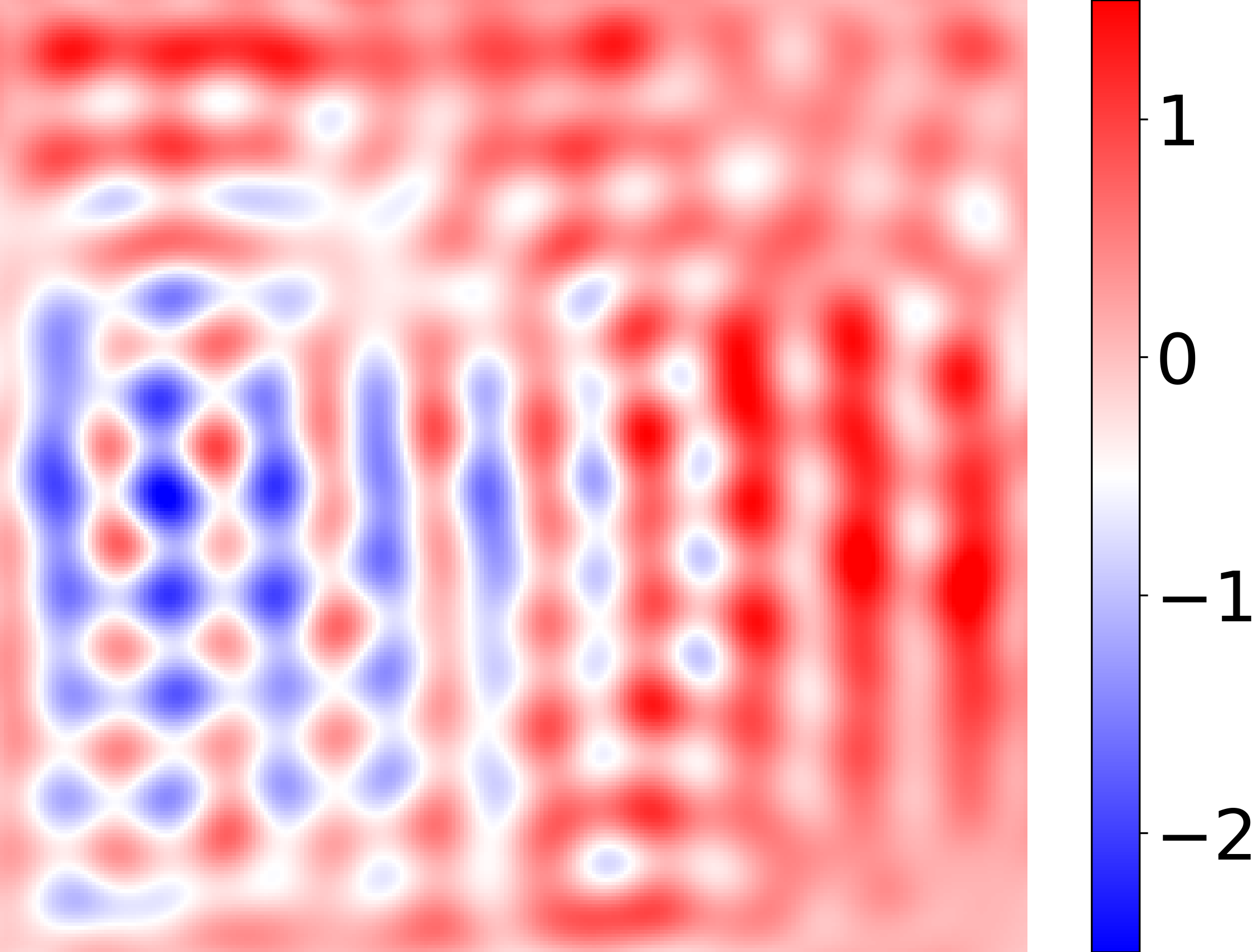}}
\hspace{0.04in}
\centering
\subfigure[NS-UNO, imaginary part\qquad]{
\includegraphics[width=0.3\textwidth]{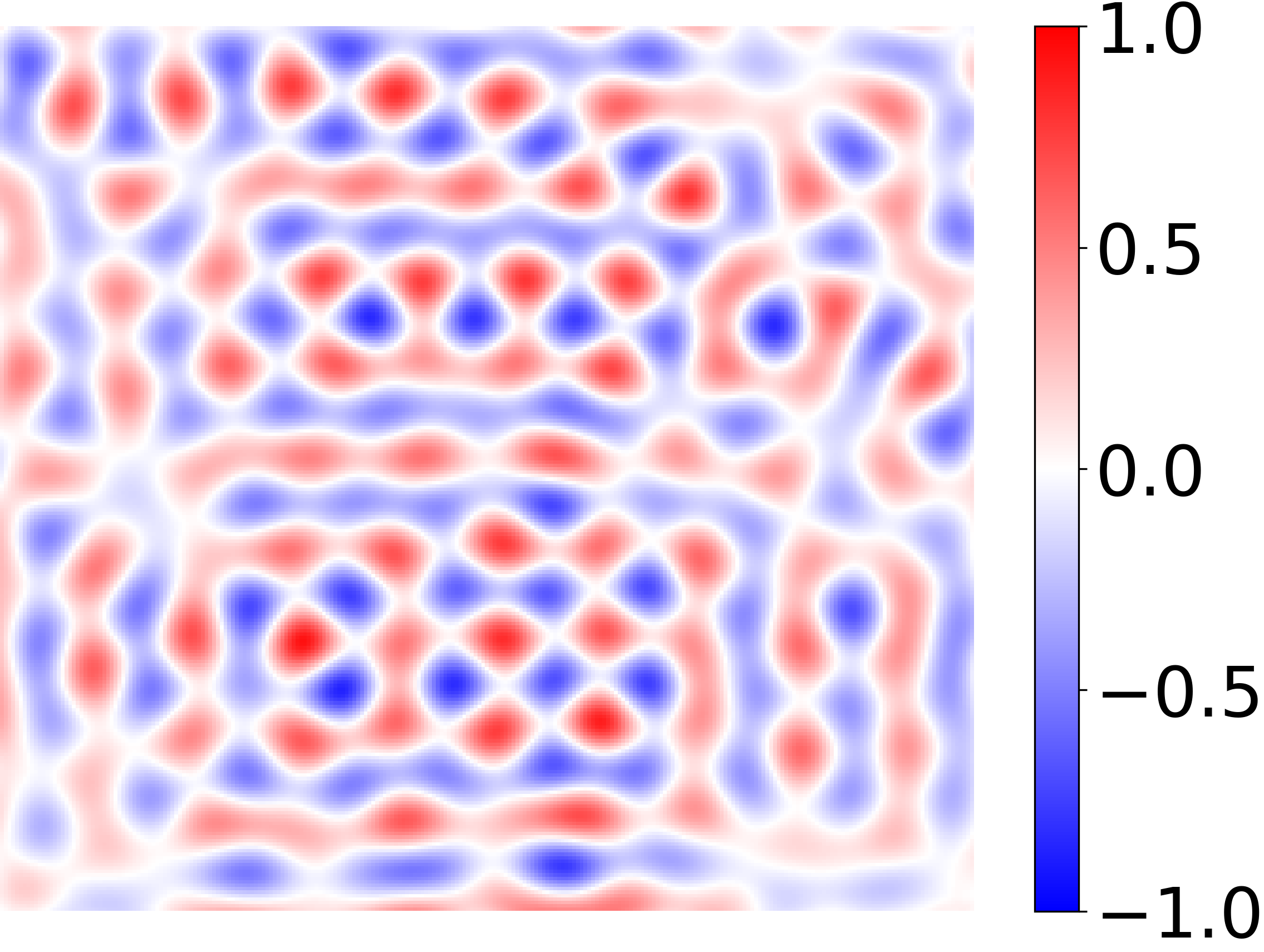}}
\hspace{0.04in}
\centering
\subfigure[NS-UNO, error, 5.40\% \qquad]{
\includegraphics[width=0.3\textwidth]{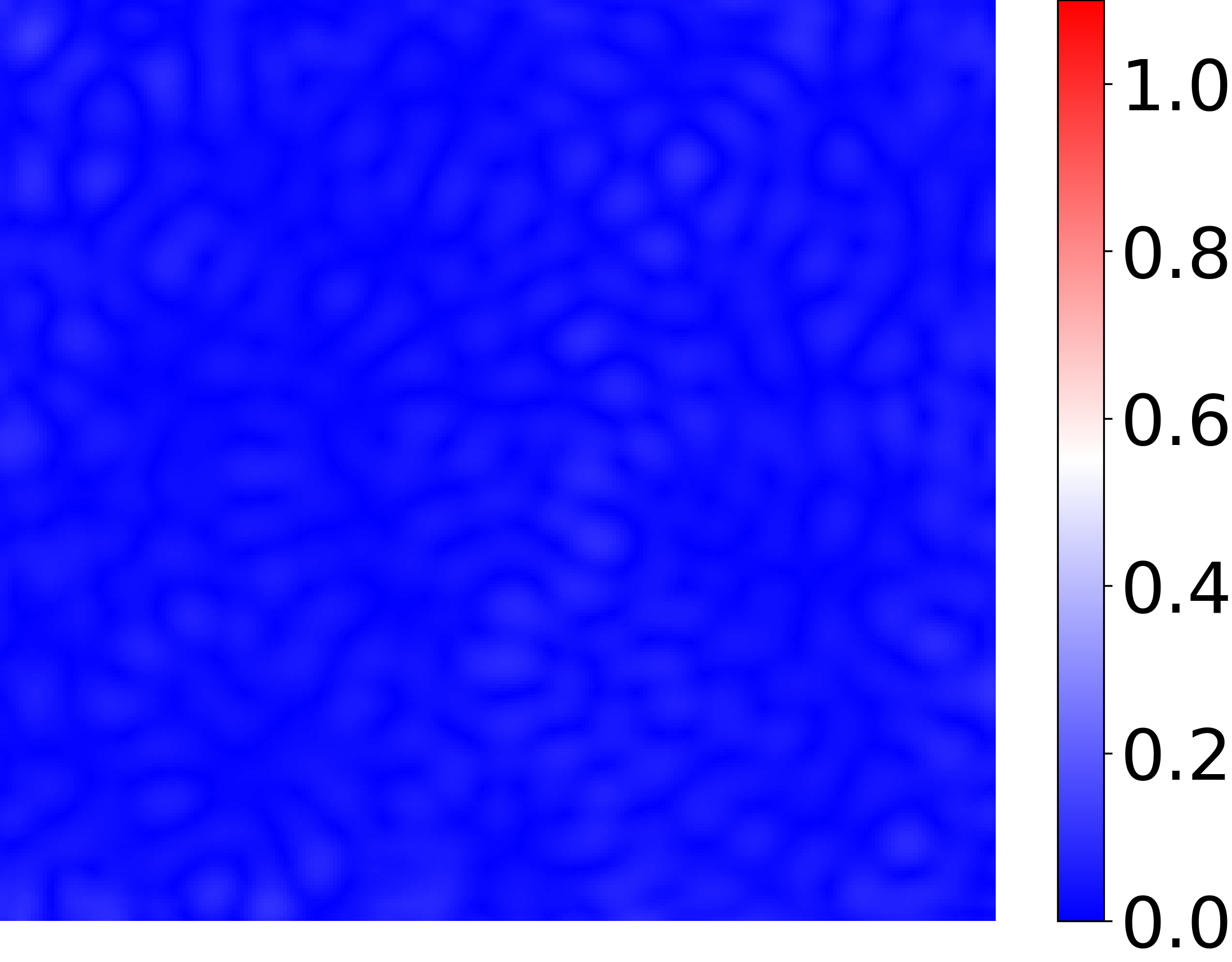}}

\caption{Example of exact solution, numerical solutions and absolute error for dataset with T-shaped $q$ and GRF $f$ when $k=60$}
\label{eg3}
\end{figure}

An example taken from the dataset with $q$ being smooth circles and $f$ being waves is shown in Fig. \ref{eg3} for $k=60$. It can be seen that compared with the $k=20$ cases shown in Fig. \ref{eg1} and Fig, \ref{eg2}, the solutions are more complicated and oscillatory. However, the proposed NS-UNO is still able to learn the solutions with uniformly small error, showing superior power in capturing the multi-scale features of the solution the Helmholtz equations.

\subsubsection{Less Training Data}

\begin{figure}
\centering
\setcounter{subfigure}{0}
\subfigure[$q$ T-shaped, $f$ Gaussian(30)]{
\includegraphics[width=0.45\textwidth]{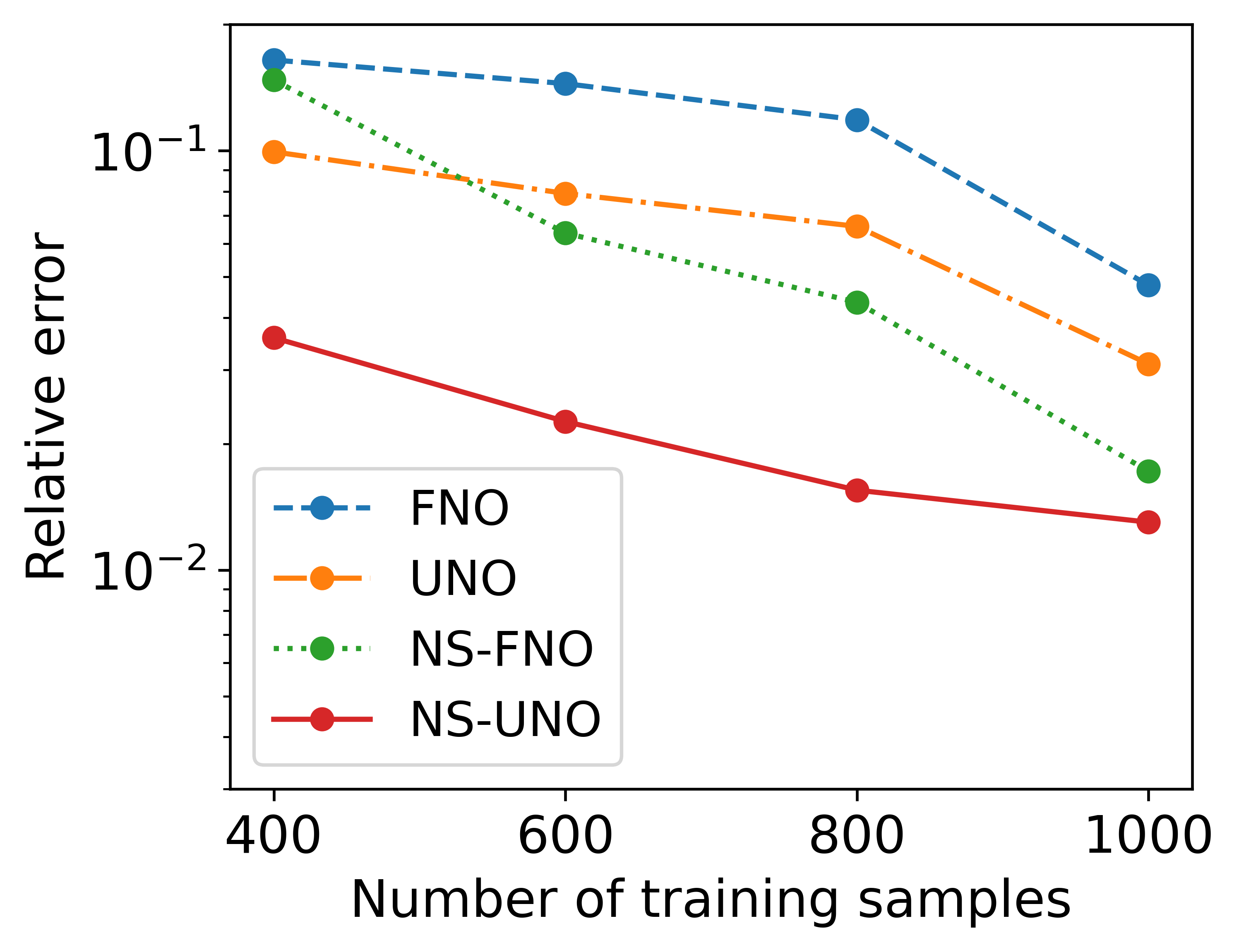}}
\hspace{0.04in}
\centering
\subfigure[$q$ smoothed circles, $f$ wave]{
\includegraphics[width=0.45\textwidth]{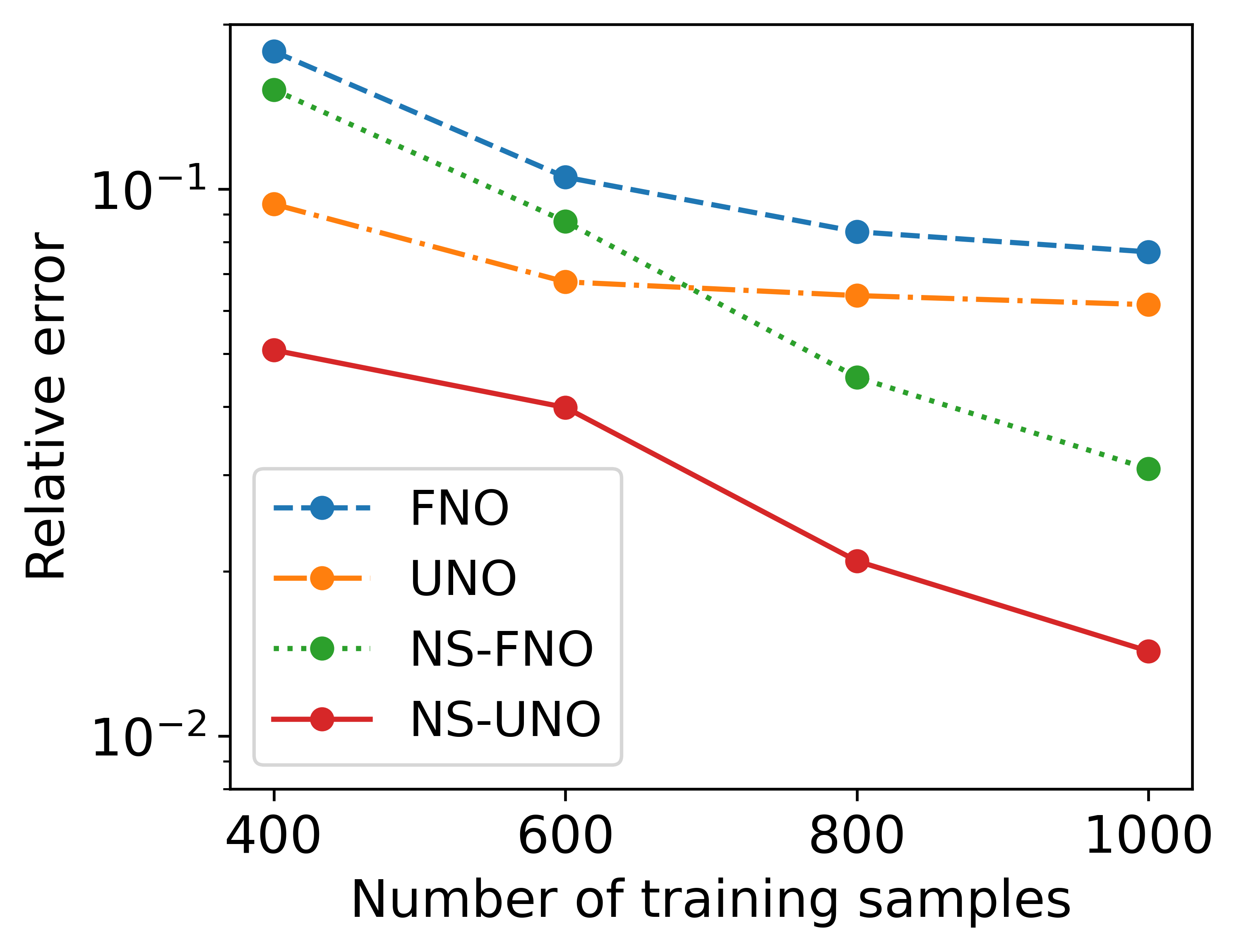}}
\caption{Relative $L^2$-error of FNO, UNO, NS-FNO and NS-UNO on two datasets specified in the subcaptions with the number of training samples $N=400, 600, 800, 1000$.}
\label{fig:data}
\end{figure}

Fig. \ref{fig:data} shows the relative error on test sets when we reduce the size of the training set to 800, 600 and 400. It can be seen that NS-UNO still significantly outperforms the other three models when less data is available, showing that NS-UNO has good data efficiency. We also observe that when there is insufficient data, UNO has lower relative error than NS-FNO, while NS-FNO outperforms UNO with the increase of training samples. This reveals that the difficulty in handling multi-scale problems can be alleviated by increasing the number of training data.

\subsubsection{Training Computational Cost}

A comparison of the training computational costs for NS-FNO and NS-UNO is given in Table \ref{cost}. It can be seen that with similar number of parameters, NS-UNO has 50\% lower computational cost than NS-FNO. This is because the number of channels in the finest level in UNO is a quarter of the number of channels in FNO, which largely reduces the size of tensors, leading to less memory usage and less training time. Therefore, NS-UNO is able to give more accurate results with less computational cost.

\begin{table}
\caption{Training computational costs for NS-FNO and NS-UNO}
\vspace{0.1cm}
\begin{tabular}{|c|ccc|}
\hline
Model  & Iters / sec   & Memory (GB) & \# of param (M)\\
\hline
NS-FNO & 26 & 12.64 & 3.56\\
\hline
NS-UNO & 12 & 5.94 & 3.55\\
\hline
\end{tabular}
\label{cost}
\end{table}

\subsection{Necessity of Physics-Informed Loss}
In this section, we conducted tests on the weights for the physics-informed loss using the dataset where $q$ is random circle and $f\sim\mathcal{N}(0, (-\Delta+9I)^{-2})$. As is shown in Table \ref{table:pinn}, the optimal $\lambda$ is 0.05, which coincides the choice of $\lambda$ given in section \ref{exsetup}. When $\lambda=0$, i.e., there is no physics-informed loss, both relative error on training and test sets are significantly higher than those with physics-informed loss, showing the necessity of physics-informed loss. 

Moreover, we observe that the gap between training error and test error is considerably narrowed with the presence of physics-informed loss, indicating that the introduction of physics-informed loss can improve the generalization ability. To further investigate on the generalization ability, we test the generalization error of the learned models on a new dataset with $f\sim\mathcal{N}(0, (-\Delta+9I)^{-\frac{3}{2}})$. It can be seen that the although the generalization error is larger than the test error, the increment from test error to generalization error is larger in the case without physics-informed loss.

\vspace{-0.1cm}
\begin{table}
\caption{Relative $L^2$-error ($\times 10^{-2}$) of NS-UNO versus weight for the physics-informed loss}
\vspace{0.1cm}
\begin{tabular}{|c|cccccc|}
\hline
$\lambda$      & 0    & 0.01 & 0.05 & 0.1  & 0.15 & 0.2  \\
\hline
Training error & 3.09 & 1.10 & 1.07 & 1.24 & 1.89 & 1.96 \\
\hline
Test error     & 8.85 & 1.73 & 1.54 & 1.68 & 2.21 & 2.29 \\
\hline
Generalization error & 11.64 & 2.57 & 2.42 & 2.65 & 3.01 & 3.07 \\
\hline
\end{tabular}
\label{table:pinn}
\end{table}

\subsection{Influence of the Number of Items in Neumann Series}\label{nitem}

\begin{figure}
\centering
\includegraphics[width=8cm]{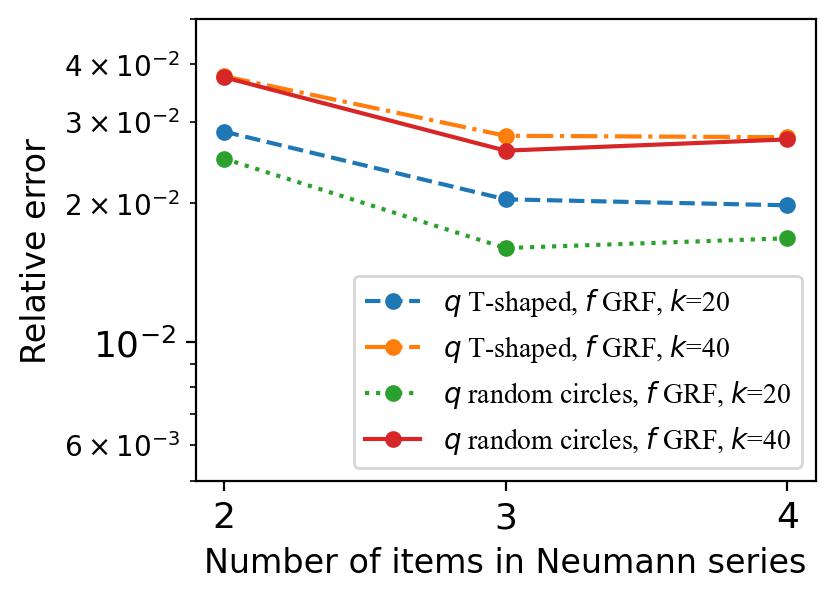}
\caption{Relative $L^2$-error versus the number of items in Neumann series}
\label{ni}
\end{figure}

The relative $L^2$-error versus the number of items to construct the Neumann series is plotted in Fig. \ref{ni}.
It can be seen that three items is sufficient to give an accurate result, and the relative error of the model with four items is almost the same as the model with three items. This is because as the depth of the network increases, the gradient backpropagation has to go through a longer path, leading to gradient vanishing problems\supercite{hochreiter1998vanishing}. The result also indicates that one Neumann series block in the neural network actually does not correspond to a single step in the exact Neumann series.

\subsection{Beyond Convergence of Neumann Series}

As discussed in \eqref{norm} in the proof of the convergence of Neumann series, the convergence rate is associated with the norm of the operator $-k^2Gq$, which is proportional to $k$ and $\Vert q\Vert_{L^{\infty}(\it\Omega)}$. Overly large $k$ and $\Vert q\Vert_{L^{\infty}(\it\Omega)}$ will lead to slow convergence and even divergence of the Neumann series. Surprisingly, we found that the proposed NS-UNO is able to break through this theoretical limitation and still gives reasonable results even the exact Neumann series diverges, as shown in Table \ref{maxq}. 

Specifically, we use the dataset where $q$ is T-shaped and $f$ is sampled from Gaussian(30). We take $k=20$ and $k=40$, and compute the exact Neumann series by MUMPS. It can be seen that for $k=20, \Vert q\Vert_{L^{\infty}(\it\Omega)}=0.35$ and $k=40, \Vert q\Vert_{L^{\infty}(\it\Omega)}=0.2$, the Neumann series exhibits a slow convergence, while for $k=20, \Vert q\Vert_{L^{\infty}(\it\Omega)}=0.4$ and $k=40, \Vert q\Vert_{L^{\infty}(\it\Omega)}=0.25$, the Neumann series diverges, which coincides with the theoretical discussion in Section \ref{NS_theory}. However, the proposed NS-UNO outperforms the exact Neumann series with 10 terms using only three iteration steps, and can still maintain a relatively low $L^2$-error even the Neumann series diverges. Besides, NS-UNO still outperforms FNO for large $k$ and $\Vert q\Vert_{L^{\infty}(\it\Omega)}$, which is not affected by this  theoretical limitation of Neumann series.



\begin{table}
\caption{Relative $L^2$-error ($\times 10^{-2}$) of NS-UNO and Exact Neumann series (NS)}
\vspace{0.1cm}
\begin{tabular}{|c|c|cc|cc|}
\hline
wavenumber & $\Vert q\Vert_{L^{\infty}(\it\Omega)}$ & NS (3 terms) & NS (10 terms) & FNO & NS-UNO \\
\hline
\multirow{2}{*}{$k=20$} & 0.35 & 24.39 & 6.85 & 10.56 & \textbf{6.42}  \\
& 0.4 & 44.11 & 57.12 & 17.90 & \textbf{9.83}  \\
\hline
\multirow{2}{*}{$k=40$} & 0.2 & 12.42 & 3.10 & 17.27 & \textbf{2.92}  \\
& 0.25 & 25.16 & 35.52 & 20.19 & \textbf{4.28} \\
\hline
\end{tabular}
\label{maxq}
\end{table}


\section{Application in Inverse Scattering Problem}\label{invp}
In this section, we solve an inverse scattering problem using the learned NSNO as the forward solver. We first demonstrate the setups for the inverse problem, including the governing equation and data measurement. Traditional and neural network-based methods for solving the inverse problem are then introduced. Following that, we show the reconstructed results using traditional finite difference method and the proposed NS-UNO as the forward solver, respectively.
\subsection{Problem Setup}
We set the incident field as the plane wave $u^i = e^{i k \mathbf{x} \cdot \mathbf{d}}$, where $\mathbf{d}$ denotes the incoming direction. In our experiments, we discrete $\it\Omega$ to a $128\times 128$ grid and generate plane waves from $M=32$ different directions uniformly. Based on \eqref{op1}, the scattered field satisfies
\begin{equation}\label{scatter field}
\Delta u^s+k^2(1+q(x))u^s=-k^2q(x)u^i.
\end{equation}
The mapping from the scatter $q$ to the scattered field $u^s$ can then be given by the solution operator defined in \eqref{operator} as $q\mapsto \mathcal{S}(q, -k^2qu^i)$.

As for data, sensors are placed on each grid on the boundary of $\it\Omega$$=[0, 1]^2$ to collect the wave field data $d_m$ for each incident wave $u_m^i,$ $m=1,2,\dots,M$. Therefore, the forward operator maps to scatter $q$ to the data collected on the boundary $\mathcal{F}_{m}(q)=T\circ\mathcal{S}(q, -k^2qu_m^i)$, where $T$ denotes the trace operator restricting the wave field on the boundary. The inverse problem is to reconstruct the scatterer $q$ from the measured data. We generate the wave field measurement with a fine grid to avoid the inverse crime.

\subsection{Solving the Inverse Problem}
We solve the inverse problem by the traditional optimization approach, which means recovering an approximation of $q$ by solving the following optimization problem
\begin{equation}
\underset{q}{\operatorname{argmin}} J(q) = \sum_{m=1}^M J_j(q) = \sum_{m=1}^M \frac{1}{2} \| \mathcal{F}_{m}(q) - d_m\|_2^2
\end{equation}

We employ the L-BFGS algorithm\supercite{liu1989limited} to address the minimization problem, initializing the iterative process with a value of 0. The gradient of the loss with respect to the model is computed using the adjoint state method\supercite{plessix2006review}, which will be detailed in the appendix. Our objective is to evaluate the effectiveness of utilizing a pre-trained neural network embedded as a forward solver within the optimization framework, as opposed to the conventional numerical solvers, for which we choose the finite difference method (FDM). For the relevant finite difference matrix, we select MUMPS as the direct solver.

The neural network training process is outlined below. The dataset consists of scatterer samples $q$, which are extracted from the large MNIST dataset\supercite{jansson2021exploring}. Specifically, these scatterer samples are initially resized to a resolution of $112\times 112$ from the original MNIST dataset\supercite{deng2012mnist} and further padded to $128\times 128$. To construct the training/test sets, we select 100/10 samples for each digit in the range 0-9. During the training process, only four plane waves with directions $0, \frac{\pi}{2}, \pi, \frac{3\pi}{2}$ are employed. The neural network is trained for a total of 1000 epochs with the learning rate initialized as 0.001 and halved every 200 epochs. The weight $\lambda$ in the loss function \eqref{loss} is set as 0.1. After the training is completed, the pre-trained neural network is evaluated on the test set and achieves an average relative $L^2$-error of 1.34\%.

\subsection{Reconstructing results}

\begin{figure}
\centering
\subfigure[Ground truth \qquad]{
\includegraphics[width=0.3\textwidth]{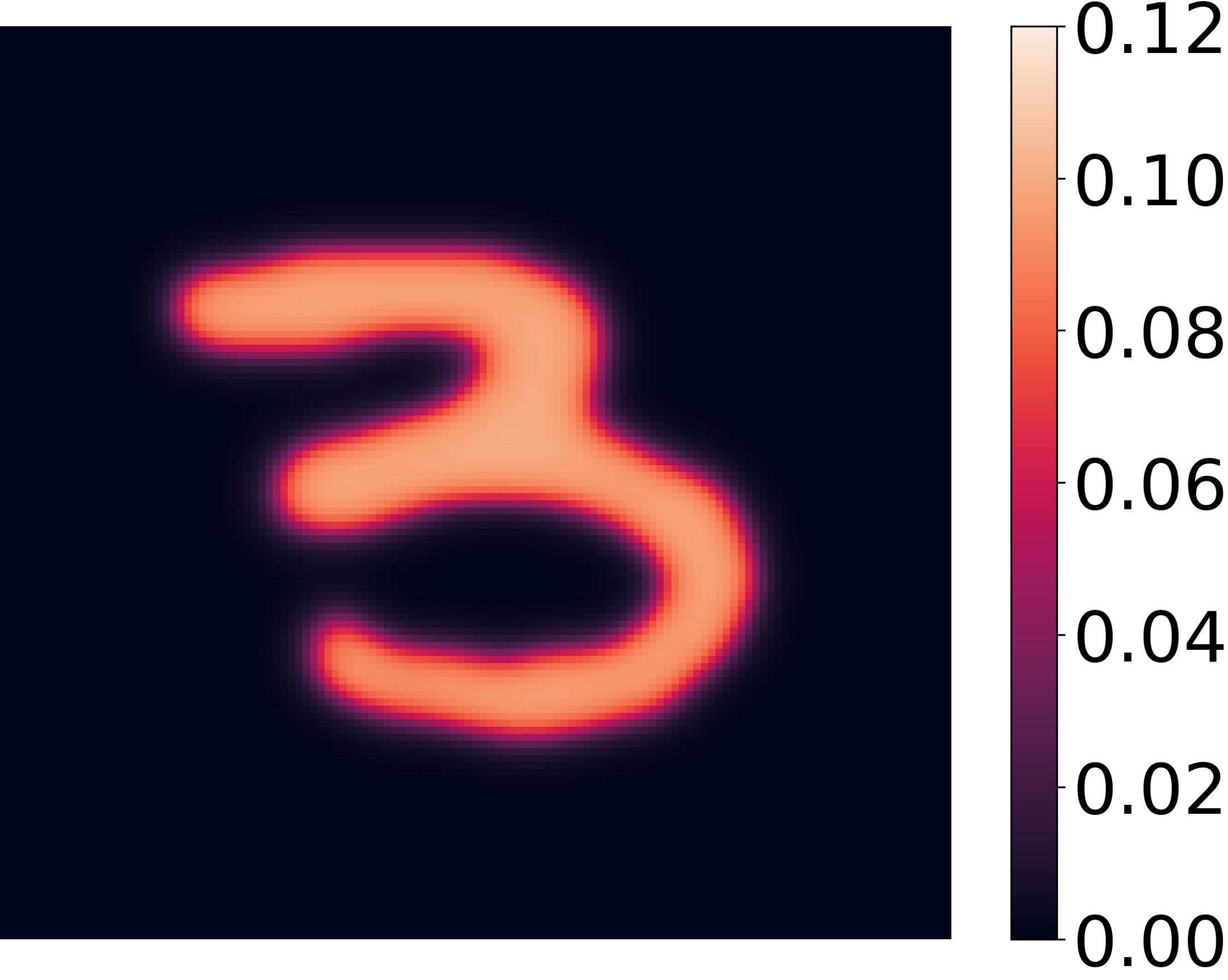}}
\hspace{0.04in}
\centering
\subfigure[FDM, 9.67\%, 25.1s \qquad]{
\includegraphics[width=0.3\textwidth]{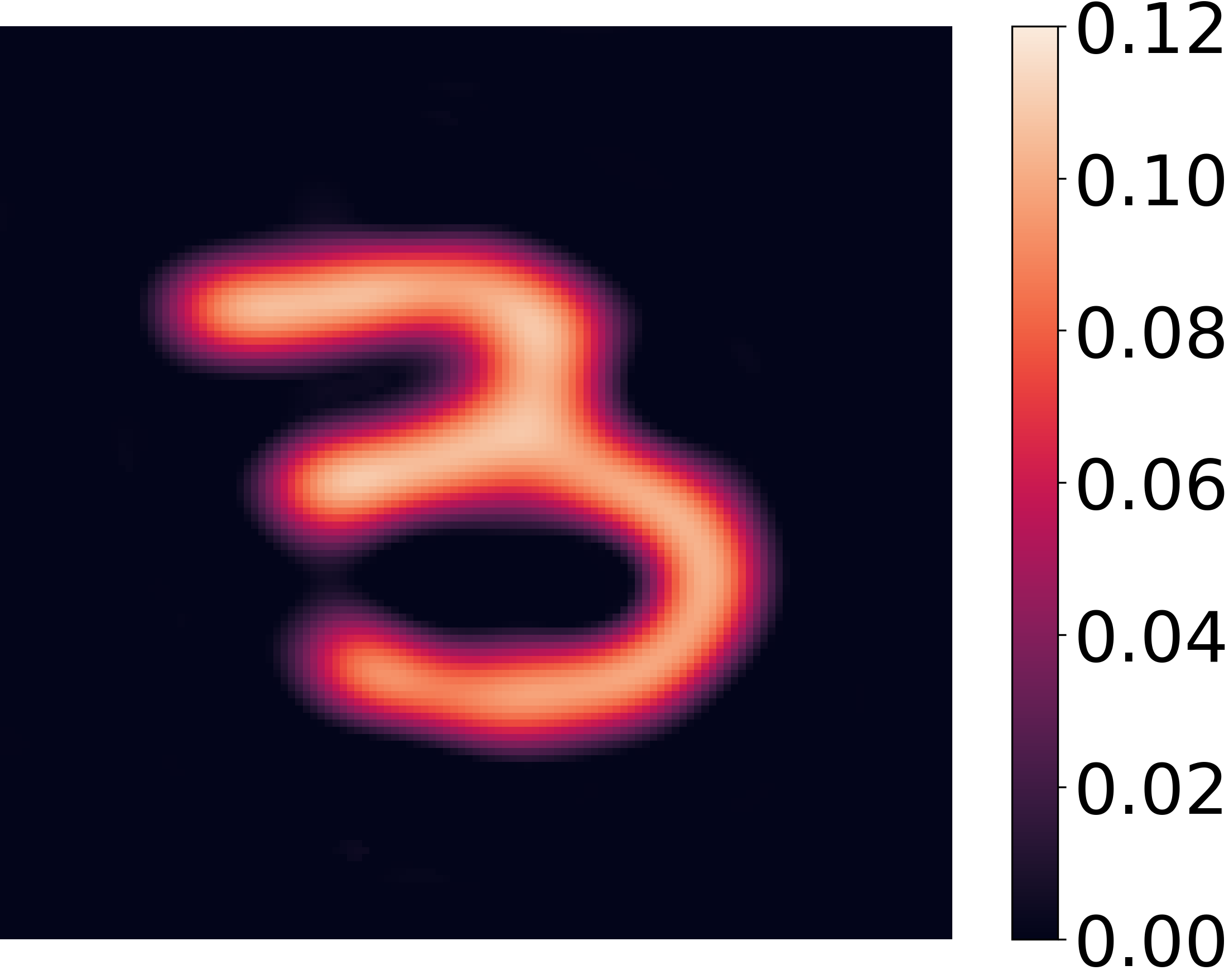}}
\hspace{0.04in}
\centering
\subfigure[NS-UNO, 11.48\%, 1.21s \qquad]{
\includegraphics[width=0.3\textwidth]{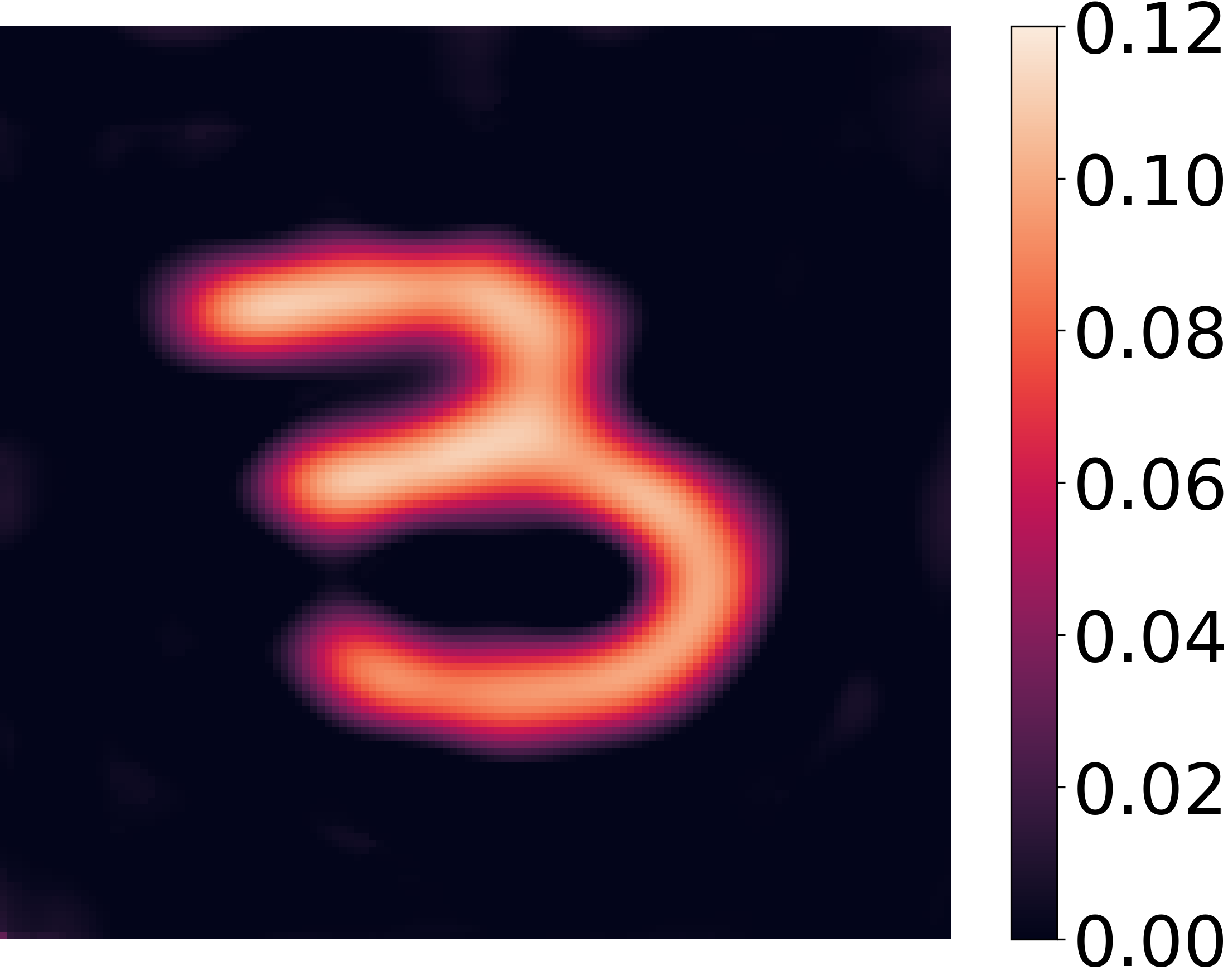}}

\subfigure[Ground truth \quad]{
\includegraphics[width=0.3\textwidth]{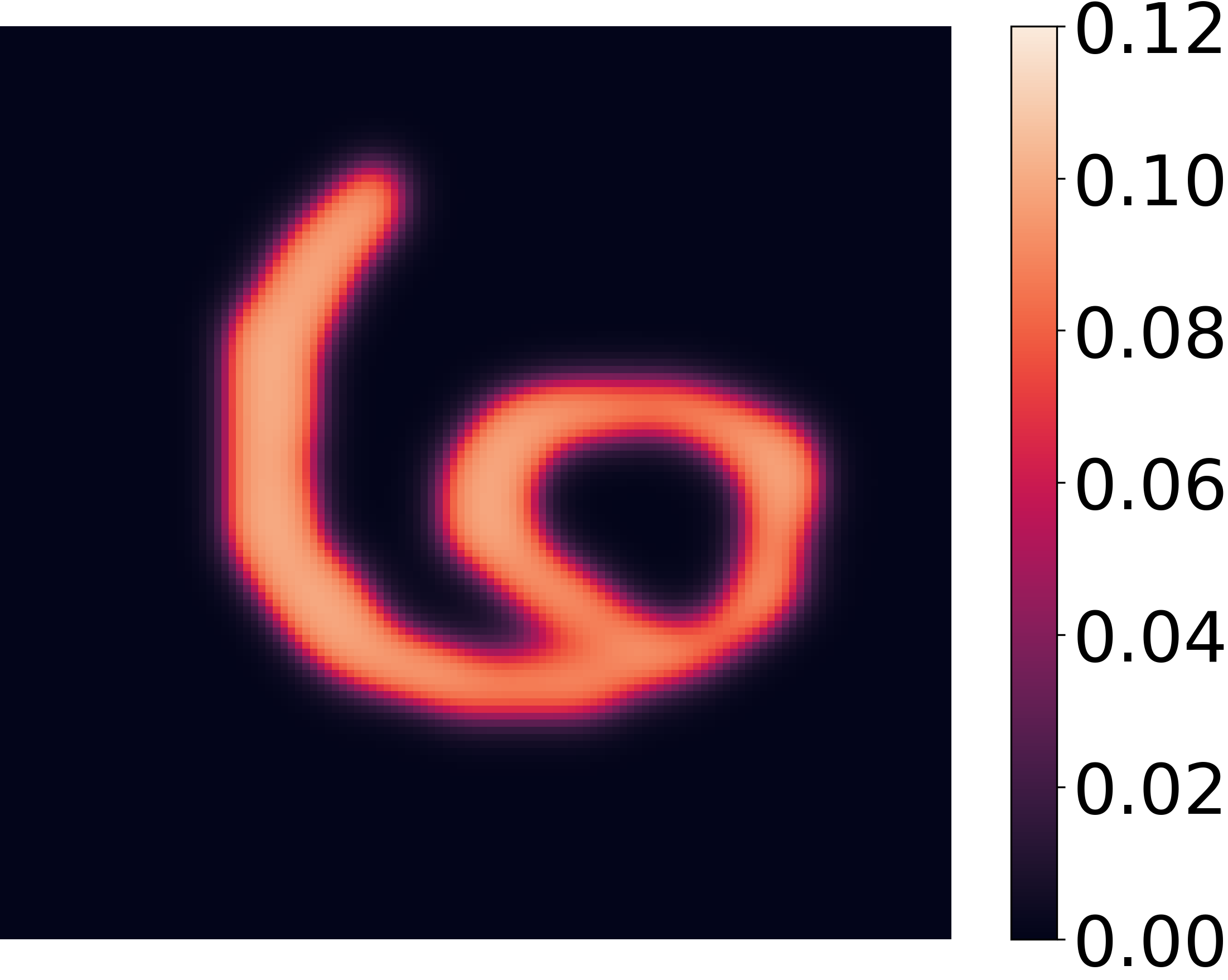}}
\hspace{0.04in}
\centering
\subfigure[FDM, 10.41\%, 38.4s\qquad]{
\includegraphics[width=0.3\textwidth]{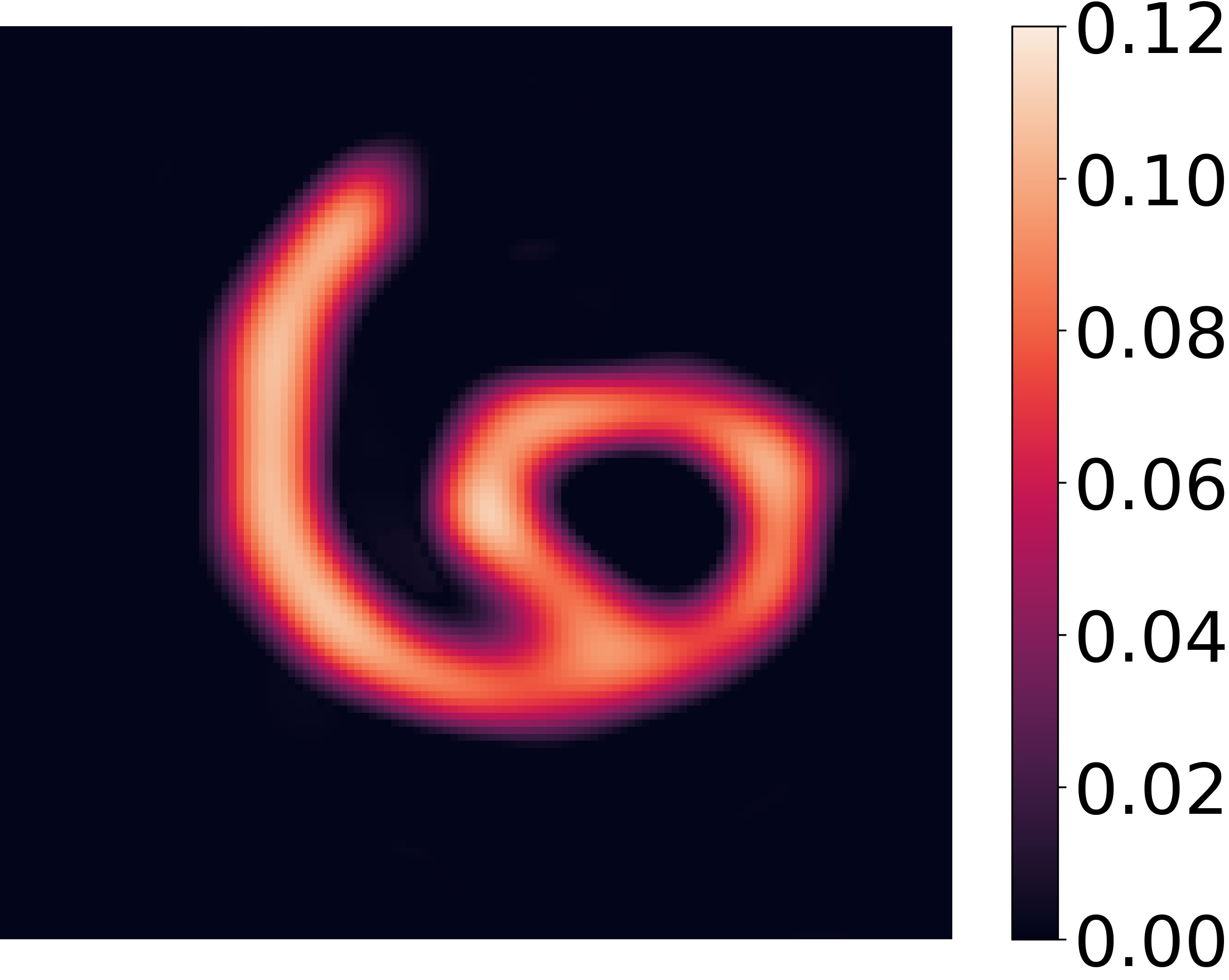}}
\hspace{0.04in}
\centering
\subfigure[NS-UNO, 12.43\%, 1.31s\qquad]{
\includegraphics[width=0.3\textwidth]{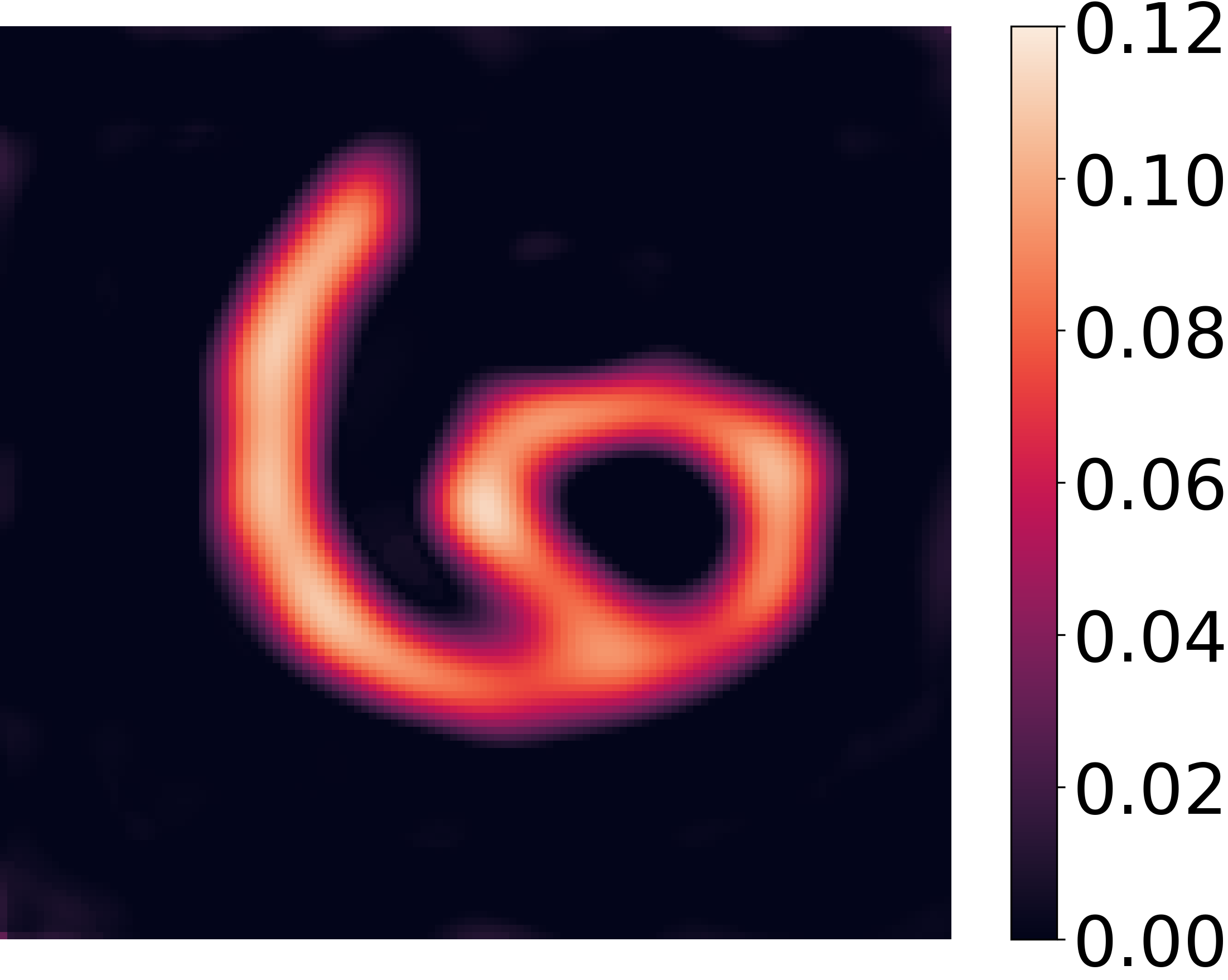}}

\subfigure[Ground truth \qquad]{
\includegraphics[width=0.3\textwidth]{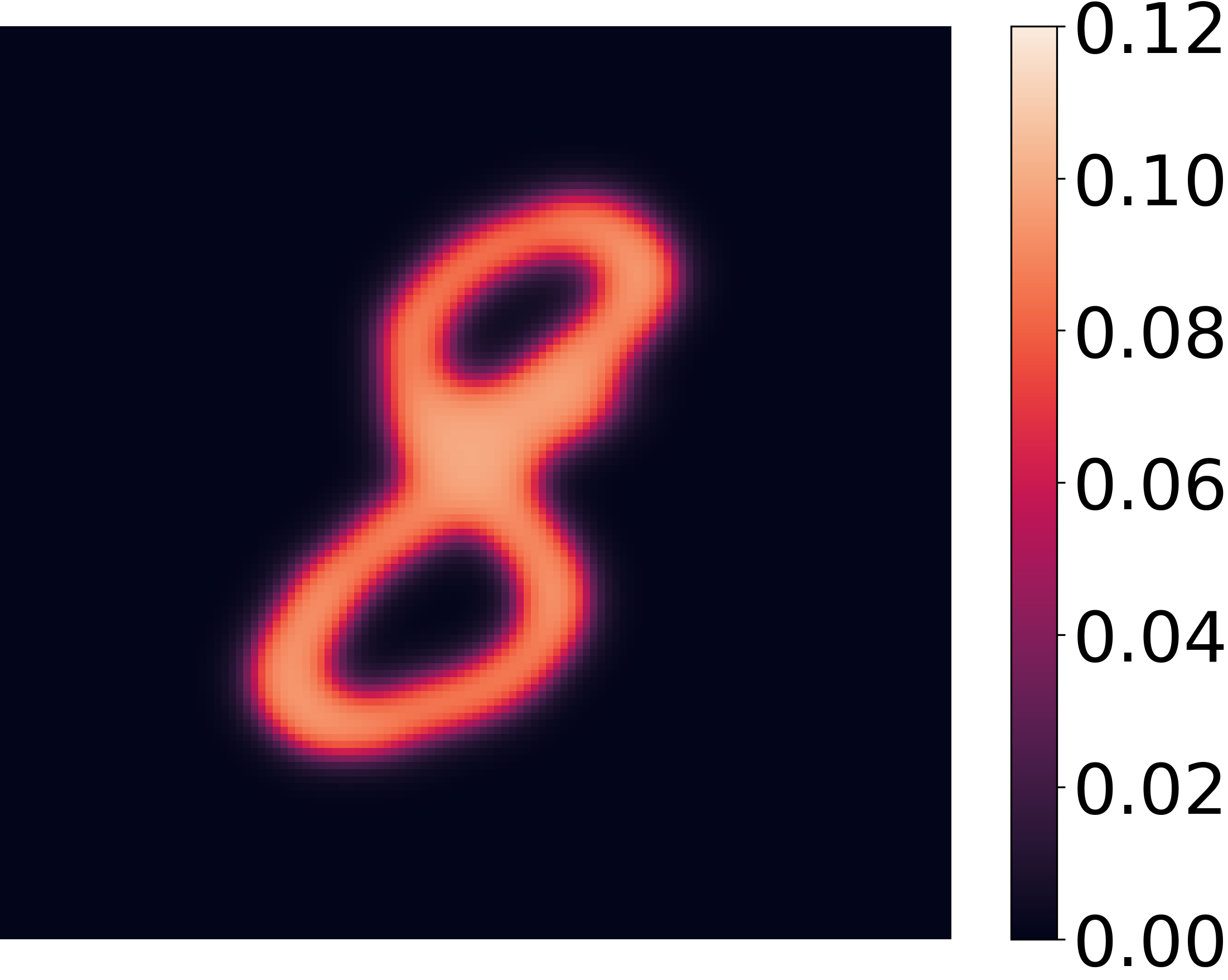}}
\hspace{0.04in}
\centering
\subfigure[FDM, 11.13\%, 51.83s\qquad]{
\includegraphics[width=0.3\textwidth]{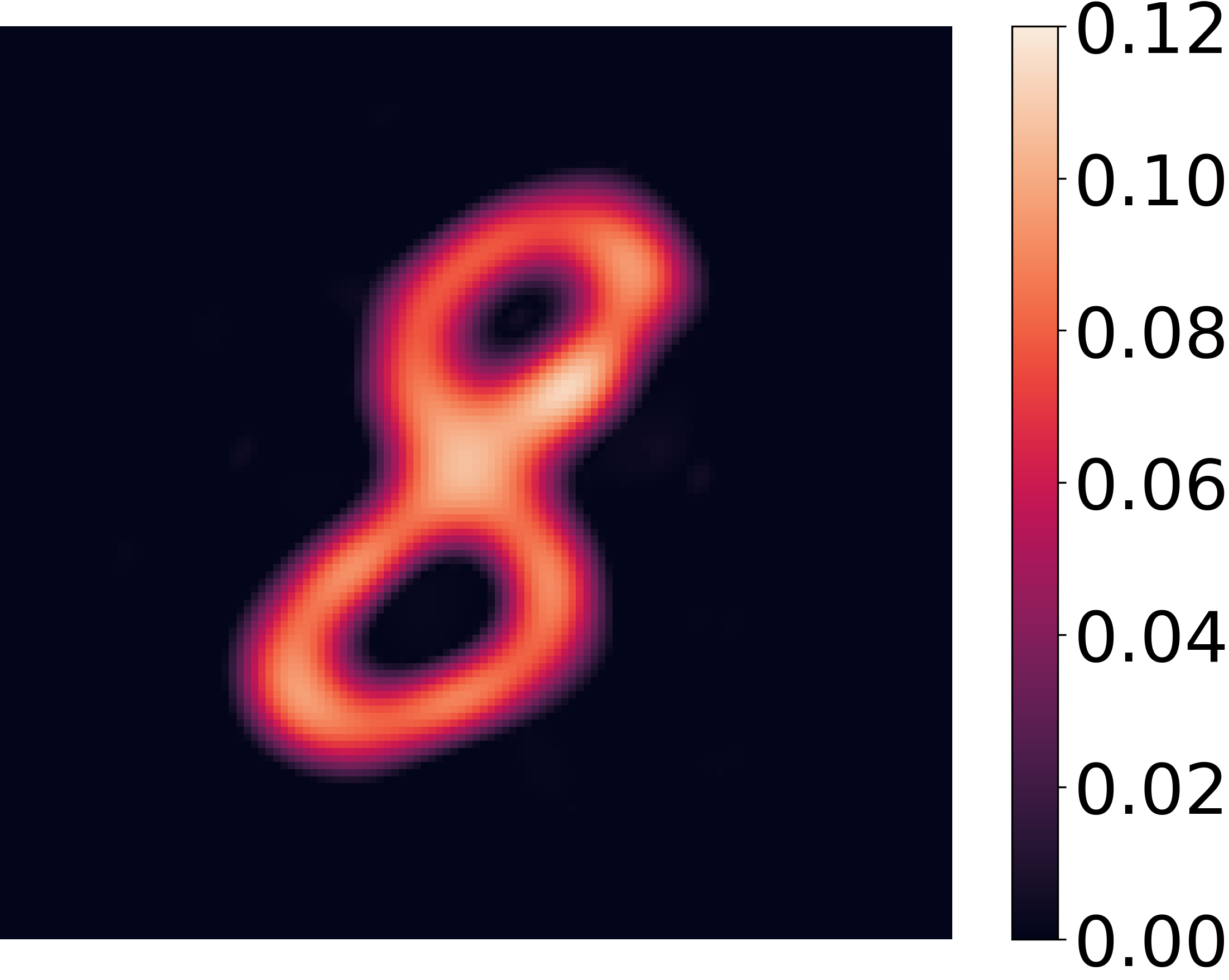}}
\hspace{0.04in}
\centering
\subfigure[NS-UNO, 13.38\%, 1.11s \qquad]{
\includegraphics[width=0.3\textwidth]{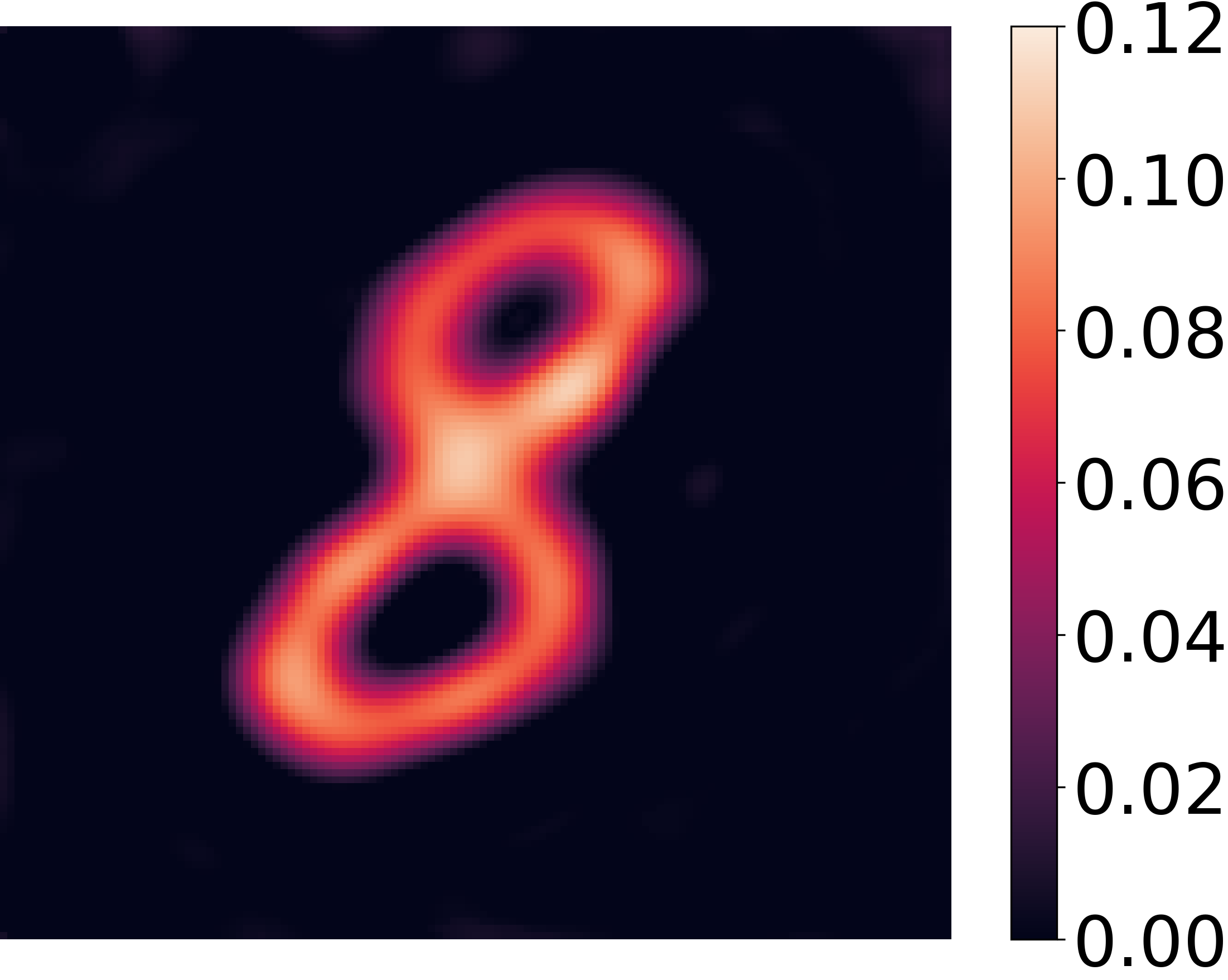}}
\caption{Relative $L^2$-error and reconstruction time of scatterer $q$ using FDM and NS-UNO as the forward solver}
\label{inverse}
\end{figure}

The reconstruction results obtained with the Finite Difference Method (FDM) and NS-UNO as the forward solver are depicted in Fig. \ref{inverse}. The visual comparison demonstrates that NS-UNO produces results on par with the finite difference method, with only a slight 2\% increase in relative $L^2$-error attributed to the inherent error of neural networks in solving forward problems. However, the key advantage of NS-UNO lies in its remarkable speed improvement. The neural network's capability to simultaneously solve all the forward problems essential for gradient computation results in NS-UNO being more than 20 times faster than the FDM with MUMPS. As a consequence, the adoption of NS-UNO as a surrogate model for the forward problem represents a significant enhancement in efficiency without causing substantial damage to accuracy.

\begin{figure}
\centering
\subfigure[Ground truth \qquad]{
\includegraphics[width=0.3\textwidth]{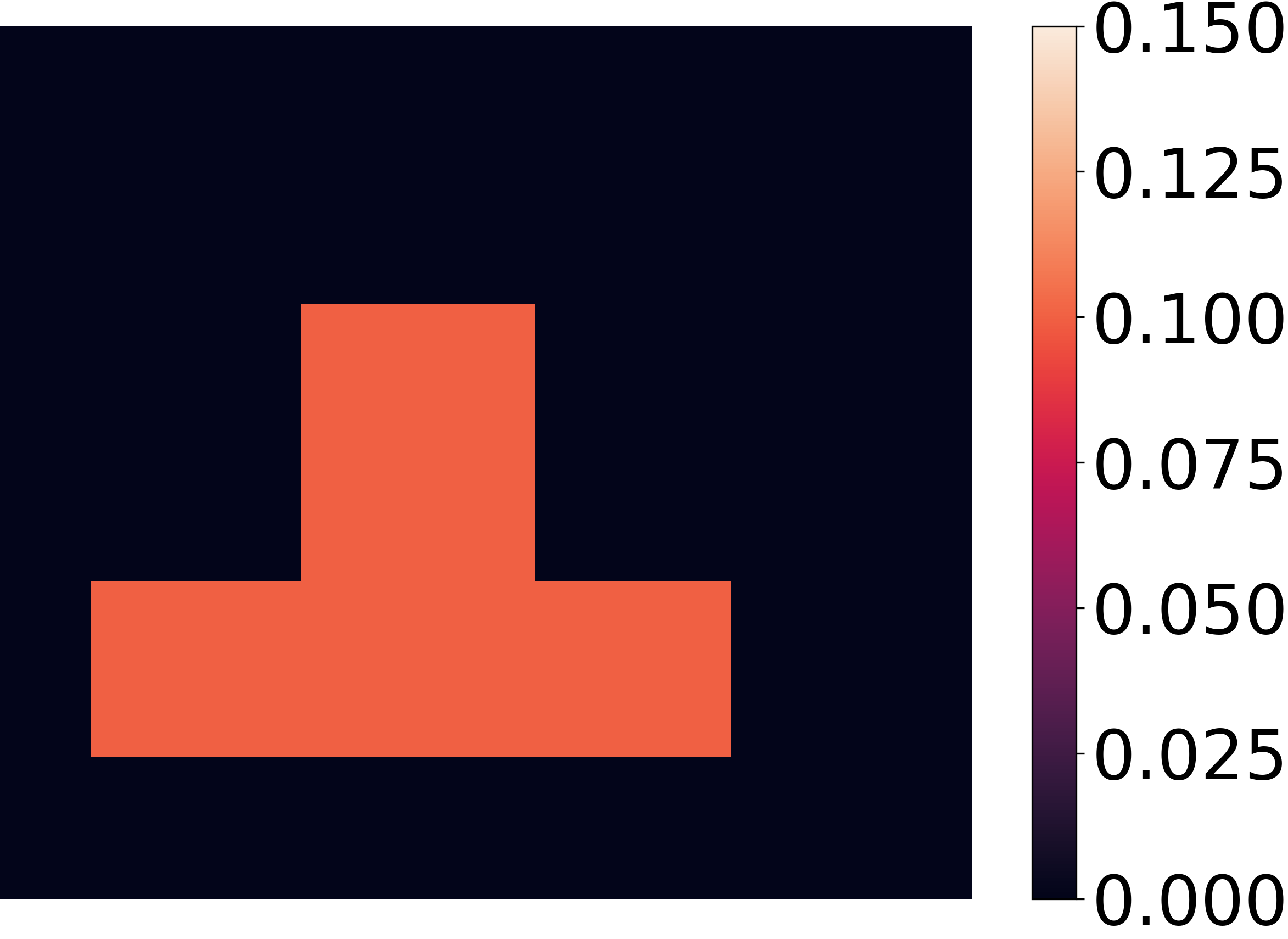}}
\hspace{0.04in}
\centering
\subfigure[FDM, 22.68\%, 29.81s \qquad]{
\includegraphics[width=0.3\textwidth]{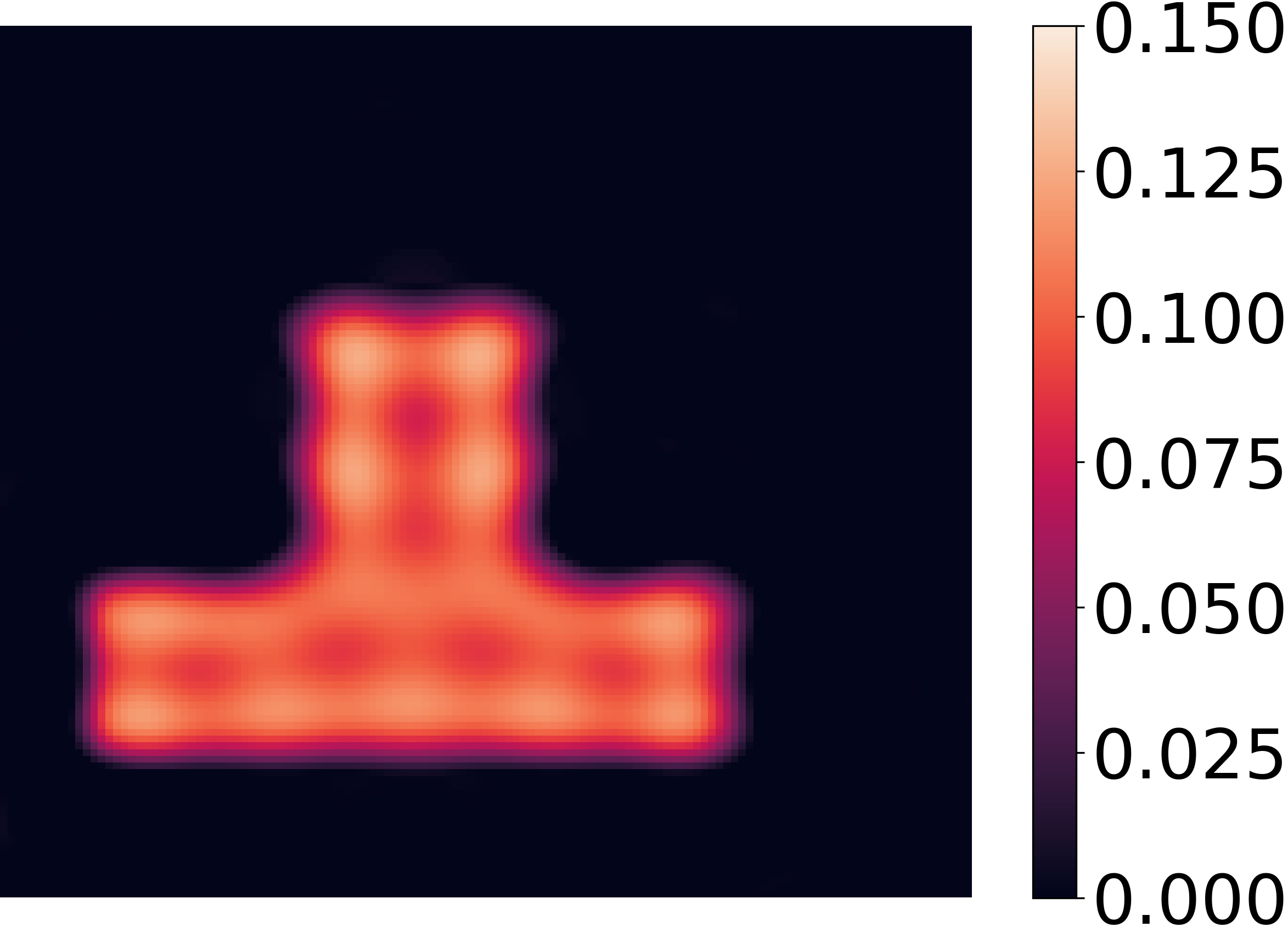}}
\hspace{0.04in}
\centering
\subfigure[NS-UNO, 25.25\%, 2.76s \qquad]{
\includegraphics[width=0.3\textwidth]{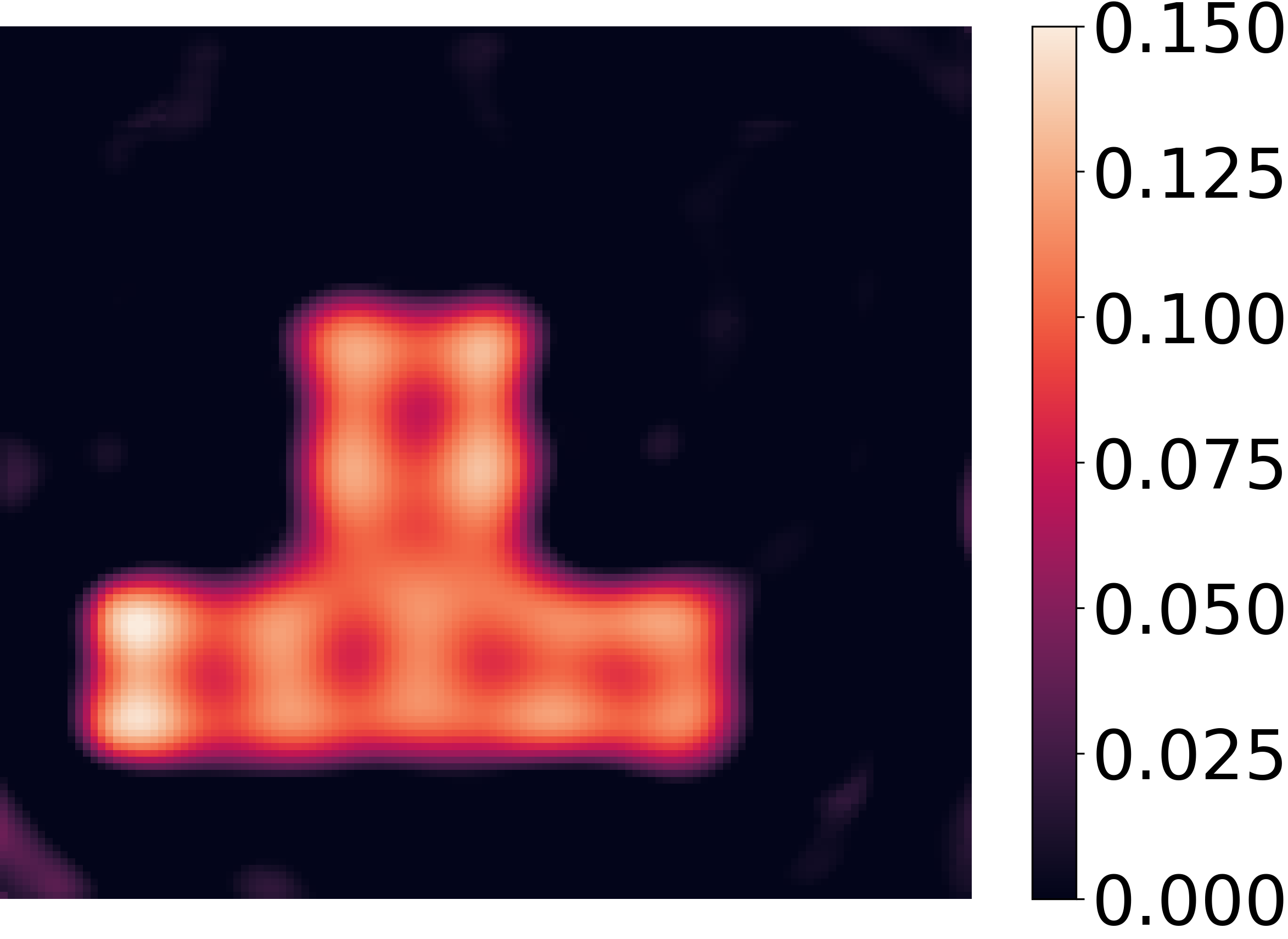}}

\subfigure[Ground truth \qquad]{
\includegraphics[width=0.3\textwidth]{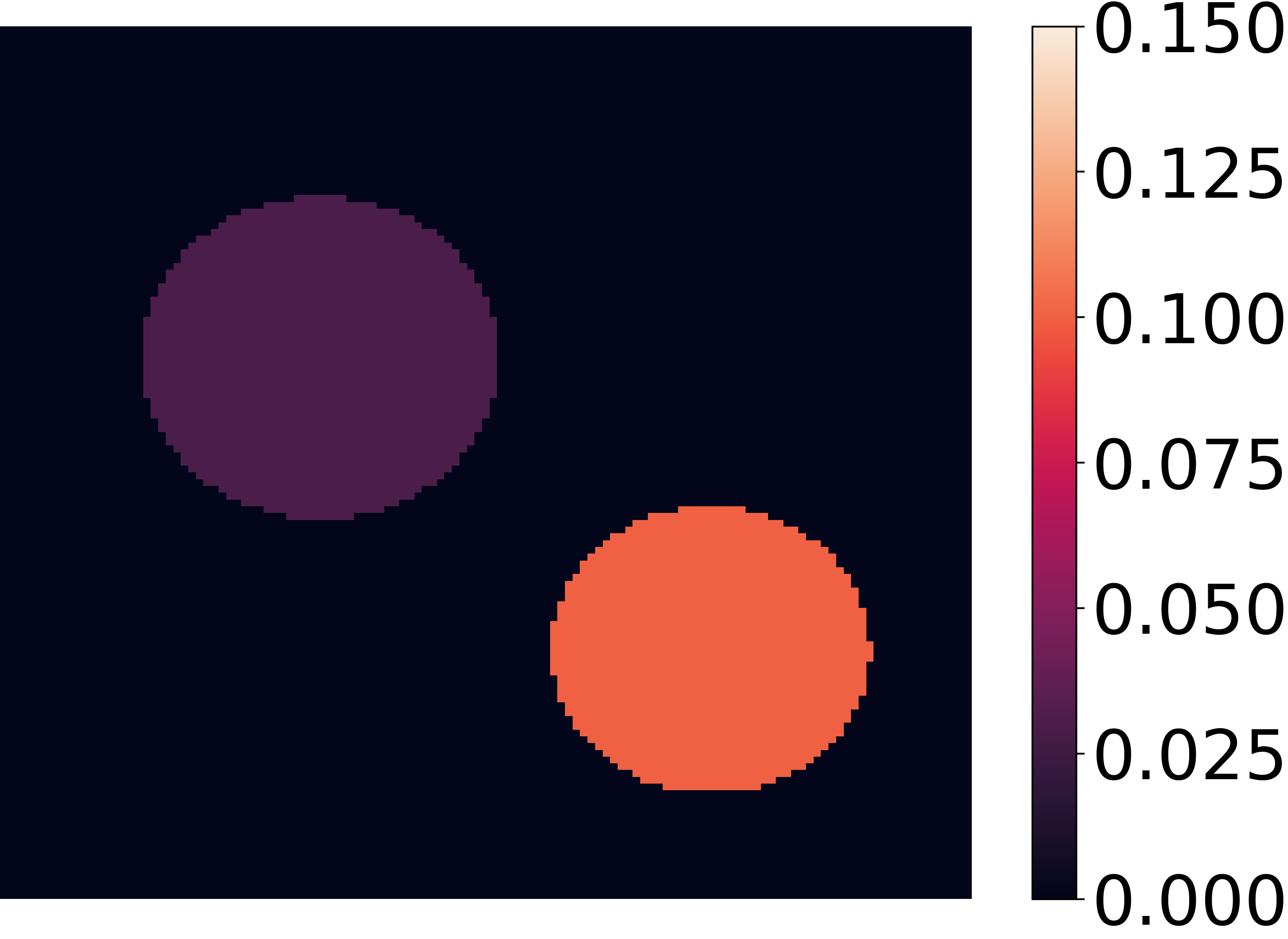}}
\hspace{0.04in}
\centering
\subfigure[FDM, 20.61\%, 41.56s \qquad]{
\includegraphics[width=0.3\textwidth]{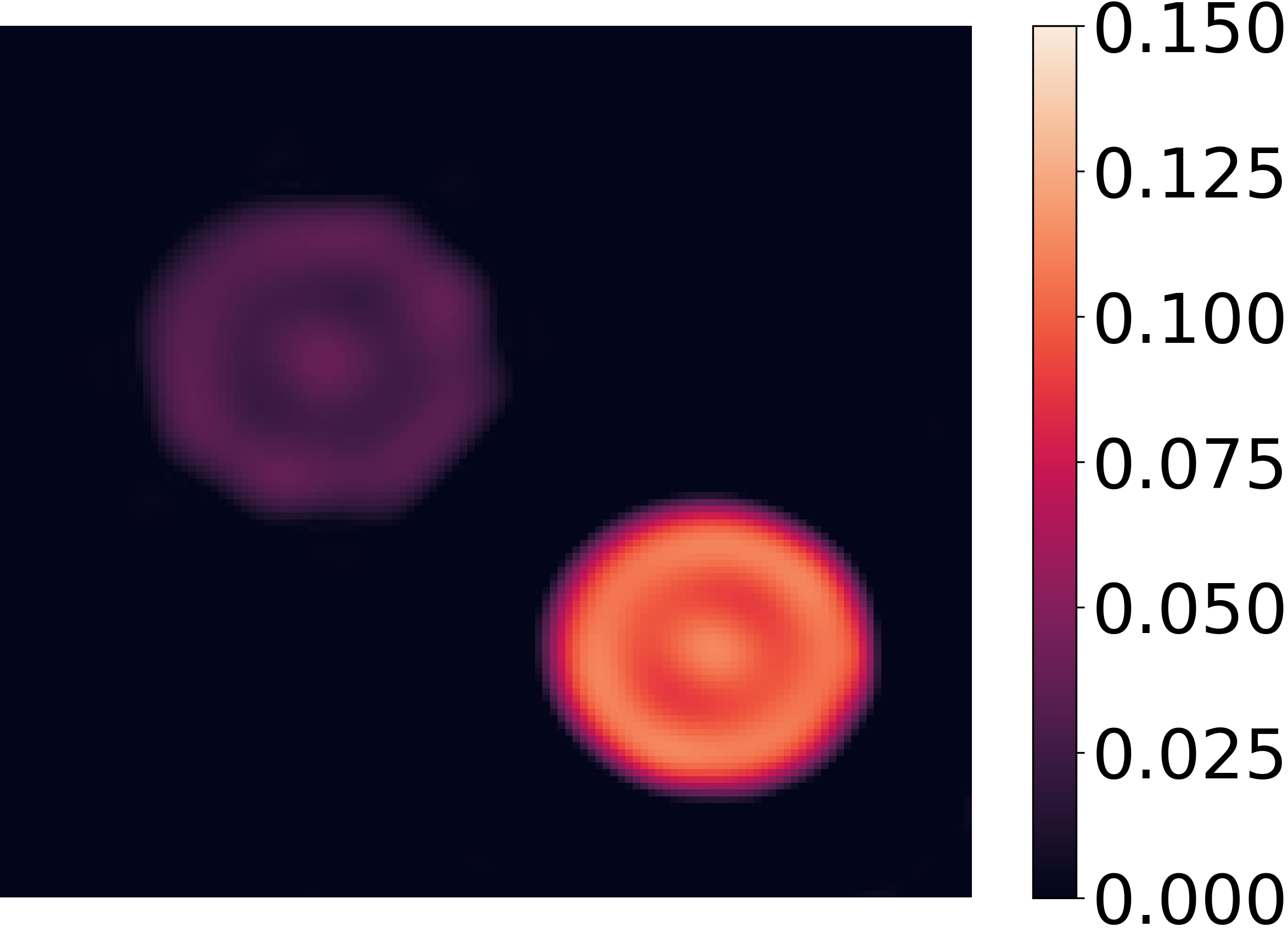}}
\hspace{0.04in}
\centering
\subfigure[NS-UNO, 26.40\%, 1.28s \qquad]{
\includegraphics[width=0.3\textwidth]{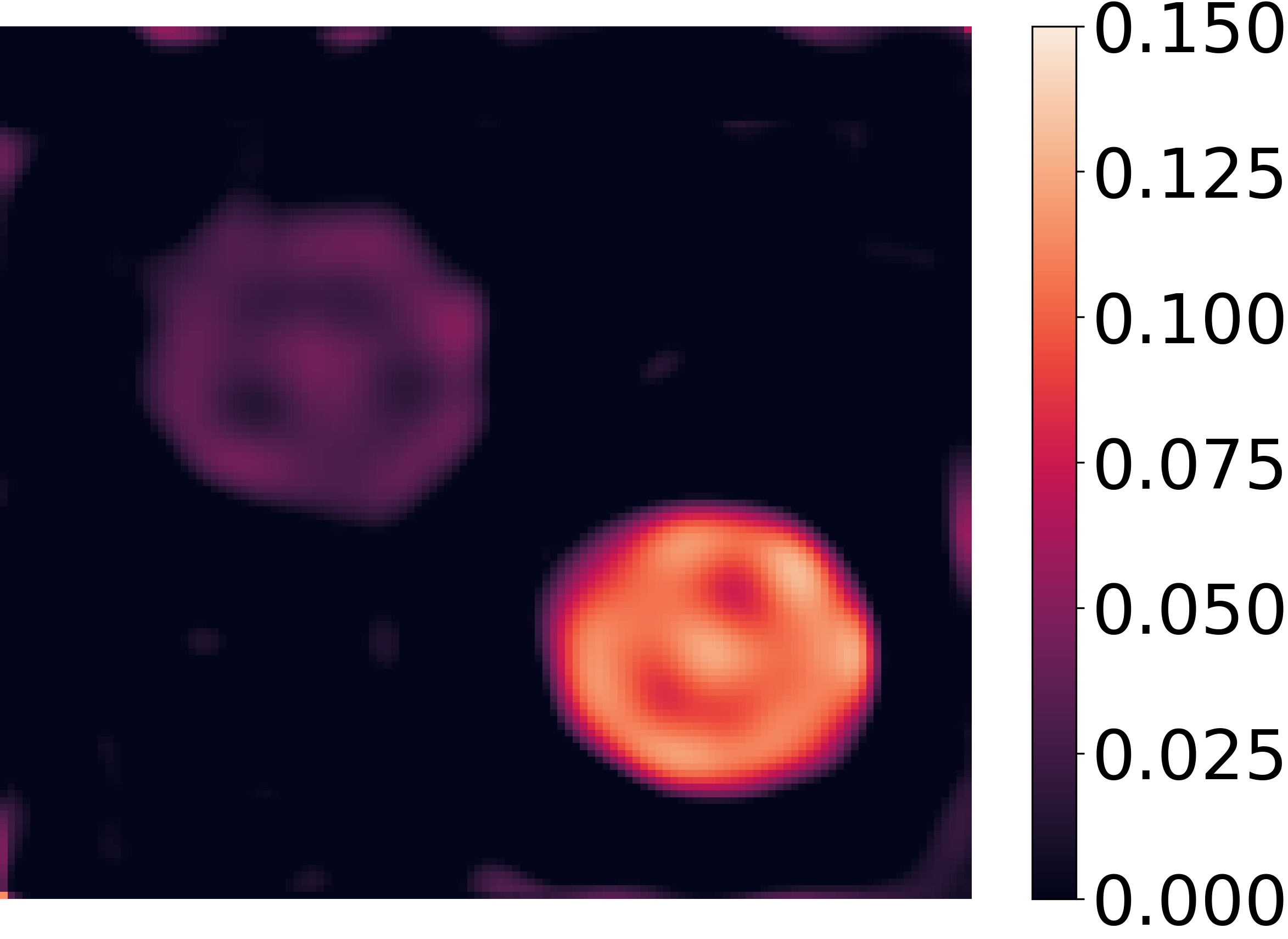}}

\caption{Relative $L^2$-error and reconstruction time of scatterer $q$ from T-shaped and random circle datasets using MUMPS and NS-UNO trained with MNIST dataset as the forward solver}
\label{inverse2}
\end{figure}

To further show the generalization ability of the proposed NS-UNO, we directly use the neural network trained with MNIST dataset to recover $q$ from T-shaped and random circle datasets, as is shown in Fig. \ref{inverse2}. It can be seen that the shape of the scatterer is accurately reconstructed by NS-UNO, and the relative error is also comparable to the finite difference method. 

\section{Conclusion}\label{con}

In this paper, we propose a novel Neumann Series Neural Operator (NSNO) to learn the solution operator of Helmholtz equation from inhomogeneity coefficients and source term to solutions. Specifically, by utilizing the Neumann series representation of the solution to the Helmholtz equation, the solution operator from both inhomogeneity coefficients and source terms to solutions can be decomposed into the solution operator from only source terms to solutions. A novel network architecture integrating U-Net structure has been designed to capture the multi-scale features in the solution to the Helmholtz equation. 

Extensive experiments have been conducted to test the performance of NSNO. The results show that NSNO has superior accuracy compared to the state-of-the-art FNO especially in the large wavenumber case, and has lower computational cost and less data requirement. The application of NSNO in inverse scattering problem have also been discussed. We have shown that the proposed NSNO is able to serve as the surrogate model and give competitive results to traditional finite difference forward solver while the computational cost is reduced significantly.



The major limitation of our work is that the restriction on the maximum value of the inhomogeneity coefficients to guarantee the convergence of Neumann series hinders the application of NSNO in high-contrast medium. In addition, the proposed Neumann series is based on the specific form of Helmholtz equation. It still remains a challenge to extend NSNO to other partial differential equations. For future work, we will try to address these issues by exploring more efficient network architecture to enlarge the application scope of NSNO, and developing alternative iterative schemes to deal with the solution operator to other partial differential equations.



\appendix
\section{Adjoint state method}
Without loss of generality, we consider the model with one incident wave denoted as $u_0$. The optimization problem is then given as
\begin{align}
\min_{q} \quad &J(q)=\frac{1}{2}\Vert Tu(q)-d\Vert_2^2,\\
\text{s.t.} \quad &\Delta u+k^2(1+q)u=-k^2qu_0, \quad \text{in} \quad\it\Omega, \label{A1}\\
&\frac{\partial 
u}{\partial n}=iku, \quad \text{on} \quad\partial\it\Omega,\label{A2}
\end{align}
where $T$ denotes the trace operator restricting the wave field on the boundary. In order to derive $\frac{\partial J}{\partial q}$, we apply the Lagrange multiplier method. Define the Lagrangian as
\begin{equation}
\mathcal{L}(u,\lambda,q)=J(q)-(\lambda, \Delta u+k^2(1+q)u+k^2qu_0),
\end{equation}
where $(f, g)=\operatorname{Re}\int_{\it\Omega}f\bar{g}$ denotes the real part of the inner product in $L^2(\it\Omega)$. After two integrations by part we obtain
\begin{equation}
\begin{split}
\mathcal{L}(u,\lambda,q)&=J(q)-\left\langle\lambda, \frac{\partial u}{\partial n}\right\rangle+\left\langle\frac{\partial \lambda}{\partial n}, u\right\rangle-(\Delta\lambda, u)-(\lambda, k^2(1+q)u)-(\lambda, k^2qu_0)\\
&=J(q)+\left\langle\frac{\partial \lambda}{\partial n}+ik\lambda, u\right\rangle-(\Delta\lambda+k^2(1+q)\lambda, u)-(\lambda, k^2qu_0),
\end{split}
\end{equation}
where $\langle f, g\rangle =\operatorname{Re}\int_{\it\Omega}f\bar{g}$ denotes the real part of the inner product in $L^2(\partial\it\Omega)$. Therefore,
\begin{equation}
\frac{\partial\mathcal{L}}{\partial q}(q)=\left(T^{\ast}(Tu-d), \frac{\partial u}{\partial q}(q)\right)+\left\langle\frac{\partial \lambda}{\partial n}+ik\lambda, \frac{\partial u}{\partial q}(q)\right\rangle-\left(\Delta\lambda+k^2(1+q)\lambda, \frac{\partial u}{\partial q}(q)\right)-(k^2\lambda, u)-(\lambda, k^2u_0).
\end{equation}
To eliminate $\frac{\partial u}{\partial q}(q)$, we choose $\lambda$ to satisfy
\begin{equation}
\begin{split}
\Delta\lambda+k^2(1+q)\lambda&=T^{\ast}(Tu-d), \quad \text{in} \quad \it\Omega,\\
\frac{\partial \lambda}{\partial n}+ik\lambda &= 0, \quad \text{on} \quad \partial\it\Omega,
\end{split}
\end{equation}
or equivalently we solve the conjugate of $\lambda$ satisfying
\begin{equation}\label{A3}
\begin{split}
\Delta\bar{\lambda}+k^2(1+q)\bar{\lambda}&=\bar{T^{\ast}(Tu-d)}, \quad \text{in} \quad \it\Omega,\\
\frac{\partial \bar{\lambda}}{\partial n}&=ik\bar{\lambda}, \quad \text{on} \quad \partial\it\Omega.
\end{split}
\end{equation}
The gradient of $J$ with respect to $q$ can then be obtained by
\begin{equation}\label{A4}
\frac{\partial J(q)}{\partial q}=\frac{\partial \mathcal{L}(q)}{\partial q}=-k^2(\bar{\lambda}, u+u_0).
\end{equation}
The process to compute the gradient of $J$ with respect to $q$ is summarized as follows, where the forward solver is called twice:
\begin{itemize}
\item Solve $u$ satisfying \eqref{A1} and \eqref{A2}.
\item Solve $\bar{\lambda}$ satisfying \eqref{A3}.
\item Compute the gradient by \eqref{A4}.
\end{itemize}

\bibliographystyle{unsrt} 
\bibliography{refs}

\end{document}